\documentclass[12pt]{amsart}
\usepackage{amsmath,amssymb}

\newcommand{\A}[3]{{}^{#1} \hskip -2.5pt a^{#2}_{#3}}
\newcommand{\W}[3]{{}^{#1} \hskip -2pt W^{#2}_{#3}}
\newcommand{\la}[3]{{}^{#1} \hskip -2pt {\lambda}^{#2}_{#3}}

\newcommand\rn{{\mathbb R}^n}

\newcommand\Ha{\mathbb H}

\newcommand\inv{^{-1}}

\newcommand\fv{\mathfrak v}
\newcommand\fh{\mathfrak h}
\newcommand\fH{\mathcal{H}}
\newcommand\fk{\mathfrak k}
\newcommand\fl{\mathfrak l}

\newcommand\fu{\mathfrak u}

\newcommand\fz{\mathfrak z}

\newcommand\fm{\mathfrak m}
\newcommand\ff{\mathfrak f}
\newcommand\ft{\mathfrak t}

\newcommand\br{\bar \rho}

\newcommand\BG{P{\backslash}G}
\newcommand\La{\mathbb L}
\newcommand\Ga{\mathbb G}
\newcommand\Pa{\mathbb P}
\newcommand\Ja{\mathbb J}
\newcommand\Aa{\mathbb A}
\newcommand\Fa{\mathbb F}
\newcommand\Ma{\mathbb M}
\newcommand\Ua{\mathbb U}
\newcommand\Ta{\mathbb T}
\newcommand\Sa{\mathbb S}

\newcommand\Ra{\mathbb R}
\newcommand\Qa{\mathbb Q}
\newcommand\Za{\mathbb Z}
\newcommand\Da{\mathbb D}

\DeclareMathOperator{\Diff}{Diff}
\DeclareMathOperator{\Rat}{Rat} \DeclareMathOperator{\Aff}{Aff}

\DeclareMathOperator{\Out}{Out} 
 \DeclareMathOperator{\Hom}{Hom}
\DeclareMathOperator{\Aut}{Aut} \DeclareMathOperator{\Ad}{Ad}
\DeclareMathOperator{\supp}{supp} \DeclareMathOperator{\pol}{pol}
 \DeclareMathOperator{\Isom}{Isom}
\DeclareMathOperator{\Zar}{Zar}
\DeclareMathOperator{\Exp}{Exp}

\newtheorem{theorem}{Theorem}[section]
\newtheorem{proposition}[theorem]{Proposition}
\newtheorem{lemma}[theorem]{Lemma}
\newtheorem{corollary}[theorem]{Corollary}
\newtheorem{defn}[theorem]{Definition}
\newtheorem{example}[theorem]{Example}

\begin{document}

\title[Local rigidity for cocycles]{Local rigidity for cocycles}
\author[D. Fisher and G.~A.~Margulis]{David Fisher and G.~A.~Margulis}
\thanks{First author partially supported by NSF grants DMS-9902411 and DMS-0226121 and a PSC-CUNY grant.
Second author partially supported by NSF grant DMS-9800607. The
authors would also like to thank the FIM at ETHZ for hospitality
and support.}

\begin{abstract}
In this paper we study perturbations of constant cocycles for
actions of higher rank semi-simple algebraic groups and their
lattices. Roughly speaking, for ergodic actions, Zimmer's cocycle
superrigidity theorems implies that the perturbed cocycle is
measurably conjugate to a constant cocycle modulo a compact valued
cocycle.  The main point of this article is to see that a cocycle
which is a continuous perturbation of a constant cocycle is
actually continuously conjugate back to the original constant
cocycle modulo a cocycle that is continuous and ``small".

We give some applications to perturbations of standard actions of
higher rank semisimple Lie groups and their lattices.  Some of the
results proven here are used in our proof of local rigidity for
affine and quasi-affine actions of these groups.

We also improve and extend the statements and proofs of Zimmer's
cocycle superrigidity.
\end{abstract}

\maketitle


\section{{\bf Introduction}}
\label{sec:intro}

Let $G$ be a connected semisimple Lie group with no compact
factors and all simple factors of real rank at least two. Further
assume $G$ is simply connected as a Lie group or simply connected
as an algebraic group. The latter means that $G=\Ga(\mathbb R)$
where $\Ga$ is a simply connected semisimple $\mathbb R$-algebraic
group.  Let $\Gamma<G$ be a lattice and $L$ be the $k$ points of
an algebraic $k$-group where $k$ is a local field of
characteristic zero. Roughly speaking, Zimmer's cocycle
superrigidity theorems imply that any cocycle into $L$ over an
ergodic action of $G$ or $\Gamma$ is measurably conjugate to a
constant cocycle, modulo some compact noise. (See below for a
precise formulation.) This theorem has many consequences for the
dynamics of smooth actions of these groups. Even stronger results
would follow if one could produce a continuous or smooth
conjugacy. The main purpose of this paper is to prove that a
perturbation of a constant cocycle is conjugate back to the
constant cocycle via a small (and often continuous) conjugacy,
modulo ``small" noise.  We also prove stronger and more general
versions of the cocycle superrigidity theorems than had previously
been known.  In particular, we do not need to pass to a finite
ergodic extension of the action and we obtain more general
statements when $k$ is non-Archimedean.

Throughout we work with a more general group $G$.  We let $I$ be a
finite index set and for each $i{\in}I$, we let $k_i$ be a local
field of characteristic zero and ${\mathbb G}_i$ be a connected
simply connected semisimple algebraic $k_i$-group. We first define
groups $G_i$, and then let $G=\prod_{i{\in}I}G_i$.  If $k_i$ is
non-Archimedean,
 $G_i={\mathbb G}_i(k_i)$ the
$k_i$-points of ${\mathbb G}_i$. If $k_i$ is Archimedean, then
$G_i$ is either ${\mathbb G}_i(k_i)$ or its topological universal
cover. (This makes sense, since when $\Ga_i$ is simply connected
and $k_i$ is Archimedean, ${\mathbb G}_i(k_i)$
 is topologically connected.) Throughout the introduction, we
 assume that the $k_i$-rank of any simple factor of any ${\mathbb G}_i$ is at least two.

We first state a version of our main result for $G$ actions and
cocycles.

\begin{theorem}
\label{theorem:localcocyclerigidityG} Let $G$ be as above, $L=\La(k)$
 where
$\La$ is an algebraic $k$-group and $k$ is a local field of
characteristic zero and let  $\pi_0:G{\rightarrow}L$ be a
continuous homomorphism. Let
$(S,\mu)$ be a standard probability measure space,  $\rho$ a
measure preserving action of $G$ on $S$, and
$\alpha_{\pi_0}:G{\times}S{\rightarrow}L$ be the constant cocycle
over the action $\rho$ given by $\alpha_{\pi_0}(g,x)=\pi_0(g)$.
Assume $\alpha:G{\times}S{\rightarrow}L$ is a Borel cocycle over
the action $\rho$ such that $\alpha$ is $L^{\infty}$ close to
$\alpha_{\pi_0}$. Then there exists a measurable map
$\phi:S{\rightarrow}L$, and a cocycle $z:G{\times}S{\rightarrow}Z$
where $Z=Z_L(\pi_0(G))$, the centralizer in $L$ of $\pi_0(G)$,
such that
\begin{enumerate}

\item we have $\alpha(g,x)=\phi(gx){\inv}\pi_0(g)z(g,x)\phi(x)$;

\item $\phi:S{\rightarrow}L$ is small in $L^{\infty}$;

\item the cocycle $z$ is $L^{\infty}$ close to the trivial cocycle

\item the cocycle $z$ is measurably conjugate to a cocycle taking
values in a compact subgroup $C$ of $Z$ where $C$ is contained in
a small neighborhood of the identity.
\end{enumerate}

\noindent Furthermore if $S$ is a locally compact topological
space, $\mu$ is a Borel measure on $S$ with $\supp(\mu)=S$ and
$\alpha$ and $\rho$ are continuous then both $\phi$ and $z$ can be
chosen to be continuous.
\end{theorem}

\noindent
{\bf Remark:}  If $k$ is Archimedean, point $(4)$ implies
that $z$ is measurably conjugate to the trivial cocycle.

Before stating the analogous theorem for $\Gamma$ actions and
cocycles, we need to recall a consequence of the superrigidity
theorems \cite{M1,Ma,M2}. We will use the notation introduced here
in the statements below. If $G$ is as above and $\Gamma<G$ is a
lattice, we call a homomorphism $\pi:\Gamma{\rightarrow}L$ {\em
superrigid} if it almost extends to a homomorphism of $G$. This
means that there is a continuous homomorphism
$\pi^E:G{\rightarrow}L$ and a homomorphism
$\pi^K:\Gamma{\rightarrow}L$ with bounded image such that
$\pi(\gamma)=\pi^E(\gamma)\pi^K(\gamma)$ and $\pi^E(\Gamma)$
commutes with $\pi^K(\Gamma)$. The superrigidity theorems imply
that any continuous homomorphism of $\Gamma$ into an algebraic
group is superrigid. This can be deduced easily from Lemma VII.5.1
and Theorems VII.5.15 and VII.6.16 of \cite{M2}.


\begin{theorem}
\label{theorem:localcocyclerigidityGamma} Let $\Gamma$ be as
above,
$L=\La(k)$ be as in Theorem \ref{theorem:localcocyclerigidityG} and
$\pi_0:\Gamma{\rightarrow}L$ be a continuous homomorphism. Let $(S,\mu)$ be a standard probability measure space,
$\rho$ be a measure preserving action of $\Gamma$ on $S$,
 and let
$\alpha_{\pi_0}:\Gamma{\times}S{\rightarrow}L$ be the constant
cocycle over the action $\rho$ given by
$\alpha_{\pi_0}(\gamma,x)=\pi_0(\gamma)$. Assume
$\alpha:\Gamma{\times}S{\rightarrow}L$ is a Borel cocycle over the
action $\rho$ such that $\alpha$ is $L^{\infty}$ close to
$\alpha_{\pi_0}$. Then there exists a measurable map
$\phi:S{\rightarrow}L$, and a cocycle
$z:\Gamma{\times}S{\rightarrow}Z$ where $Z=Z_L(\pi_0^E(G))$ such
that
\begin{enumerate}

\item we have
$\alpha(\gamma,x)=\phi({\gamma}x){\inv}\pi_0^E(\gamma)z(\gamma,x)\phi(x)$;

\item $\phi:S{\rightarrow}L$ is small in $L^{\infty}$;

\item the cocycle $z$ is $L^{\infty}$ close to the constant
cocycle defined by $\pi_0^K$

\item  $z$ is measurably conjugate to a cocycle taking values in a
compact subgroup $C$ of $Z$ where $C$ is contained in a small
neighborhood of $\pi_0^K(\Gamma)$.
\end{enumerate}

\noindent Furthermore if $S$ is a locally compact topological
space, $\mu$ is a Borel measure on $S$ with $\supp(\mu)=S$ and
$\alpha$ and $\rho$ are continuous then both $\phi$ and $z$ can be
chosen to be continuous.
\end{theorem}

{\noindent}
{\bf Remark:} If $k$ is Archimedean, point $(4)$
implies that $z$ is measurably conjugate to a cocycle taking
values in the closure of $\pi_0^K(\Gamma)$.

To prove Theorems \ref{theorem:localcocyclerigidityG} and
\ref{theorem:localcocyclerigidityGamma}, we prove a very general
result about perturbations of cocycles over measure preserving
actions of groups with property T.  The result shows that any
perturbation of a cocycle taking values in a compact group also
takes values in a compact group, see Theorem
\ref{theorem:Tlocalrigidity}.

Use of (an extension and modification of) Zimmer's cocycle
superrigidity theorems is a key step in the proof of Theorem's
\ref{theorem:localcocyclerigidityG} and
\ref{theorem:localcocyclerigidityGamma}.  The cocycle
superrigidity theorems are generalizations of the second author's
superrigidity theorems.  Our strongest results require an
integrability condition on the cocycles considered.

\begin{defn}
\label{defn:quasi-integrable} Let $D$ be a locally compact group,
$(S,\mu)$ a standard probability measure space on which $D$ acts
preserving $\mu$  and $L$ be a normed topological group. We call a
cocycle $\alpha:D{\times}S{\rightarrow}L$ over the $D$ action {\em
$D$-integrable} if for any compact subset $M\subset{D}$, the
function $Q_{M,\alpha}(x)=\supp_{m{\in}M}\ln^+\|\alpha(m,x)\|$ is
in $L^1(S)$.
\end{defn}

Any continuous cocycle over a continuous action on a compact
topological space is automatically $D$-integrable.  We remark that
a cocycle over a cyclic group action is $D$-integrable if and only
if $\ln^+\|(\alpha(\pm{1},x)\|$ is in $L^1(S)$.

\begin{theorem}
\label{theorem:Gsuperrigidity} Let $G,S,\mu,L$ be as in Theorem
\ref{theorem:localcocyclerigidityG}. Assume $G$ acts ergodically
on $S$ preserving $\mu$. Let $\alpha:G{\times}S{\rightarrow}L$ be
a $G$-integrable Borel cocycle.  Then $\alpha$ is cohomologous to
a cocycle $\beta$ where $\beta(g,x)=\pi(g){c(g,x)}$. Here
${\pi:G{\rightarrow}L}$ is a continuous homomorphism and
$c:G{\times}S{\rightarrow}C$ is a cocycle taking values in a
compact group centralizing $\pi(G)$.
\end{theorem}

\begin{theorem}
\label{theorem:Gammasuperrigidity} Let $G,\Gamma,S,L$ and $\mu$ be
as Theorem \ref{theorem:localcocyclerigidityGamma}. Assume
$\Gamma$ acts ergodically on $S$ preserving $\mu$. Assume
$\alpha:\Gamma{\times}S{\rightarrow}L$ is a $\Gamma$-integrable,
Borel cocycle. Then $\alpha$ is cohomologous to a cocycle $\beta$
where $\beta(\gamma,x)=\pi(\gamma){c(\gamma,x)}$. Here
${\pi:G{\rightarrow}L}$ is a continuous homomorphism of $G$ and
$c:\Gamma{\times}X{\rightarrow}C$ is a cocycle taking values in a
compact group centralizing $\pi(G)$.
\end{theorem}

\noindent The principal improvements over earlier results is that
we do not need to pass to a finite ergodic extension of the action
and that we consider the case where $k$ is a non-Archimedean
fields of characteristic $0$. This builds on work of the second
author, Zimmer, Stuck, Lewis, Lifschitz, Venkataramana and others
\cite{L,Li,M1,Ma,M2,Z1,ZA,Z6,Stu,V}. In the case where $S$ is a
single point, Theorem \ref{theorem:Gammasuperrigidity} is
equivalent to the fact  that all homomorphisms of $\Gamma$ to
algebraic groups are superrigid. Theorem
\ref{theorem:Gsuperrigidity} is equivalent to the same fact when
applied to $S=G/{\Gamma}$.

\noindent
{\bf Remark:}When $k$ is non-Archimedean, it is not
always the case that the algebraic hull of the cocycle is
reductive unlike the case $k={\mathbb R}$ treated in \cite{Z6}.

\noindent{\bf Remark:} We also prove a result showing uniqueness
of the homomorphism $\pi$ occurring in Theorems
\ref{theorem:Gsuperrigidity} and \ref{theorem:Gammasuperrigidity}.
See subsection \ref{subsection:uniqueness} for details.

\noindent
{\bf Remark:}   Most of the results here should be true
for $k_i$ and $k$ of positive characteristic as well, though
additional arguments, similar to those in \cite{V,Li} are
apparently required. Some partial results in this direction are in
\cite{Li}.

The main applications of our results on perturbations of constant
cocycles are to studying perturbations of affine actions of $G$
and $\Gamma$.

\begin{defn}
\label{definition:affine}
 {\bf 1} Let $A$ and $D$ be topological
groups, and $B<A$ a closed subgroup. Let
$\rho:D{\times}A/B{\rightarrow}A/B$ be a continuous action. We
call $\rho$ {\em affine}, if, for every $d{\in}D$ there is a
continuous automorphism $L_d$ of $A$ and an element $t_d{\in}A$
such that $\rho(d)[a]=[t_d{\cdot}(L_d(a))]$.

\smallskip
\noindent {\bf 2} Let $A$ and $B$ be as above.  Let $C$ and $D$ be
two commuting groups of affine diffeomorphisms of $A/B$, with $C$
compact. We call the action of $D$ on $C{\backslash}A/B$ a {\em
generalized affine action}.

\smallskip
\noindent {\bf 3} Let $A$, $B$, $D$ and $\rho$ be as in $1$ above.
Let $M$ be a compact Riemannian manifold and
$\iota:D{\times}A/B{\rightarrow}\Isom(M)$ a $C^1$ cocycle.  We
call the resulting skew product $D$ action on $A/B{\times}M$ a
{\em quasi-affine action}. If $C$ and $D$ are as in $2$, and
$\alpha:D{\times}C{\backslash}A/B{\rightarrow}\Isom(M)$ is a $C^1$
cocycle, then we call the resulting skew product $D$ action on
$C{\backslash}A/B{\times}M$ a {\em generalized quasi-affine
action}.
\end{defn}

Our notion of {\em generalized affine action} is from \cite{F}.
The main application of our results on local rigidity of constant
cocycles is as part of our work on local rigidity of volume
preserving quasi-affine actions of $G$ and $\Gamma$ on compact
manifolds. We believe that volume preserving generalized
quasi-affine actions on compact manifolds are locally rigid as
well. As evidence for this, we have the following local entropy
rigidity result.  For any measure preserving action $\rho$ of $D$,
we denote by $h_{\rho}(d)$ the entropy of $\rho(d)$.  Let $\Ha$ be
an algebraic group defined over $\Ra$.  We will refer to the
connected component of the identity in $\Ha(\Ra)$ as a connected
real algebraic group.

\begin{corollary}
\label{corollary:qz}  Let $H$ be a connected real algebraic group,
$\Lambda<H$ a cocompact lattice and $K<H$ a compact subgroup. Let
$D=G$ or $\Gamma$ be as above and let $\rho$ be a $C^2$
generalized affine action of $D$ on $K{\backslash}H/{\Lambda}$.
Let $\rho'$ be any $C^2$ action sufficiently $C^1$ close to
$\rho$. Then $h_{\rho}(d)=h_{\rho'}(d)$ for all $d{\in}D$.
\end{corollary}

This result generalizes the one in \cite{QZ}.  Given the
description of generalized standard affine actions below, the
proof in \cite{QZ} actually applies.  We will prove Corollary
\ref{corollary:qz} as a corollary of (part of) the proof of
Theorems \ref{theorem:localcocyclerigidityG} and
\ref{theorem:localcocyclerigidityGamma}.

We note here that our techniques prove local rigidity results for
perturbations of more general cocycles over actions of $G$ and
$\Gamma$ than those in Theorems
\ref{theorem:localcocyclerigidityG} and
\ref{theorem:localcocyclerigidityGamma}. We can prove an analogous
theorem for perturbations of cocycles that are products of compact
valued cocycles with constant cocycles. More generally, the
original cocycle and the perturbed cocycle need not be cocycles
over the same action, but only over actions that are ``close".
For example if $S$ is a topological space, then the actions being
$C^0$ close is sufficient.  (Since constant cocycles are cocycles
over any action, one need only consider a single action in the
formulations of Theorems \ref{theorem:localcocyclerigidityG} and
Theorem \ref{theorem:localcocyclerigidityGamma}.) The proof of
Corollary \ref{corollary:qz} then implies a local entropy rigidity
result for generalized quasi-affine actions of $G$ and $\Gamma$.
The interested reader is welcome to adjust the proofs below to
cover these situations, but for the sake of clarity we have
restricted to the generality that we need for our next set of
applications.

We now state a theorem which is used in our work on local rigidity
of quasi-affine actions \cite{FM1,FM2}.  This theorem shows that
any perturbation of any quasi-affine action is continuously
semi-conjugate back to the original action, at least ``along
hyperbolic directions".

Let $H$ be a connected real algebraic group and $\Lambda<H$ a
discrete cocompact subgroup and let $D$ be either $G$ or $\Gamma$.
Let $\rho$ be a quasi-affine action of $D$ on
$H/{\Lambda}{\times}M$ which lifts to $H{\times}M$.  By the
discussion in section \ref{section:cocycles} there is a unique
subgroup $Z$ in $H$ which is the maximal subgroup of $H$ such that
the derivative of $\rho$ on $Z$ cosets is an isometry for an
appropriate choice of metric on $H/{\Lambda}$. The description
given there shows that the lift of $\rho$ to $H{\times}M$ descends
to an action $\bar \rho$ on $Z{\backslash}H$.  For example, if
$G<H$ acts on $H/{\Lambda}$ by left translations, then $Z=Z_H(G)$.

\begin{theorem}
\label{theorem:semiconjugacy}   Let
$H/{\Lambda}{\times}M,\rho,D,Z$ and $\br$ be as in the preceding
paragraph. Given any action $\rho'$ sufficiently $C^1$ close to
$\rho$, there is a continuous $D{\times}{\Lambda}$ equivariant map
$f:(H{\times}M, \rho'){\rightarrow}(Z{\backslash}H, \br)$, and $f$
is $C^0$ close to the natural projection map.
\end{theorem}

\noindent For actions by left translations this follows from
Theorems \ref{theorem:localcocyclerigidityG} and
\ref{theorem:localcocyclerigidityGamma}. To prove Theorem
\ref{theorem:semiconjugacy} as stated here, we need a stronger
result which is Theorem \ref{theorem:localcocyclerigidity} in
section \ref{section:generalresults}.  Theorem
\ref{theorem:semiconjugacy} holds more generally for any skew
product action of $D$ on $H/{\Lambda}$ which is affine on
$H/{\Lambda}$ and given by a cocycle
$\iota:D{\times}H/{\Lambda}{\rightarrow}\Diff^1_{\omega}(M)$ where
$\omega$ is a volume form on $M$ and $M$ is compact. The version
stated here is what is needed in \cite{FM1}.  We note that, by
Theorems \ref{theorem:describingactionsG} and
\ref{theorem:describingactionsgamma} below any quasi-affine $D$
action on $H/\Lambda{\times}M$ lifts to $H{\times}M$ on a finite
index subgroup $D'<D$.

Theorem \ref{theorem:localcocyclerigidityG}, Theorem
\ref{theorem:localcocyclerigidityGamma} and their applications
hold in a wider setting than the groups $G$ and $\Gamma$ discussed
above. The proof uses only that the cocycles we are considering
satisfy the conclusion of the cocycle superrigidity theorems and
that the group $G$ has ``few" representations. For example, for
$Sp(1,n)$, $F_4^{-20}$ and their lattices, our techniques can be
combined with the results of \cite{CZ} to obtain local rigidity
theorems for certain perturbations of certain cocycles of these
groups.  If variants of Theorems \ref{theorem:Gsuperrigidity} and
\ref{theorem:Gammasuperrigidity} hold for $Sp(1,n)$ and
$F_4^{-20}$ and their lattices, then Theorems
\ref{theorem:localcocyclerigidityG},
\ref{theorem:localcocyclerigidityGamma} and
\ref{theorem:localcocyclerigidity} hold for these groups as well.

In section \ref{section:preliminaries} we collect various standard
definitions used throughout the paper.  Section
\ref{section:superrigidity} concerns superrigidity for cocycles.
Section \ref{section:reporbits} proves that certain orbits in
representation varieties are closed. Section
\ref{section:generalresults} contains the proof of our main
results.  The final section of the paper contains the proofs of
Corollary \ref{corollary:qz} and Theorem
\ref{theorem:semiconjugacy}.  This section also contains a
detailed description of all affine actions of $G$ and $\Gamma$ as
above.

\section{\bf Preliminaries}
\label{section:preliminaries}

We now collect various definitions that will be used in the course of the paper.

\subsection{\bf Algebraic groups}
\label{section:prelim3} In this paper the words ``algebraic group"
mean a linear algebraic group defined over a local field $k$ in
the sense of \cite{B}.  Unless otherwise noted, throughout this
paper $k$ will be a local field of characteristic zero.
For background on algebraic groups particularly relevant to what
follows, we refer the reader to \cite[I.1-2]{M2}.

\subsection{\bf Cocycles and ergodic theory}
\label{section:prelim1}

Given a group $D$, a space $X$ and an action
$\rho:D{\times}X{\rightarrow}X$, we define a {\em cocycle over the
action} as follows.  Let $L$ be a group, the cocycle is a map
$\alpha:D{\times}X{\rightarrow}L$ such that
$\alpha(g_1{g_2},x)=\alpha(g_1, g_2x)\alpha(g_2,x)$ for all
$g_1,g_2{\in}D$ and all $x{\in}X$. The regularity of the cocycle
is the regularity of the map $\alpha$. If the cocycle is
measurable, we only insist on the equation holding almost
everywhere in $X$.  Note that the cocycle equation is exactly what
is necessary to define a skew product action of $D$ on
$X{\times}L$ or more generally an action of $D$ on $X{\times}Y$ by
$d(x,y)=(dx,\alpha(d,x)y)$ where $Y$ is any space with an $L$
action.

We say two cocycles $\alpha$ and $\beta$ are {\em cohomologous} if
there is a map $\phi:X{\rightarrow}L$ such that
$\alpha(d,x)={\phi(dx)}{\inv}\beta(d,x){\phi(x)}$.  Again we can
define the cohomology relation in any category, depending on how
much regularity we seek or can obtain on $\phi$. A cocycle is
called {\em constant} if it does not depend on $x$, i.e.
$\alpha_{\pi}(d,x)=\pi(d)$ for all $x{\in}X$ and $d{\in}D$.  One
can easily check from the cocycle equation that this forces the
map $\pi$ to be a homomorphism $\pi:D{\rightarrow}L$.  When
$\alpha$ is cohomologous to a constant cocycle $\alpha_{\pi}$ we
will often say that $\alpha$ is cohomologous to the homomorphism
$\pi$.  The cocycle superrigidity theorems imply that many
cocycles are cohomologous to constant cocycles, at least in the
measurable category.

A measurable cocycle $\alpha:D{\times}S{\rightarrow}L$ is called
{\em strict} if it is defined for all points in $D{\times}S$ and
the cocycle equation holds everywhere instead of almost
everywhere.  For a dictionary translating facts about strict
cocycles on homogeneous $D$-spaces to facts about homomorphisms of
subgroups of $D$, see \cite[Section 4.2]{Z1}.

An action of a group $D$ on a topological space $X$ is called {\em
tame} if the quotient space $D{\backslash}X$ is $T_0$, i.e. if for
any two points in $D{\backslash}X$, there is an open set around
one of them not containing the other.

Given a locally compact group $D$ and a discrete subgroup
$\Gamma<D$, there is a particularly important strict cocycle
$\beta_X:D{\times}D/{\Gamma}{\rightarrow}{\Gamma}$.  We define
this by choosing a fundamental domain $X$ for the $\Gamma$ action
on $D$.  By this we mean that there is a unique representation
$d=\omega(d)\tau(d)$ where $\omega(d)$ is in $X$ and $\tau(d)$ is
in $\Gamma$.  Identifying $D/{\Gamma}$ with $X{\subset}D$, we
define $\beta_X(d,x)=\tau(dx){\inv}$.  This cocycle is of
particular interest when $\Gamma<D$ is a lattice.  We call
$\beta_X$ the {\em strict cocycle corresponding to the fundamental
domain X.}

Let $D$ be a compactly generated group, with compact generating
set $K$.  Let $A$ be a metrizable, locally compact group and fix a
distance function $d:A{\times}A{\rightarrow}\Ra$. Given two
measurable cocycles $\alpha,\beta:D{\times}S{\rightarrow}A$ into a
locally compact group $A$, we can define a measurable function on
$S$ by $d(\alpha(d,x),\beta(d,x))$. We say that $\alpha$ and
$\beta$ are {\em $L^{\infty}$ close} if there exists a small
$\varepsilon>0$ such that
$\|d(\alpha(k,x),\beta(k,x))\|_{\infty}<{\varepsilon}$ for any
$k{\in}K$.

 Let $\La$ be an algebraic $k$-group and $L=\La(k)$.
Let a group $D$ act ergodically on a measure space $S$ and let
$\alpha:D{\times}S{\rightarrow}L$ be a cocycle. There is a unique
(up to conjugacy), minimal algebraic subgroup $\Ha$ in $\La$ such
that $\alpha$ is cohomologous to a cocycle taking values in
$H=\Ha(k)$. The group $H$ is referred to as the {\em algebraic
hull} for the cocycle. This is a generalization the Zariski
closure of a subgroup of an algebraic group. For more details, see
chapter 9 of \cite{Z1}.

We recall that given any group $D$ acting on a compact metric
space $X$ preserving a Borel measure $\mu$, there is an {\em
ergodic decomposition} of $\mu$.  That is, there are Borel
measures $\mu_i$ on $X$, where each $\mu_i$ is an invariant
ergodic measure for the action of $D$, and the measure $\mu$ is
obtained as an integral of the $\mu_i$ over a specific measure
$\tilde \mu$ on the space of measures on $X$. Furthermore, the
measures $\mu_i$ are mutually singular.

\subsection{\bf The space of actions}
\label{section:prelim2}

In the introduction, some statements are made about actions being
$C^k$ close.  Let $D$ be a locally compact topological group and
$X$ a smooth manifold.
 Since an action is a map $D{\rightarrow}\Diff^k(X)$ we can topologize the space of actions by taking the compact
open topology on $\Hom(D, \Diff^k(X))$.  Two actions are $C^k$
close if they are close with respect to this topology.  If $D$ is
compactly generated with compact generating set $K$, this means
that $\rho$ and $\rho'$ are $C^k$ close if and only if
$\rho(d){\circ}\rho'(d){\inv}$ is in a small neighborhood of the
identity in $\Diff^k(X)$ for all $d{\in}K$.  Given a manifold or a
space $X$ equipped with an action $\rho$, we often write
$(X,\rho)$ to denote the space with the action. Similarly a map
written $(X,\rho){\rightarrow}(X', \rho')$ is a map of $D$-spaces
or a $D$ equivariant map.

\section{{\bf Superrigidity for Cocycles}}
\label{section:superrigidity}

In this section we prove Theorems \ref{theorem:Gsuperrigidity} and
\ref{theorem:Gammasuperrigidity} as well as some related results.
Our integrability condition allows us to use Oseledec'
Multiplicative Ergodic Theorem to obtain our general result.  Some
partial results below do not require the integrability condition.
Theorem \ref{theorem:Gammasuperrigidity} is deduced from Theorem
\ref{theorem:Gsuperrigidity}.  The proof of Theorem
\ref{theorem:Gsuperrigidity} requires that one first argue the
case where $L$ is semi-simple and then use the result in that case
to prove the more general result.

Theorems \ref{theorem:Gsuperrigidity} and
\ref{theorem:Gammasuperrigidity} imply a general result on the
algebraic hull of the cocycles considered. In fact, at least for
$G$ cocycles, this result is a step in the proof of Theorem
\ref{theorem:Gsuperrigidity}, see Theorem
\ref{theorem:algebraichull}.  It is proved in \cite{M2} that for
any field $k$ and any homomorphism
$\pi:\Gamma{\rightarrow}\La(k)$, the Zariski closure of
$\pi(\Gamma)$ is semisimple.  This is equivalent to saying that
the algebraic hull of the cocycle
$\pi{\circ}{\beta}:G{\times}G/{\Gamma}{\rightarrow}L$ is
semisimple.  In \cite{Z6}, it is shown that if $k={\mathbb R}$,
any $G$-integrable cocycle $\alpha:G{\times}X{\rightarrow}L$ has
algebraic hull reductive with compact center.  If $k$ is
non-Archimedean, it is no longer the case that the algebraic hull
is reductive.  The following example shows that our results on the
algebraic hull are sharp.
\begin{example}
\label{example:algebraichullnonreductive}{\em We let $J$ be a
finite index set and for each $j{\in}J$, we let $k_j$ be a local
field of characteristic zero and ${\mathbb H}_j$ be a connected
simply connected semisimple algebraic $k_j$-group. We let
$H_j={\mathbb H}_j(k_j)$ the $k_j$-points of ${\mathbb H}_j$ and
$H=\prod_{j{\in}J}H_j$. We further assume that there is an
irreducible lattice $\Lambda<H$. For many examples where
irreducible $\Lambda$ exist, we refer the reader to
\cite[IX.1.7]{M2}. Let $\pi:G{\rightarrow}H$ be a homomorphism and
assume that $Z_H(\pi(G))$ contains a non-trivial unipotent
subgroup $U<H_l$ for some $l{\in}J$ where $k_l$ is
non-Archimedean. (We leave the easy construction of explicit
examples to the reader.) Let $K<U$ be a Zariski dense compact
subgroup and consider the $G$ action on $K{\backslash}H/{\Lambda}$
and $H/{\Lambda}$.  Choosing a measurable trivialization of the
$K$-bundle $H/{\Lambda}{\rightarrow}K{\backslash}H/{\Lambda}$
defines a cocycle
$\alpha:G{\times}K{\backslash}H/{\Lambda}{\rightarrow}K$, which we
view as $\alpha:G{\times}K{\backslash}H/{\Lambda}{\rightarrow}U$
via the inclusion of $K<U$.  Standard arguments using Mautner's
Lemma show that the $G$ actions on $H/{\Lambda}$ and
$K{\backslash}H/{\Lambda}$ are ergodic.  A simple argument using
the fact that the Mackey range of the cocycle $\alpha$ is
$H/{\Lambda}$ and ergodicity of the $G$ action on $H/{\Lambda}$
shows that $U$ is the algebraic hull of $\alpha$.  See
\cite[4.2.24]{Z1} for definitions and discussion of the Mackey
range.

The reader should note the following
\begin{enumerate}
\item the above construction yields the same results when applied
to the restriction of the actions and cocycles to any lattice
$\Gamma<G$;

\item the construction gives non-trivial examples even when
$G=\Ga(\Ra)$;

\item one can take products of cocycles constructed as above with
constant cocycles to obtain cocycles whose algebraic hull is
neither unipotent nor reductive;

\item  the argument above works for more general subgroups
$Z<Z_H(\pi(G)){\cap}H_l$ where $K<Z$ is a Zariski dense compact
subgroup.  One can construct examples where $Z=F{\ltimes}U$ is a
Levi decomposition and the $F$ action on $U$ is non-trivial.
\end{enumerate}
\noindent Let $\La$ be an algebraic group over $k_l$ and
$L=\La(k_l)$ and $D=G$ or $\Gamma$.  The above outline constructs
cocycles $\alpha:D{\times}S{\rightarrow}L$ of the form
$\alpha=\pi{\cdot}c$ where $\pi:G{\rightarrow}L$ is a continuous
homomorphism and $c:D{\times}S{\rightarrow}C$ is a cocycle taking
values in a compact group $C<Z_L(\pi(G))$. We can construct
$\alpha$ with algebraic hull $L$ for any $\La$ provided we choose
$\pi$ so that $\pi(G)$ commutes with the unipotent radical of
$\La$.}
\end{example}
We now briefly indicate the plan of this section. Subsection
\ref{subsection:superrigiditywarmup} fixes notation for all of
section \ref{section:superrigidity} and contains some technical
lemmas used throughout. In subsection
\ref{subsection:keyreduction} we prove a key technical result
which shows that certain cocycles are cohomologous to constant
cocycles. Subsection \ref{subsection:semisimplehull} applies the
results of subsection \ref{subsection:keyreduction} to prove a
variant of Theorem \ref{theorem:Gsuperrigidity} where the
algebraic hull of the cocycle is assumed to be semisimple.
Subsection \ref{subsection:characteristicmaps} proves some
conditional results on $G$-integrable cocycles, again using the
results from subsection \ref{subsection:keyreduction}. We show how
to use property T to control cocycles into amenable and reductive
groups in subsection \ref{subsection:useofT} and then prove
Theorem 1.4 in subsection \ref{subsection:proofs}. Theorem 1.5 is
also proven in subsection \ref{subsection:proofs} modulo some
facts concerning $G$-integrability of certain induced cocycles.
These facts are then proven in subsection
\ref{subsection:integrability}.  Subsection
\ref{subsection:uniqueness} concerns the uniqueness of the
homomorphism $\pi$ in Theorems \ref{theorem:Gsuperrigidity} and
\ref{theorem:Gammasuperrigidity}.  These results are used in
subsection \ref{subsection:funnycocycles} to prove some results on
cocycles with constrained projections. The result on cocycles with
constrained projections is required to prove Theorem
\ref{theorem:localcocyclerigidity} which is used in the proof of
Theorem \ref{theorem:semiconjugacy}.

\subsection{Notations and reductions.}
\label{subsection:superrigiditywarmup}

In this subsection, we fix notations and definitions for all of
section \ref{section:superrigidity}.  We also prove some technical
lemmas that are used throughout this section.

The group $G$ will be as specified in the introduction, but we
both weaken the rank assumption and make some preliminary
reductions. Let $\mathcal{S}$ be the union of primes of $\mathbb
Z$ and $\{\infty\}$ and let ${\mathbb Q}_p$ be the $p$-adic
completion of $\mathbb Q$, where as usual, ${\mathbb
Q}_{\infty}={\mathbb R}$. By application of restriction of
scalars, we can assume that each $k_i={\mathbb Q}_{p_\alpha}$,
where the $p_i$ are distinct elements of the set $\mathcal{S}$. As
before, for the Archimedean factor, we can replace ${\mathbb
G}_i(\mathbb R)$ by it's topological universal cover. Actually
this can be done or not done for each simple factor independently,
though we simplify exposition by ignoring this nuance.  Instead of
assuming that each simple factor of ${\mathbb G}_i(k_i)$ has
$k_i$-rank at least two, we let $r_i=k_i$-rank$({\mathbb
G}_i(k_i))$ and define the {\em rank of $G$} as
$\sum_{i{\in}I}r_i$ and assume that the rank of $G$ is at least
two and that $G$ has no non-trivial compact factors (or,
equivalently, that every simple factor of $\Ga_i$ has $k_i$-rank
at least one).

 We specify a certain compact homogeneous
$G$ space, often called a {\em boundary} for $G$.  Let
$\Pa_i<\Ga_i$ be a minimal parabolic subgroup.  we define $P_i$ to
be $\Pa_i(k_i)$ if $G_i=\Ga_i(k_i)$. If $G_i$ is the topological
universal cover of $\Ga_i(k_i)$, we define $P_i$ to be the
pre-image of $\Pa_i(k_i)$ under the covering map from $G_i$ to
$\Ga_i(k_i)$. We let $P=\prod_{i{\in}I}P_i$ and the homogeneous
space we consider is $\BG$.  We note that the $G$ action on $\BG$
factors through the projection to $\prod_{i{\in}I}\Ga_i(k_i)$ and
the space $\BG$ can be identified with
$\prod_{i{\in}I}{\Pa_i}(k_i){\backslash}\Ga_i(k_i)$ which can be
identified as a variety with
$\prod_{i{\in}I}({\Pa_i}{\backslash}\Ga_i)(k_i).$

We fix $(S,\mu)$ to be a standard probability measure space. Also
$\La$ will denote an algebraic $k$-group and $L=\La(k)$.  We
denote by $\La^0$ the connected component of $\La$ and let
$L^0=\La^0(k)$.  As above, we apply restriction of scalars and
assume that $k={\mathbb Q}_p$ for some $p{\in}\mathcal{S}$.

By a simple factor of $G$, we mean a subgroup $F<G$ which is
either $\Fa(k_i)$ or its topological universal cover, where $\Fa$
is almost simple.  We note that under our hypothesis, $G$ is the
direct product of all of its simple factors. We say a simple
factor $F_i$ has rank one if the $k_i$ rank of $\Fa_i$ is one.  If
$F_i$ is a simple factor of $G$ then there is a group $F_i^c<G$
such that $G=F_i{\times}F_i^c$.  We call $F_i^c$ the complement of
$F_i$.

\begin{defn}
\label{defn:weaklyirreducible} Let $(S,\mu)$ be a finite measure
space.  Given a group $G$ acting ergodically on $S$ preserving
$\mu$, we call the action {\em weakly irreducible} if for any rank
one simple factor $F<G$, the complement $F^c$ acts ergodically on
$S$.
\end{defn}


If no simple factor of $G$ has rank $1$, this is equivalent to the
ergodicity of the $G$ action.  This is weaker than the standard
definition of irreducibility where it is assumed that all simple
factors act ergodically \cite{Z1}.  The definition of an
irreducible action is motivated by properties of irreducible
lattices.  We call a lattice $\Gamma<G$ {\em weakly irreducible}
if the projection of $\Gamma$ to any rank 1 factor of $G$ is
dense.  Standard arguments using the generalized Mautner
phenomenon, see \cite[II.3.3]{M2}, show that a lattice is weakly
irreducible if and only if the action of $G$ on $G/{\Gamma}$ is
weakly irreducible.

We will use the following elementary lemmas repeatedly. The first
is obvious.

\begin{lemma}
\label{lemma:ratioiscocycle} Let $A$ be a group and let
$\alpha:D{\times}S{\rightarrow}A$ and
$\beta:D{\times}S{\rightarrow}A$ be cocycles over the action of a
group $D$ on a set $S$.  Assume that $\beta(D{\times}X)$ is
contained in a subgroup $B<A$ and let $Z=Z_A(B)$.  Let
$\eta:D{\times}S{\rightarrow}Z$ be a map.   If
$\alpha(d,x)=\beta(d,x)\eta(d,x)$, then $\eta$ is a cocycle over
the $D$ action.
\end{lemma}

We let $\tau_i:G{\rightarrow}{\Ga_i(k_i)}$ be the natural
projection.

\begin{lemma}
\label{lemma:liftingreps} Given a non-trivial continuous
homomorphism $\pi:G{\rightarrow}L$ there is $i{\in}I$ such that
$k=k_i$ and a $k$-rational homomorphism
$\pi_i:\Ga_i{\rightarrow}\La$ such that $\pi=\pi_i{\circ}\tau_i.$
From this we can deduce:
\begin{enumerate}
\item the Zariski closure of $\pi(G)$ is semisimple and connected
and; \item if $\La'{\rightarrow}\La$ is a $k$-isogeny, then $\pi$
lifts to a continuous homomorphism $\pi':G{\rightarrow}\La'(k)$
\end{enumerate}
\end{lemma}

\begin{proof}
We first give the proof where all $G_i$ are $k_i$-points of
algebraic $k_i$-groups and then describe the modifications
necessary when $G_i$ is the topological universal of such a group.

Let the projection from $G$ to $G_i$ be ${\tau}_i$. Since
$k={\mathbb Q}_p$, by \cite[I.2.6]{M2} any continuous homomorphism
of any $G_i$ into $L$ is the restriction of rational map from
$\Ga_i$ to $\La$. This implies there is an $i$ and a rational
homomorphism $\bar \pi:\Ga_i{\rightarrow}\La$ such that $\pi$ is
the restriction
 of ${\tau}_i{\circ}{\bar \pi}$.

Since $\Ga_i$ is connected and semisimple and the characteristic
of $k$ is zero, it follows that the Zariski closure of ${\bar
\pi}(\Ga_i)$ is connected and semisimple.   If
$\La'{\rightarrow}\La$ is an isogeny, then $\bar \pi$ lifts to a
map ${\bar \pi}':\Ga{\rightarrow}\La'$ since $\Ga_i$ is simply
connected.

Now assume that $G_i$ is the topological universal cover of
$\Ga_i(\mathbb R)$. If $k{\neq}{\mathbb R}$ then any continuous
homomorphism from $G_i$ to $L$ is trivial, so we are done by the
discussion above.  If $k={\mathbb R}$ then $\pi$ factors through a
continuous homomorphism ${\bar \pi}:G_i{\rightarrow}L$.  The image
of $\bar \pi$ is a closed subgroup of $L$ and so is the real
points of a real algebraic subgroup.  This implies that $\bar \pi$
factors through the covering map $G_i{\rightarrow}\Ga_i(\mathbb
R)$.  The conclusions of the lemma now follow as before.
\end{proof}

\subsection{$\alpha$-invariant maps into algebraic varieties}
\label{subsection:keyreduction}

Given two $G$-spaces $S$ and $Y$, an $L$ space $R$ and a cocycle
$\alpha:G{\times}S{\rightarrow}L$, we call a map
$f:Y{\times}S{\rightarrow}R$ {\em $\alpha$-invariant} if
$f(gy,gs)=\alpha(g,s)f(y,s)$ for all $g$ and almost every $(y,s)$.
Note that this definition differs slightly from the one in
\cite{Z1}, where this map would be called $\tilde
\alpha$-invariant where $\tilde \alpha$ is the pullback of
$\alpha$ to $G{\times}Y{\times}S$.

The following theorem will play a key role in all proofs in this
section.  The assumption on the rank of $G$ is only used to be
able to apply this theorem.

\begin{theorem}
\label{theorem:rationalimpliesrep} Assume $G$ acts weakly
irreducibly on $S$ preserving $\mu$. Let $M$ be the $k$ points of
an algebraic variety $\Ma$ defined over $k$ on which $\La$ acts
$k$ rationally. Assume that $\alpha:G{\times}S{\rightarrow}L$ is a
Borel cocycle whose algebraic hull is $L$ and that there exists a
measurable $\alpha$-invariant map
$\phi:\BG{\times}S{\rightarrow}M$ such that the essential image of
$\phi$ is not contained in the set of $L$-fixed points of $M$.
Then there is a normal $k$-subgroup $\Ha<\La$ of positive
codimension such that:

\begin{enumerate}


\item  $p_H{\circ}\alpha$ is cohomologous to a continuous
homomorphism $\pi_H:G{\rightarrow}L/H$, where
$p_H:L{\rightarrow}L/H$ and $H=\Ha(k)$;

\item $\La/\Ha$ is semisimple and connected.

\end{enumerate}
\end{theorem}

\begin{proof}

Let $\phi_s(x)=\phi(x,s)$ where $x{\in}\BG$ and $s{\in}X$.  First
one shows that either $\phi_s$ is rational for almost every  $s$
or $\phi_s$ is constant for almost every $s$.  By rational we mean
that there is $i{\in}I$ such that $k=k_i$, the map $\phi$ factors
through the projection
$p_i:\BG{\rightarrow}{\Pa}_i(k_i){\backslash}{\Ga}_i(k_i)$ which
means that $\phi=p_i{\circ}{\bar \phi}$ where that $\bar \phi$ is
a $k$ rational map ${\Pa}_i{\backslash}{\Ga}_i{\rightarrow}\Ma$.
Rationality of $\phi$ was shown by Zimmer in \cite{ZA} for
irreducible actions with each $G_i=\Ga_i(k_i)$ using an adaptation
of an argument due to the second author \cite{M1,Ma}. The proof
goes through almost verbatim for weakly irreducible actions, as
well as for the case where one $G_i$ is the universal cover of
$\Ga_i(\Ra)$. See also pages 104-5 of \cite{Z1} or \cite{Fu} for
accessible presentations of special cases. Our definition of weak
irreducibility is motivated by the ergodicity needed at this step
of the proof. We now assume that $\phi_s$ is rational and proceed
in this case, the case of $\phi_s$ constant is discussed at the
end of the proof.

Secondly, one sees that the map $\Phi:S{\rightarrow}\Rat(\BG,M)$
defined by $\Phi(s)=\phi_s$ takes values in a single orbit.  This
follows from tameness of the $G{\times}L$ action on $\Rat(\BG,M)$
and the ergodicity of the $G$ action on $S$, see the ``Proof of
Step 3" on pages 105-6 and also Proposition 3.3.2 of \cite{Z1}.

One now picks a rational map $\psi$ in this orbit and defines a
map $l:S{\rightarrow}L$ such that $\phi_s=l(s)\psi$. Letting
$\beta(g,s)=l(gs){\inv}\alpha(g,s)l(s)$ we have that
$\beta(g,s)\psi(x)=\psi(gx)$.  Let $H$ denote the point-wise
stabilizer of $\psi(\BG)$ in $M$.  Since $M=\Ma(k)$ and $\La$ acts
rationally on $\Ma$, $H=\Ha(k)$ where $\Ha<\La$ is an algebraic
subgroup defined over $k$.   Since $\beta(g,s)\psi(x)=\psi(gx)$,
the Zariski closure of $\psi(\BG)$ is invariant under
$\beta(G{\times}S)$ and since the algebraic hull of $\beta$ is
$L$, the Zariski closure of $\psi(\BG)$ is $L$-invariant.
Therefore $H$ is normal in $L$, and $\Ha$ is normal in $\La$.
Fixing (almost any) $s$, and writing $\beta_s(g)=\beta(g,s)$, we
have that
$\beta_s(g_1g_2)\psi(x)=\psi(g_1g_2x)=\beta_s(g_1)\psi(g_2x)=\beta_s(g_1)\beta_s(g_2)\psi(x)$.
Therefore $\beta_s(g_1g_2)\beta(g_2){\inv}\beta(g_1){\inv}$ fixes
$\psi(\BG)$ pointwise.  It follows that
$p_H{\circ}{\beta}_s:G{\rightarrow}L/H$ is a homomorphism. That
$\pi=p_H{\circ}\beta_s$ is continuous follows from a result of
Mackey, see \cite[B.3]{Z1}. The remaining conclusions of the
theorem follow from Lemma \ref{lemma:liftingreps}.

If $\phi$ is constant for almost every $s{\in}S$, we have an
$\alpha$-invariant map $\Phi:S{\rightarrow}M$.  The image of this
map is contained in a single $H$ orbit since the $L$ action on $M$
is tame and the $G$ action on $S$ is ergodic.  Since the $L$
action on $M$ is defined by an algebraic action of $\La$ on $\Ma$,
the stabilizer of this orbit is $H={\Ha}(k)$ where $\Ha<\La$ is an
algebraic subgroup.  This means that we have an $\alpha$ invariant
map ${\phi}:S{\rightarrow}L/H$.  By \cite[Lemma 5.2.11]{Z1}, this
implies that the cocycle $\alpha$ is equivalent to one taking
values in $H$, which contradicts our assumption on the algebraic
hull of the cocycle unless $H=L$ in which case we contradict our
assumption that the essential image of $\phi$ is not contained in
the set of $L$-fixed points.
\end{proof}

The reader should note that essentially the same result can be
proven by considering equivariant measurable maps from
$G{\times}X$ to vector spaces, as in \cite[VII.1-4]{M2}.  The
argument there requires some modification since, in the language
of that text, one needs to consider maps that are not strictly
effective.

\subsection{\bf Algebraic hull semisimple}
\label{subsection:semisimplehull}

We now prove Theorems \ref{theorem:Gsuperrigidity} and Theorem
\ref{theorem:Gammasuperrigidity} in the case where the algebraic
hull of the cocycle is semisimple.

\begin{theorem}
\label{theorem:Gsuperrigidity1} Let $G$ act weakly irreducibly on
$S$ preserving $\mu$ and let $\alpha:G{\times}S{\rightarrow}L$ be
a Borel cocycle with algebraic hull $L$. Further assume that $\La$
is semisimple. Then $\alpha$ is cohomologous to a cocycle
$\beta=\pi{\cdot}c$. Here ${\pi:G{\rightarrow}L}$ is a continuous
homomorphism and $c:G{\times}S{\rightarrow}C$ is a cocycle taking
values in a compact group $C$ centralizing $\pi(G)$.
\end{theorem}

\begin{theorem}
\label{theorem:Gammasuperrigidity1}  Let ${\Gamma}<G$ be a weakly
irreducible lattice. Assume $\Gamma$ acts ergodically on $S$
preserving $\mu$. Let $\alpha:\Gamma{\times}S{\rightarrow}L$ be a
Borel cocycle with algebraic hull $L$.  Further assume $\La$ is
semisimple.  Then $\alpha$ is cohomologous to a cocycle $\beta$
where $\beta(\gamma,x)=\pi(\gamma){c(\gamma,x)}$. Here
${\pi:G{\rightarrow}L}$ is a continuous homomorphism and
$c:\Gamma{\times}S{\rightarrow}C$ is a cocycle taking values in a
compact group centralizing $\pi(G)$.
\end{theorem}

To reduce Theorem \ref{theorem:Gsuperrigidity1} to Theorem
\ref{theorem:rationalimpliesrep} we need to find a $k$ variety $M$
on which $H$ acts $k$-rationally and an $\alpha$-invariant map
$f:\BG{\times}S{\rightarrow}M$. To produce $\alpha$ invariant
maps, one uses the following modification of a lemma of
Furstenberg from \cite{Fu2} which can be deduced from Propositions
$4.3.2, 4.3.4$ and $4.3.9$ of \cite{Z1}. The lemma holds under
more general circumstances than those needed here.  For the lemma,
$G$ can be a locally compact, $\sigma$-compact group, $P$ a closed
amenable subgroup and $L$ any topological group.

\begin{lemma}
\label{lemma:boundarymap} Assume $G$ acts on $S$ preserving $\mu$.
Let $\alpha:G{\times}S{\rightarrow}L$ be a Borel cocycle. Let $B$
be any compact metrizable space on which $L$ acts continuously and
$\mathcal{P}(B)$ the space of Borel regular probability measures
on $B$. Then there is an $\alpha$-invariant map
$f:\BG{\times}S{\rightarrow}\mathcal{P}(B)$.
\end{lemma}

We note here that we give a proof using only amenability of $P$,
without reference to the notion of an amenable actions, though we
do rely on ideas of Zimmer's to construct a convex compact space
on which $P$ acts affinely and continuously.

\begin{proof}
Let $B$ be any compact $L$-space.   Via $\alpha$ we can define a
skew product action of $G$ on $S{\times}B$.  We consider the
diagonal $G$ action on $G{\times}S{\times}B$ given by the right
$G$ action on $G$ and the skew product action on $S{\times}B$,
which we note commutes with the left $G$ action on $G$. Let
$\mu_G$ be Haar measure on $G$ and
$\mathcal{M}(G{\times}S{\times}B)$ be the space of regular Borel
measures on $G{\times}S{\times}B$ which are invariant under the
diagonal action and project to $\mu_G{\times}\mu$ on $G{\times}S$.
We want to topologize $\mathcal{M}(G{\times}S{\times}B)$ so the
left $G$ action is continuous and the space is a compact convex
affine $G$-space. Using disintegration of measures, we can
identify $\mathcal{M}(G{\times}S{\times}B)$ with
$F(G{\times}S,\mathcal{P}(B))$ the space of measurable maps from
$G{\times}S$ to $\mathcal{P}(B)$. Let $C(B)$ be the Banach space
of continuous functions on $B$. We identify
$F(G{\times}S,\mathcal{P}(B))$ as a subset of
$L^{\infty}(G{\times}S, C(B)^*)$ and give $L^{\infty}(G{\times}S,
C(B)^*)$ the weak topology coming from the identification
$L^{\infty}(G{\times}S, C(B)^*)=L^1(G{\times}S, C(Y))^*$.  In this
topology the action of $G$ is continuous and
$F(G{\times}S,\mathcal{P}(B))$ is a closed, convex subset of the
unit ball in $L^{\infty}(G{\times}S, C(B)^*)$.  See \cite[Section
4.3]{Z1} for more discussion of this and related constructions. It
follows that there is a fixed point
$\mu^P{\in}\mathcal{M}(G{\times}S{\times}B)$ for the left $P$
action.  By applying disintegration of measures, this is left
$P$-invariant, $\alpha$-invariant map ${\tilde
f}:G{\times}S{\rightarrow}\mathcal{P}(B)$ or equivalently an
$\alpha$-invariant map $f:\BG{\times}S{\rightarrow}\mathcal{P}(B)$
\end{proof}

We will also need the following lemma essentially due to
Furstenberg.

\begin{lemma}
\label{lemma:measuresupports} Let $J<GL(V)$, where $V=k^n$ and $k$
is a local field of characteristic zero.  Let $J$ act on the
projective space $P(V)$ preserving a measure $\mu$.  Then either
$J$ is projectively compact (i.e. the image of $J$ in $PGL(V)$ is
compact) or there is a proper subspace $W<V$ with $\mu(W)>0$.
\end{lemma}

For a proof, we refer the reader to \cite[Lemma 3.2.2]{Z1} or the
original article of Furstenberg \cite{Fu1} in the case where
$k={\mathbb R}$.

\begin{proof}[Proof of Theorem \ref{theorem:Gsuperrigidity1}]
We call a representation of an algebraic group {\em almost
faithful} if the kernel of the representation is finite. We choose
an almost faithful irreducible $k$-rational representation
$\sigma$ of $\La$ on $V$ such that the restriction of $\sigma$ to
any $\La^0$ invariant subspace is almost faithful, where as usual
$\La^0$ denotes the connected component. (This can be done by
inducing an almost faithful irreducible $\La^0$ representation.)
Let $B=\Pa(V)$ be the corresponding projective space.

Since $P<G$ is an amenable subgroup Lemma \ref{lemma:boundarymap}
provides an $\alpha$ equivariant map
$f:\BG{\times}X{\rightarrow}\mathcal{M}(B)$.  In fact this map
takes values in a single $H$ orbit $\mathcal{O}$ in
$\mathcal{P}(B)$. This is deduced from ergodicity of the $G$
action on $S$ and the tameness of the action of $L$ on
$\mathcal{P}(B)$ which is \cite[Corollary 3.2.12]{Z1}.

Let $J$ be the stabilizer of a point $\mu$ for the $L$ action on
$\mathcal{O}$.  We prove that either $J$ is compact or $J$ is
contained in an algebraic subgroup of positive codimension in
$\La$.  If $\La$ is connected, this is Proposition 3.2.15 of
\cite{Z1}. By Lemma \ref{lemma:measuresupports}, if $J$ is not
projectively compact, then there is a proper subspace $W<V$ such
that $\mu([W])>0$. Since $\La$ is semisimple and the
representation on $V$ is almost faithful, the map from $\La$ to
$PGL(V)$ has finite kernel and only compact subgroups of $L$ are
projectively compact. Assuming $J$ is non-compact, we choose $W$
of minimal dimension among subspaces with positive $\mu$ measure.
Since the measure of $W$ is positive, $J{\cdot}W$ must be a finite
union of disjoint subspaces $\cup_{l=1}^n{W_l}$, and we let $\Fa$
be the stabilizer of the $J$ orbit $J{\cdot}W$.  Let $F=\Fa(k)$.
The stabilizer $J^W$ in $J$ of $W$ is of finite index in $J$ and,
by minimality of $W$ and Lemma \ref{lemma:measuresupports}, acts
on $P(W)$ via a homomorphism to a compact subgroup of $PGL(W)$. If
$\dim(\Fa)=\dim(\La)$, then the connected component of $\La$
preserves $\cup_{l=1}^n{W_l}$ and by connectedness preserves $W$.
Since $J^W<L$ and acts compactly on $W$, we then have that
$J^W{\cap}\La^0(k)$ is projectively compact. But, since we have
that $\La^0$ is semisimple and the representation of $\La^0$ on
$W$ is almost faithful, the map from $\La^0$ to $PGL(W)$ has
finite kernel.  This implies that $J^W{\cap}\La^0(k)$ is compact.
The group $J^W{\cap}\La^0(k)$ has finite index in $J^W$ which has
finite index in $J$, so in this case, $J$ is compact. Therefore
either $J$ is compact or $\Fa$ is of positive codimension in
$\La$.


If $J$ is compact, then Lemma $5.2.10$ of \cite{Z1} applies and
shows that the cocycle $\alpha$ is cohomologous to one with
bounded image. If $J<\Fa$ an algebraic subgroup of positive
codimension, then we compose
$\phi:\BG{\times}X{\rightarrow}\mathcal{O}$ with the projection
$p:\mathcal{O}{\rightarrow}L/F{\subset}({\La}/{\Fa})(k)$.  We note
that the set of $L$ fixed points in $L/F$ is empty so we can apply
Theorem \ref{theorem:rationalimpliesrep}.  This theorem produces a
normal $k$-subgroup of positive codimension $\Ha<\La$, such that
the projection of $\alpha$ to $\La(k)/\Ha(k)$ is cohomologous to a
continuous homomorphism $\pi$ of $G$. That theorem also implies
that $\La/\Ha$ is semisimple and connected. Since $\La$ is
semisimple, there is a connected normal subgroup $\Ha^c<\La$ such
that the map $\Ha^c{\rightarrow}\Ha/\La$ is an isogeny. Now
$\La=\Ha{\cdot}\Ha^c$ where $\Ha{\cap}\Ha^c$ is finite.  Because
$[\Ha,\Ha^c]{\subset}\Ha{\cap}\Ha^c$ which is finite and $\Ha^c$
is connected,  $\Ha$ and $\Ha^c$ commute.  By Lemma
\ref{lemma:liftingreps} we can lift $\pi$ to a homomorphism
$\pi':G{\rightarrow}\Ha^c$. It then follows that $\alpha$ is
cohomologous to $\pi'{\cdot}{\alpha'}$ where $\alpha'$ takes
values in $\Ha(k)$ and is a cocycle by Lemma
\ref{lemma:ratioiscocycle}.
 One can now replace $\alpha$ by $\alpha'$ and complete the proof of
the theorem by induction on the dimension of $L$.
\end{proof}

\begin{proof}[Proof of Theorem \ref{theorem:Gammasuperrigidity1}]
This is proved by inducing actions and cocycles, exactly as in
\cite[Theorem 9.4.14]{Z1}. We briefly outline the argument.  Given
a $\Gamma$ action on $S$ and a $\Gamma$ cocycle
$\alpha:\Gamma{\times}S{\rightarrow}L$, we consider the induced
$G$ action on $(G{\times}S)/{\Gamma}$ and a cocycle ${\tilde
\alpha}:G{\times}G/{\Gamma}{\times}S{\rightarrow}L$.  We define
$\tilde \alpha$ by taking a fundamental domain $X$ for $\Gamma$ in
$G$ and the strict cocycle
$\beta_X:G{\times}G/{\Gamma}{\rightarrow}{\Gamma}$ corresponding
to $X$ and letting ${\tilde
\alpha}(g,[g_0,x])=\alpha(\beta(g,[g_0]),x)$.  It is
straightforward to verify that weak irreducibility of $\Gamma$ and
ergodicity of the $\Gamma$ action imply that the induced $G$
action is weakly irreducible.  One then shows that the algebraic
hull of $\tilde \alpha$ is $L$ and applies Theorem
\ref{theorem:Gsuperrigidity1} to $\tilde \alpha$. Straightforward
manipulation allows one to deduce the desired conclusions
concerning $\alpha$.
\end{proof}

\subsection{\bf Conditional results using $G$-integrability of $\alpha$}
\label{subsection:characteristicmaps}

In this section we prove a conditional result concerning the
algebraic hull of $G$-integrable cocycles.  The assumption of
$G$-integrability is only used here.

Before stating our result, we fix some notation and assumptions.
We assume that $G$ has property T and that $G$ acts weakly
irreducible on $(S,\mu)$. Let $\alpha:G{\times}S{\rightarrow}L$ be
a $G$-integrable Borel cocycle and assume that $L$ is the
algebraic hull of the cocycle. We can write $\La=\Fa{\ltimes}\Ua$
where $\Fa$ and $\Ua$ are $k$-subgroups, $\Ua$ is the unipotent
radical of $\La$ and $\Fa$ is reductive. Let
$p_F:\La{\rightarrow}\Fa$ and be the natural projection. We assume
that the cocycle $p_F{\circ}{\alpha}$ is cohomologous to a cocycle
of the form $\pi{\cdot}c$ where $\pi:G{\rightarrow}F$ is a
continuous homomorphism and $c$ is a cocycle taking values in a
compact subgroup $C<Z_F(\pi(G))$.    We note that $\pi$ can be
viewed as defining a homomorphism of $G$ into $L$, and we let
$\alpha'$ be the cocycle cohomologous to $\alpha$ that projects to
$\pi{\cdot}c$.

\begin{theorem}
\label{theorem:algebraichull} Under the hypotheses discussed in
the preceding paragraph, $U$ commutes with $\pi(G)$.
\end{theorem}

We prove the theorem by contradiction.  The general scheme is as
follows.  If $U$ does not commute with $\pi(G)$ there exists a
$k$-rational action of $\La$ on a variety $\Ma$ and an
$\alpha$-invariant map $\phi$ into $\Ma(k)$ such that the
pointwise stabilizer $H=\Ha(k)$ of the image does not contain all
of $U$. Applying Theorem \ref{theorem:rationalimpliesrep} we
obtain a contradiction, since number $2$ of that theorem implies
that $\La/\Ha$ is semisimple  and this implies that $\Ua<\Ha$.

We will construct a measurable map $\phi$ that satisfies the
hypotheses of Theorem \ref{theorem:rationalimpliesrep} by using
Oseledec' multiplicative ergodic theorem. We will give an argument
that is close to the one in \cite[Section V.3-4]{M2}, but also
refer the reader to \cite{Z6} for a somewhat different approach.


Let $\fl$ be the Lie algebra of $L$. Let $Gr_j(\fl)$ be the
Grassmann variety of $j$ planes in $\fl$.  We have an action of
$L$ on $\fl$ by the adjoint representation which also defines an
action of $L$ on $Gr_j(\fl)$.

We look at the representation $\Ad_{\fl}{\circ}\pi$.

\begin{theorem}
\label{theorem:characteristicmap} Assume $\pi(G)$ does not commute
with $U$. Then there is an integer $0<m< \dim(\fl)$ and an
$\alpha$-equivariant measurable map
$\phi:P{\backslash}G{\times}S{\rightarrow}Gr_{m}(\fl)$ such that
the pointwise stabilizer of the image does not contain all of $U$.
\end{theorem}

\noindent Before proving Theorem \ref{theorem:characteristicmap},
we show how that result implies Theorem
\ref{theorem:algebraichull}.

\begin{proof}[Proof of Theorem \ref{theorem:algebraichull}]
We now apply Theorem \ref{theorem:rationalimpliesrep} to the map
$\phi$ from Theorem \ref{theorem:characteristicmap}. This is
possible since the stabilizer of the essential image does not
contain $U$ and so the essential image is not contained in $L$
fixed points. If $H$ is the stabilizer of the essential image of
$\phi$ this implies that $\La/\Ha$ is semisimple and therefore
that $U<H$. But this implies that $U$ is contained in the
stabilizer of the essential image of $\phi'$, a contradiction.
\end{proof}

Before proving Theorem \ref{theorem:characteristicmap} we recall
several facts from \cite{M2}.  The following, which is
\cite[I.4.6.2]{M2}, is a simple corollary of the Poincar\'{e}
recurrence theorem.

\begin{lemma}
\label{lemma:recurrenceestimate} Let $A$ be an automorphism of
$(S,\mu)$ and $f$ a non-negative measurable function on $X$.  Then
for almost all $x{\in}X$
$$\liminf_{m{\rightarrow}\infty}\frac{1}{m}f(A^m(x))=0$$
and
$$\liminf_{m{\rightarrow}-\infty}\frac{-1}{m}f(A^m(x))=0.$$
\end{lemma}

Let $A$ be an ergodic automorphism of $(S,\mu)$ and $W$ be a $k$
vector space. Let $u:{\mathbb Z}{\times}S{\rightarrow}GL(W)$ be a
${\mathbb Z}$-integrable cocycle over the $\mathbb Z$-action
generated by $A$.   If we define
$$\chi_+(u,x,w)=\lim_{m{\rightarrow}{\infty}}\frac{1}{m}\ln\|u(A^m,x)w)\|$$
and
$$\chi_-(u,x,w)=\lim_{m{\rightarrow}{-\infty}}\frac{1}{m}\ln\|u(A^m,x)w)\|$$
it follows from Oseledec multiplicative ergodic theorem that
$\chi_+(u,x,w)$ and $\chi_-(u,xw)$ exist for almost all $x{\in}X$
and all $w{\in}W$.  Furthermore, that theorem shows that there
exists a finite set $J$ and real numbers $\chi_j(u)$ and maps
$\omega_j(u,x):S{\rightarrow}Gr_{l(j)}(W)$ such that
\begin{enumerate}
\item for almost all $x{\in}S$, the space $W$ is the direct sum
$\oplus_J\omega_j(u,x)$; \item for almost all $x{\in}S$ the
sequence $\{\frac{1}{m}\ln\|u(m,x)\|/\|w\|\}_{m{\in}{\mathbb
N}^+}$ converges to $\chi_j(u)$ uniformly in
$w{\in}\omega_j(u,x)-\{0\}$.
\end{enumerate}
Furthermore
$$\{0\}{\cup}\{w{\in}\fh-\{0\}|\chi_+(u,x,w){\leq}a\}={\oplus}_{\chi_j{\leq}a}\omega_j(u,x)$$
and
$$\{0\}{\cup}\{w{\in}\fh-\{0\}|\chi_-(u,x,w){\geq}a\}={\oplus}_{\chi_j{\geq}a}\omega_j(u,x).$$
All of this is contained in \cite[V.2.1]{M2}.  We refer to
$\omega_j(u,x)$ as the {\em characteristic subspace} for $u$ with
{\em characteristic number} $\chi_j$.

We call a cocycle $v:G{\times}X{\rightarrow}H$ {\em
$G$-quasi-integrable} if $v$ is cohomologous to a $G$-integrable
cocycle $u$.  If $v$ is a $\mathbb Z$-quasi-integrable cocycle
then $v(g,x)=\psi(gx)u(g,x)\psi(x){\inv}$ for some $\mathbb
Z$-integrable cocycle $u$.  We then define
$\omega_j(v,x)=\psi(x)\omega_j(u,x)$.  It is easy to verify, using
Lemma \ref{lemma:recurrenceestimate}, that if $v$ is in fact
$\mathbb Z$-integrable, our two definitions of $\omega_j(v,x)$
agree and so $\omega_j(v,x)$ is well-defined and independent of
the choice of $u$.  Though $\omega_j(v,x)$ does not satisfy the
dynamical condition $2$ above, it can be shown to satisfy weaker
dynamical conditions, see Definition 2.4 and Remark 2.4 in
\cite{Z6}.

In the proof below we will need some functorial properties of
characteristic subspaces.  For $\Za$-integrable cocycles these are
\cite[V.2.3]{M2} and it follows easily from the definition that
they hold for $\Za$-quasi-integrable cocycles as well. Let
$u:\Za{\times}S{\rightarrow}GL(W)$ be a $\Za$-quasi-integrable
cocycle and let $Q<W$ be a subspace such that $u(m,x)Q=Q$.  Let
$V$ be the quotient $W/Q$ and $p:W{\rightarrow}V$ the projection.
We have two additional cocycles $u^Q(m,x)=u(m,x)|_Q$ and
$u^V(m,x)=p{\circ}u(m,x)$ both of which can easily be seen to be
quasi-integrable.  Then for any characteristic subspace
$\omega_l(u^V,x)$ (respectively $\omega_l(u^Q,x)$) there is a
characteristic subspace $\omega_{l(j)}(u,x)$ such that
$\omega_l(u^V,x)=p(\omega_{l(j)}(u,x))$ (respectively
$\omega_l(u^Q,x)=\omega_{l(j)}(u,x){\cap}Q$).

For certain cocycles it is easy to compute characteristic
subspaces and numbers.  Let $\sigma:\Za{\rightarrow}GL(W)$ be a
homomorphism, let $M=\sigma(1)$ and let
$c:\Za{\times}S{\rightarrow}GL(W)$ be a cocycle taking values in a
compact subgroup of $GL(W)$.  We let $u(m,x)=\sigma(m)c(m,x)$.  We
let $\Omega(M)$ be the set of all eigenvalues of $M$,
$W_\lambda(M)$ the eigenspace corresponding to
$\lambda{\in}\Omega(M)$
 and $W_d(M)=[\oplus_{\ln|\lambda|=d}W_{\lambda}(M)]_k$.
We also let $W_+(M)=\oplus_{d>0}W_d(M)$ and
$W_{{\leq}0}(M)=\oplus_{d{\leq}0}W_d(M)$. We will call $d$ a {\em
characteristic number} of $M$ and $W_d(M)$ a {\em characteristic
subspace} of $M$.  Then the characteristic numbers of $u$ are the
characteristic numbers of $M$ and the characteristic subspaces for
$u$ are the characteristic numbers for $M$. Furthermore the space
$\oplus_{\chi_j{\leq}0}\omega_j(u,x)=W_{{\leq}0}(M)$ and
$\oplus_{\chi_j{>}0}\omega_j(u,x)=W_{+}(M)$

\begin{proof}[Proof of Theorem \ref{theorem:characteristicmap}]
As the proof is very involved, we divide it into several steps.
The basic idea is to choose an element $t$ of $G$ and use Oseledec
theorem to construct characteristic maps from
$S{\rightarrow}Gr_m(\fl)$ for $\alpha$ and each $g{\inv}tg$. This
gives an $\alpha$-invariant map
$\phi:G{\times}S{\rightarrow}Gr_m(\fl)$, which we show descends to
an $\alpha$-invariant map
$\phi:\BG{\times}S{\rightarrow}Gr_m(\fl)$.  We then pass to
characteristic subspaces for the cocycle $\alpha'$ which is
cohomologous to $\alpha$ and where
$p_F{\circ}{\alpha'}=\pi{\cdot}c'$.  We use the functoriality of
characteristic subspaces and the form of $\alpha'$ to compute the
characteristic subspaces quite explicitly. Finally using the
assumption that $\pi(G)$ does not commute with $U$, we show that
the stabilizer of the essential image does not contain $U$.

{\bf Step One:  Choosing $t$.}

We call a subgroup diagonalizable if it can be conjugated to a
subgroup of the group of diagonal matrices. Recall that a subgroup
$\Sa_i<\Ga_i$ is called a {\em maximal torus} if it is maximal
diagonalizable subgroup of $\Ga_i$.  We fix a torus $\Sa_i$ in
each $\Ga_i$.  We let $\Ta_i<\Sa_i$ be the {\em maximal split
torus}, i.e. the maximal subgroup of $\Sa_i$ that is
diagonalizable over $k_i$. We let $X(\Ta_i)$ be the group of $k_i$
characters of $\Ta_i$, and
 $T_i=\Ta_i(k_i)$

\begin{lemma}
\label{lemma:splitelement} There exists an element
$t{\in}\prod_{i{\in}I}T_i$ such that
\begin{enumerate}
\item the group generated by $t$ is not contained in any proper
normal subgroup of $G$

\item for any $\chi{\in}X(\Ta_i)$ where $\chi(\tau_i(t))$ has
modulus one it follows that $\chi(\tau_i(t))=1$.

\end{enumerate}
\end{lemma}

\begin{proof}
To satisfy $1$, it suffices to choose $t$ such that it projects to
an element which generates an infinite subgroup in each simple
factor of each $\Ga_i$.

For $k_i=\mathbb R$ it suffices to assume that $\chi(\tau_i(t))$
is positive for every $\chi{\in}X(\Ta_i)$.   If $k_i$ is
non-Archimedean, we identify $T_i$ with $(k^*_i)^{l(i)}$ where
$l(i)$ is the $k_i$-rank of $\Ga_i$. We choose  $\pi$ a
uniformizer of $k_i$. We assume that the projection of $\tau_i(t)$
to each copy of $k^*_{i}$ in $(k^*_{\alpha_0})^{i(\alpha_0)}$ is
the product of a unit of $k^*_{i}$ with a non-zero power of $\pi$.
\end{proof}

\noindent {\bf Remark:} Let $\pi$ be  finite dimensional
representation of $G$ on a vector space $V$.  It follows from our
choice of $t$ that if $\pi(t)$ has all eigenvalues of modulus one
then $\pi$ is trivial.

{\bf Step two: Oseledec theorem and characteristic maps.}

 We will
construct the map $\phi:G{\times}S{\rightarrow}Gr_{k}(\fh)$ by
applying Oseledec theorem to certain cocycles over the action of
$g{\inv}tg$ on $S$.

Since $G$ acts ergodically on $S$, it follows from the Mautner
phenomenon that $t$ and therefore $g{\inv}tg$ does as well
\cite[II.3.3]{M2}. We define a map
$\phi':G{\times}S{\rightarrow}Gr_k(\fl)$. The element $g{\inv}tg$
generates a $\Za$ action on $S$.  We define a cocycle
$u_g:\Za{\times}S{\rightarrow}GL(\fl)$ over this $\Za$ action by
$u_g(m,x)=\Ad_{\fl}\alpha(g{\inv}t^mg,x)$ and apply Oseledec
theorem to each cocycle $u_g$. Since different choices of
$g{\in}G$ define cohomologous cocycles over conjugate actions, it
follows that the characteristic numbers $\chi_j(u_g)$ do not
depend on $g$ nor do the dimensions $l(j)$ of the subspaces
$\omega_j(u_g,x)$. We $\chi_j=\chi_j(u_g)$ and
$\omega_j(g,x)=\omega(u_g,x)$. We now have maps
$\omega_j:G{\times}S{\rightarrow}Gr_{l(j)}(\fl)$. If we let $G$ on
$G{\times}X$ by $h(g,x)=(gh{\inv},hx)$, the map
$\omega_j:G{\times}X{\rightarrow}H$ is $\alpha$-invariant. To show
this one uses the cocycle identity to see that
$$\alpha(hg{\inv}s^mgh{\inv},
hx)=\alpha(h,g{\inv}s^mgx)\alpha(g{\inv}s^mg,x)\alpha(h,x){\inv}$$
and notes that
$$\liminf_{m{\rightarrow}{\infty}}\frac{1}{m}\ln^+\|\alpha(h,g{\inv}s^mgx)\|=0$$
for almost every $x$ by Lemma \ref{lemma:recurrenceestimate}.

{\bf Step 3:  $P$-invariance of characteristic maps.}

We now show that there is a minimal parabolic $P$ such that the
map
$\omega_{{\leq}\chi}(g,x)={\oplus}_{\chi_j{\leq}\chi}\omega_i(g,x)$
descends to a map from $\BG{\times}S{\rightarrow}Gr_{k}(\fh)$.
This follows as in \cite[Theorem V.3.3]{M2}.  If we let $P$ be the
set of elements in $G$ such that the set
$M=\{s^mps^{-m}|m{\in}{\mathbb N}^+\}$ is relatively compact in
$G$ then $P$ is a minimal parabolic in $G$ as discussed in
\cite[VI.4.8 and Lemma II.3.1(b)]{M2}. One then can compute that

$$\alpha((pg){\inv}s^mpg,x)=\alpha((pg){\inv}s^mps^{-m}gg{\inv}s^mg,x)$$
$$=\alpha((pg){\inv},s^mpgx)\alpha(s^mps^{-m},s^mgx)\alpha(g,g{\inv}s^mgx)\alpha(g{\inv}s^mg,x).$$

\noindent Because $\alpha$ is $G$-integrable and $M$ is
precompact, it then follows that
$Q_{M,\alpha}(s^mx)=\supp_{m{\in}{\mathbb
N^+}}\ln^+\|\alpha(s^mps^{-m},y)\|<{\infty}$ for almost all
$y{\in}X$ and therefore that
$$\liminf_{m{\rightarrow}{\infty}}\frac{1}{m}\ln^+\|\alpha(s^mps^{-m},s^mgx)\|=0$$
almost everywhere by Lemma \ref{lemma:recurrenceestimate}.
Similarly $\ln^+\|\alpha((pg){\inv}, y)\|{<}{\infty}$ and
$\ln^+\|\alpha(g,y)\|<{\infty}$ for almost all $y{\in}X$ and
therefore
$$\liminf_{m{\rightarrow}{\infty}}\frac{1}{m}\ln^+\|\alpha((pg){\inv},s^mpgx)\|=0$$
and
$$\liminf_{m{\rightarrow}{\infty}}\frac{1}{m}\ln^+\|\alpha(g,g{\inv}s^mgx)\|=0$$
almost everywhere by Lemma \ref{lemma:recurrenceestimate}. It
follows by direct computation that
$\chi^+(pg,x,w){\leq}\chi^+(g,x,w)$ and the reverse inequality
follows by replacing $p$ by $p{\inv}$. See \cite[V.3.3]{M2} for
more detailed computations for certain choices of cocycle.

Let
$\phi'={\oplus}_{\chi_j{\leq}0}\omega_i(g,x):\BG{\times}S{\rightarrow}Gr_k(\fh)$.

{\bf Step 4: Modifying the map.}

Let $\alpha'$ be the cocycle cohomologous to $\alpha$ which
projects via $p_F$ to $\pi{\cdot}c$, and let $\psi$ be the map
implementing the cohomology between $\alpha$ and $\alpha'$. We
will let $\phi([g],x)=\psi(x)\phi'([g],x)$.  Since $\alpha'$ is
$G$-quasi-integrable, it follows from the definitions that
$\phi([g],x)$ is  $\oplus_{\chi_j{\leq}0}\omega'_j([g],x)$ where
$\omega'_j([g],x)=\psi(x)\omega_j([g],x)$ are the characteristic
subspaces for the cocycle $u'(m,x)=\alpha'(g{\inv}t^mg,x)$ over
the $\Za$ action on $S$ generated by $g{\inv}tg$.

{\bf Step 5: Application of functoriality.}

Since $L$ is an algebraic group with unipotent radical $U$ and
Levi complement $F$, we have that $\fl=\ff{\oplus}\fu$. Now $\ff$
and $\fu$ are invariant under $\Ad_{\fl}$ and so are invariant
under the cocycle $\alpha':G{\times}S{\rightarrow}GL(\fl)$.   By
the functoriality of characteristic subspaces discussed above, we
have that $\phi([g],x){\cap}\ff$ is
$\oplus_{\chi_j{\leq}0}\omega^{\ff}_j([g],x)$ where
$\omega^{\ff}_j([g],x)$ are the characteristic subspaces for the
cocycle ${u'}^{\ff}_g(m,x)={\alpha'}(g{\inv}t^mg,x)|_{\ff}$. Since
$\Ad_{\fl}|_{\ff}$ factors through the map $p_F:L{\rightarrow}F$,
and $p_F{\circ}\alpha'(g,x)=\pi(g)c(g,x)$ where $c$ is a cocycle
taking values in a compact group, it follows that
$\phi([g],x){\cap}{\ff}=\Ad_{\fh}(\pi(g){\inv})W_{{\leq}0}(\Ad_{\ff}{\circ}\pi(t))$.

The intersection of $\phi([g],x)$ with $\fu$ is more complicated
to describe. To do this, we let $\fu_0=\fu$ and
$\fu_i=[\fu,\fu_{i-1}]$. Since $\fu$ is unipotent, there is a
number $k$ such that $\fu_k$ is the center of $\fu$ and $\fu_l=0$
for all $l>k$.  Furthermore, $\fu_{i+1}$ is an ideal in $\fu_i$
and we have a sequence of quotients $\fu_i/\fu_{i+1}$.  A key fact
for what follows is that $\Ad(u)(\fu_i){\subset}\fu_{i+1}$ for any
$u{\in}U$. Since $\fu_i$ is $\Ad_{\fl}$ invariant, it follows that
the cocycle $\Ad_{\fl}{\circ}{\alpha'}$ leaves each $\fu_i$
invariant.  Let $p_i:\fu_i{\rightarrow}\fu_{i+1}$ be the
projection.  It follows that
\begin{equation}
\label{equation:unipotent}p_i((\Ad_{\fl}{\circ}{\alpha'}(g,x)|_{\fu_i})v)=p_i(\Ad_{\fl}|_{\fu_i}(\pi(g)c(g,x))v)
\end{equation}
for any $v{\in}\fu_i$.

Since the representation of $G$ on $\fh$ is defined by
$\Ad_{\fh}{\circ}\pi$, and each $\fu_i$ is $\Ad_{\fh}$ invariant,
the $G$ action leaves each $\fu_i$ invariant . Therefore,since $G$
is semisimple, there are $G$ invariant subspaces $\fv_i$ in
$\fu_i$ such that $\fu_i=\fv_i{\oplus}\fu_{i+1}$.  Since $c$ takes
values in a compact subgroup $C<F$, we can assume that $\fv_i$ is
also $c$ invariant. We identify $\fv_i$ and $\fu_i/{\fu_{i+1}}$ as
$G$ modules in the following paragraphs.

Now for any $v{\in}\fv_i$, we can rewrite equation
\ref{equation:unipotent} to
\begin{equation}
\label{equation:actsoncomplement}
p_i(\Ad_{\fl}{\circ}{\alpha'}(g,x)|_{\fu_i})v)=\Ad_{\fl}((\pi(g)c(g,x)))v).
\end{equation}

\noindent It follows from the definition of $\phi([g],x)$ in terms
of characteristic subspaces, the functoriality of characteristic
subspaces and equation \ref{equation:actsoncomplement} above, that
$p_i(\phi([g],x){\cap}\fu_i)=p_i(\Ad_{\fl}(\pi(g){\inv})W_{{\leq}0}((\Ad_{\ff}{\circ}\pi)|_{\fv_i}(t)))$.

{\bf Step 6: $U$ does not stabilize the essential image of
$\phi$.}

Since $\pi(G)$ does not commute with $U$ by assumption, the
representation of $G$ on $\fv_i$ is non-trivial for some $i$. We
fix one such $i$ for what follows.  Since $G$ is semisimple
$\Ad_{\fl}{\circ}\pi(G)|_{\fv_i}<SL(\fv_i)$.  By our choice of $t$
and $i$ this implies that the decomposition
$$\fv_i=\big(\fv_i{\cap}W_+((\Ad_{\fl}{\circ}{\pi})_{\fv_i}(t))\big){\oplus}\big(\fv_i{\cap}
W_{{\leq}0}((\Ad_{\fl}{\circ}{\pi})_{\fv_i}(t))\big)$$
 is non-trivial.
We let
$$V_+= \fv_i{\cap}W_+((\Ad_{\fl}{\circ}{\pi})|_{\fv_i}(t))$$
and
$$V_-=\fv_i{\cap}
W_{{\leq}0}((\Ad_{\fl}{\circ}{\pi})|_{\fv_i}(t)).$$ For any choice
of $g{\in}G$ we denote $(\Ad_{\fl}{\circ}{\pi}(g){\inv})V_+$ by
$V^g_+$ and we denote $(\Ad_{\fl}{\circ}{\pi}(g){\inv})V_-$ by
$V^g_-$. Note that $\fv_i=V_+^g{\oplus}V_-^g$ is non-trivial for
all choices of $g$.  Letting $\bar V^g_+$ (resp. $\bar V^g_-$) be
the projection of $p_i(V^g_+)$ (resp. $p_i( V^g_-)$) to and ${\bar
\phi}([g],x)=p_i(\phi([g],x){\cap}\fu_i)$ we have that $\bar
V^g_-={\bar \phi}([g],x)$ so $\bar V^g_+{\cap}{\bar
\phi}([g],x)=0$.

Let $\ft$ be the Lie algebra of $\pi(T)$ and write ${}^g\ft$ for
$\Ad_{\fl}(\pi(g){\inv}(\ft)$. Then
${}^g\ft{\subset}\phi'([g],x){\cap}\ff$ for almost every
$([g],x)$. Now $V^g_+{\subset}[{}^g\ft,V^g_+]$ since
$\Ad_{\fh}(\pi(g){\inv}\pi(t)\pi(g))V^g_+=V^g_+$. We let $U_g$ be
the collection of elements of $U$ of the form $\Exp(V^g_+)$ where
$\Exp:\fu{\rightarrow}U$ is the Lie group exponential map. That
$V^g_+{\subset}[{}^g\ft,V^g_+]$ implies that, for some
$u_0{\in}U_g$ the projection of
$(\Ad_{\fl}(u_0)(\phi([g],x){\cap}\ff)){\cap}\fu_i$ contains
non-zero vectors in $\bar V^g_+$ and
$p_i((\Ad_{\fl}(u_0)(\phi([g],x){\cap}\ff)){\cap}\fu_i)$ is not
contained in ${\bar \phi([g],x)}=\bar V^g_-$. This implies that,
for $g=g_0$ fixed, the subgroup of $U$ generated by $u_0$ does not
stabilize the essential image of $\phi([g_0],x)$. We want the same
conclusion for a set of $g$ of positive measure.  We note that the
spaces
$(\phi([g],x){\cap}\ff)=\Ad_{\fl}(\pi(g{\inv}))W_{{\leq}0}((\Ad_{\fl}{\circ}{\pi})_{\ff}(t))$
and $\bar V^g_-$ depend continuously on $g$, and that the action
of $u_0$ via $\Ad_{\fl}$ is continuous.  So there exist a small
$\epsilon>0$ such that
$$p_i((\Ad_{\fl}(u_0)(\phi([g_0],x){\cap}\ff)){\cap}\fu_i){\nsubseteq}\bar
V^{g_0}_-$$ implies
$$p_i((\Ad_{\fl}(u_0)(\phi([g],x){\cap}\ff)){\cap}\fu_i){\nsubseteq}\bar
V^{g}_-$$ for all $g$ in $B(g_0,\epsilon)$. This immediately
implies that $u_0$ does not stabilize
$\phi([B(g_0,\epsilon)]{\times}S)$ which suffices to see that
$u_0$ is not in the stabilizer of the essential image of $\phi$.
\end{proof}

\subsection{Property T and cocycles into amenable and reductive groups}
\label{subsection:useofT}

In this section we note some results for cocycles for groups with
property T of Kazhdan .

\begin{theorem}
\label{theorem:T} Let $D$ be a group with property T of Kazhdan
and $A$ be an amenable group. Assume $D$ acts ergodically on $S$
preserving $\mu$. Let $\alpha:D{\times}S{\rightarrow}A$ be a
cocycle. Then $\alpha$ is cohomologous to a cocycle taking values
in a compact subgroup of $A$.
\end{theorem}

This is \cite[Theorem 9.1.1]{Z1}.  Our first application of this
result is to cocycles with reductive target.  We will also require
the following algebraic lemma.

\begin{lemma}
\label{lemma:centralizinglifts} Let $\La$ be a reductive group and
$p_Z:\La{\rightarrow}\La/Z(\La^0)$ the natural projection. Let
$\Fa<[\La^0,\La^0]$ be a connected subgroup.  Let $g{\in}\La$, and
assume that $[p_Z(g),p_Z(\Fa)]$ is trivial.  Then $[g,\Fa]$ is
also trivial.
\end{lemma}

\begin{proof}
The commutator $[g,\Fa]$ is contained in $Z(\La^0)$ by assumption.
Since $[\La^0,\La^0]$ is normal in $\La$ it follows that
$[g,\Fa]<[\La^0,\La^0]$. Therefore $[g,\Fa]$ is a connected
subgroup of the finite group $Z(\La^0){\cap}[\La^0,\La^0]$.
\end{proof}

\begin{theorem}
\label{theorem:Gsuperrigidity2} Let $G$ act weakly irreducibly on
$S$ preserving $\mu$. Let $\alpha:G{\times}S{\rightarrow}L$ be a
cocycle with algebraic hull $L$.  Assume in addition that $G$ has
property T of Kazhdan and that $\La$ is reductive. Then $\alpha$
is cohomologous to a cocycle $\beta=\pi{\cdot}c$. Here
${\pi:G{\rightarrow}L}$ is a continuous homomorphism and
$c:G{\times}S{\rightarrow}C$ is a cocycle taking values in a
compact group $C$ centralizing $\pi(G)$.
\end{theorem}

\begin{proof}  Since $\La$ is reductive, the connected component $\La^0$ is
reductive.  This implies that $\La^0=[\La^0,\La^0]Z(\La^0)$ where
$[\La^0,\La^0]$ is semisimple, $Z(\La^0)$ is the center of $\La^0$
and $[\La^0,\La^0]{\cap}Z(\La^0)$ is finite.  Since $Z(\La^0)$ is
characteristic in $\La^0$, it is normal in $\La$ and the quotient
$\Ja=\La/Z(\La^0)$ is semisimple. We let $p_J:\La{\rightarrow}\Ja$
and look at the cocycle $p_J{\circ}{\alpha}$.  By Theorem
\ref{theorem:Gsuperrigidity1}, we have that $p_J{\circ}{\alpha}$
is cohomologous to a $\pi'{\cdot}c'$ where $\pi':G{\rightarrow}J$
is a continuous homomorphism and $c'$ is a cocycle taking values
in a compact subgroup $C'$ of $J$ which commutes with $\pi'(G)$.
By Lemma \ref{lemma:liftingreps}, the Zariski closure $\Ja_E$ of
$\pi'(G)$ is semisimple and connected. Letting $\Ja_K$ be the
Zariski closure of $c'(G{\times}S)$, it is clear that $\Ja$ is the
almost direct product $\Ja_E{\times}\Ja_K$. The map
$[\La,\La]{\rightarrow}\Ja$ is an isogeny, so again by Lemma
\ref{lemma:liftingreps}, we can lift $\pi'$ to a homomorphism
$\pi$ from $G$ to $[\La,\La]<\La$.  Then the Zariski closure of
$\pi(G)$ is a connected semisimple subgroup of $\La$ whose
projection to $\Ja$ is $\Ja_E$.  Let $\La_K<\La$ be the pre-image
of $\Ja_K$.  Then by Lemma \ref{lemma:centralizinglifts}, the
Zariski closure of $\pi(G)$ commutes with $\La_K$.  We now see
that $\alpha$ is cohomologous to $\pi{\circ}{\alpha'}$ where
$\alpha'$ takes values in a subgroup $H$ of $\La_K(k)<L$ and is a
cocycle by Lemma \ref{lemma:ratioiscocycle}.  In fact $H$ is the
pre-image in $\La_K(k)$ of $C'$ and so is a compact extension of
$Z(\La^0)(k)$ and is therefore an amenable group. Since $G$ has
property T of Kazhdan, it follows from \ref{theorem:T} that
$\alpha'$ is cohomologous to a cocycle taking values in a compact
group.
\end{proof}

\noindent {\bf Remark:} In Theorem \ref{theorem:Gsuperrigidity2}
$G$ can be replaced by $\Gamma$. This can be proven in two ways.
It follows from Theorem \ref{theorem:Gsuperrigidity2} exactly as
Theorem \ref{theorem:Gammasuperrigidity1} follows from Theorem
\ref{theorem:Gsuperrigidity1}. Or, using the Borel density theorem
and Theorem \ref{theorem:Gammasuperrigidity1}, one can modify the
proof of Theorem \ref{theorem:Gsuperrigidity2} to prove the same
result for $\Gamma$.

\subsection{Proofs of Theorems \ref{theorem:Gsuperrigidity} and \ref{theorem:Gammasuperrigidity}}
\label{subsection:proofs}

We note here that Theorem \ref{theorem:Gsuperrigidity} holds under
weaker hypotheses.  Namely, it holds for $G$-integrable cocycles
over weakly irreducible actions for $G$ as in Theorem
\ref{theorem:Gsuperrigidity2}. However, the existence of ergodic
decompositions of measures makes the formulation of Theorem
\ref{theorem:Gsuperrigidity} more useful in applications.

\begin{proof}[Proof of Theorem \ref{theorem:Gsuperrigidity}]
We first verify that we can use Theorems
\ref{theorem:algebraichull} and \ref{theorem:Gsuperrigidity2}
under the assumptions of Theorem \ref{theorem:Gsuperrigidity}.

For a proof that $G$ has property T, we refer the reader to
\cite[III.5]{M1} for the case where all factors of $G$ are
algebraic.  When we replace some factors by their topological
universal covers, the resulting $G$ is a central extension of a
$T$ group which does not split with respect to any non-trivial
subgroup of the center.  It follows that $G$ has $T$ by an
argument due to Serre presented in \cite[Proposition.3.d.17]{HV}.
The $G$ action is weakly irreducible since it is ergodic and $G$
contains no rank one simple factors.

Let $L=J{\ltimes}U$ be a Levi decomposition as above, and
$p_J:L{\rightarrow}J$ be the natural projection. Combining Theorem
\ref{theorem:algebraichull} with Theorem
\ref{theorem:Gsuperrigidity2}, we have the following:
\begin{enumerate}

\item a cocycle $\beta$ cohomologous to $\alpha$;

\item a homomorphism $\pi:G{\rightarrow}J<L$ such that $\pi(G)$
commutes with $U$

\item a cocycle $c:G{\times}S{\rightarrow}C$ where $C<Z_J(\pi(G))$
 is compact

\item $p_J{\circ}{\beta}(g,x)=\pi(g)c(g,x)$

\end{enumerate}

 We can define a map $\tilde
\alpha(g,x)=\pi(g){\inv}\beta(g,x)$. That $\tilde \alpha$ is a
cocycle follows by Lemma \ref{lemma:ratioiscocycle} since $U$
commutes with $\pi(G)$. It suffices to show that $\tilde \alpha$
is cohomologous to a cocycle taking values in a compact group. But
$p_J{\circ}\tilde \alpha=c$ so $\tilde \alpha$ takes values in
$K{\ltimes}U$, an amenable group. Since $G$ has property T of
Kazhdan, we are done by \ref{theorem:T}.
\end{proof}

\begin{proof}[Proof of Theorem \ref{theorem:Gammasuperrigidity}]
We prove Theorem \ref{theorem:Gammasuperrigidity} from Theorem
\ref{theorem:Gsuperrigidity}  by inducing cocycles and actions
exactly as in the proof of Theorem
\ref{theorem:Gammasuperrigidity1}. It is easy to see that the
induced action is ergodic. The only additional difficulty is
verifying that the induced cocycle is $G$-integrable. If $\Gamma$
is cocompact in $G$, the argument is straightforward. Since,
$G/{\Gamma}$ is compact, the strict cocycle
$\beta:G{\times}G/{\Gamma}{\rightarrow}{\Gamma}$ can be defined
using a precompact fundamental domain.  This forces $\beta(g,S)$
to be finite for any $g{\in}G$.  In this case it is easy to verify
that $\tilde \alpha$ is $G$-integrable.

For $\Gamma$ non-uniform, we need to prove that the induced
cocycle is $G$-integrable.  Since these arguments take us rather
far afield, we relegate the proof to the next subsection, see
Proposition \ref{proposition:integrability}.
\end{proof}

\subsection{\bf Integrability of Cocycles}
\label{subsection:integrability}

 When $\Gamma$ is non-uniform, we
need to be careful to verify that the induced cocycle is
$G$-integrable before we can apply Theorem
\ref{theorem:Gsuperrigidity} to the induced action and cocycle.
For $G$ algebraic, the solution is the close to the content of
Lewis' note \cite{L}. Since Lewis' note continues to circulate in
draft form and not appear and since our integrability assumption
for the $\Gamma$ cocycle is weaker than his and our groups $G$ and
$\Gamma$ are more general than his, we give an argument here.

\begin{proposition}
\label{proposition:integrability} Let $G$ be as in the
introduction and $\Gamma<G$ a lattice. Assume $\Gamma$ acts on a
standard measure space $(S,\mu)$ preserving $\mu$.  Let
$\alpha:\Gamma{\times}S{\rightarrow}L$ be a $\Gamma$-integrable
Borel cocycle.  Then there is a choice of fundamental domain
$X{\subset}G$ for $\Gamma$ such that the induced cocycle ${\tilde
\alpha}:G{\times}(G{\times}S)/{\Gamma}{\rightarrow}L$ is
$G$-integrable.
\end{proposition}

As noted above, the proposition is only non-trivial when $\Gamma$
is non-uniform.  For now we restrict our attention to the case
where $G=\prod_i{\Ga}_i(k_i)$ and consider the case where a factor
of $G$ is the topological universal cover of $\Ga_i(\mathbb R)$
only at the end of the proof.  We will work with a fundamental
domain $X$ that is a contained in a finite union of (generalized)
Siegel domains.

Recall that the choice of a fundamental domain $X$ allows us to
write any element of $G$ as $\omega(g)\tau(g)$ where
$\omega(g){\in}X$ and $\tau(g){\in}{\Gamma}$.  We fix an embedding
of each $\Ga_i(k_i)$ in $GL(n,k_i)$ and use this to define a norm.
We then take the supremum norm over factors to define a norm on
$G$.  In the case where $\Ga_i(\mathbb R)$ is replaced by it's
topological universal cover, we define a norm as in \cite{F},
section 7.2.

\begin{proposition}
\label{proposition:fundamentaldomain} There exists a fundamental
domain $X$ for $\Gamma$ in $G$
such that:
\begin{enumerate}
\item $\int_X\ln^+\|x\|d\mu_G(x)<{\infty}$;

\item for any compact set $M{\subset}G$, there is a constant $C_M$
 such that $\sup_{g{\in}M}\|w(gx)\|{\leq}C_M\|x\|$ for all
 $x{\in}X$.
\end{enumerate}
\end{proposition}

\noindent
{\bf Remark:}  Proposition \ref{proposition:fundamentaldomain} is
true with no assumption on the rank of $\Ga_{\alpha}$.  For
uniform lattices, the proposition is trivial.  It can easily be
reduced to the case where $\Gamma$ is non-uniform and irreducible.
For irreducible $\Gamma$ there are two cases, one where the rank
of $G$ is one and the other where the rank of $G$ is at least two.
Though we do not use the first case, the fundamental domain
constructed by Garland and Ragunathan for such $\Gamma$ and $G$
can easily be seen to satisfy the proposition \cite{GR}.  In the
second case, by the second authors arithmeticity theorems,
$\Gamma$ is arithmetic in $G$.  We give a proof for the case
$G=\Ga(\mathbb R)$ and $\Gamma$ a non-uniform irreducible
arithmetic lattice, in which case $\Gamma$ is commensurable to
$\Ga(\mathbb Z)$.  The case of non-uniform irreducible arithmetic
$\Gamma$ in more general $G$ is analogous though the notation
becomes more involved.  We note that for irreducible $\Gamma$,
Proposition \ref{proposition:integrability} is true as long as
$G{\neq}SL_2(\mathbb R)$.  This requires different estimates from
the ones given below.

\noindent {\bf Remark:}  The careful reader may have noted that
Proposition \ref{proposition:fundamentaldomain} follows from
\cite[VIII.1.2]{M2}.  Though there may exist a fundamental domain
for which that proposition is correct, the proof indicated there,
using a fundamental domain $X$ contained in a finite union of
Siegel sets is not.  More precisely, for a fundamental domain $X$
of this type, part $a$ of the conclusions there is true as stated,
but part $b$ is true only if $X$ is contained in a single Siegel
domain.  All applications of \cite[VIII.1.2.b]{M2} both in that
text and in the articles \cite{F,L} can be replaced by the
estimate in part $2$ of Proposition
\ref{proposition:fundamentaldomain}.

\begin{proof}[Proof of proposition
\ref{proposition:fundamentaldomain}]
 Recall that we are proving
the proposition in the case where $G=\Ga(\mathbb R)$ and $\Gamma$
commensurable to $\Ga(\mathbb Z)$. Let $T$ be a maximal $\mathbb
Q$ split torus in $G$.
 Let $F$ be a root system for $G$ with respect to $A$, and let ${\Delta}$
be a set of simple roots. Let $A={\{}a{\in}A|{\alpha}(a){\leq}1
{\forall} {\alpha}{\in}{\Delta}{\}}$ be a Weyl chamber. By
standard reduction theory, we can assume there is a finite set
$F{\in}\Ga(\mathbb Q)$  and a bounded set $L{\subset}G$ such that
$X{\subset}LAF$. The first conclusion of the proposition is
standard, and proofs can be found in \cite{B1} or \cite{PlRa}. The
proof of the second conclusion depends on standard facts from
reduction theory.  Our principle reference for these facts is
\cite{B1}, particularly Section 14.4. Though the discussion there
is restricted to real groups, analogous statements are known for
more general $G$ as above.

 Let $l=\dim(A)$.   As in \cite[Section 14]{B1}, we can
find a finite collection $\rho_1,{\cdots},\rho_l$ of $\mathbb Q$
representations from $\Ga$ into $GL(W_i)$ and vectors
$w_i{\in}W_i(\mathbb Z)$ such that
$$\Pa=\{g{\in}\Ga| \rho_i(g)w_i=\chi_i(g)w_i\}$$
where $\chi_i:\Pa{\rightarrow}{\mathbb R}_{{\geq}0}^*$ restricted
to the split torus $\Aa$ is the highest weight of $\rho_i$ and
$\chi_i|_{\Aa}=d_i\alpha_i$ for some simple root $\alpha_i$ and an
integer $d>0$. This implies that for any other weight
$\chi{\neq}\chi_i$ of $\rho_i$ we have
$|\chi(s)|>|\alpha_i(s)|^{-d}|\chi_i(s)|$ where $d>0$ is an
integer.

Given two real valued functions $f$ and $g$ on $x$, we write
$f\prec{g}$ if $f(x){\leq}Cg(x)$ and $f\asymp{g}$ if $f\prec{g}$
and $g\prec{f}$.

Given $x{\in}X$, we write $x=laf$. We take a compact set
$M{\subset}G$ and $x{\in}LA_{1}F$. We write $gx=glaf=l'a'f'\gamma$
for any $g{\in}M$ where $\gamma=\tau(gx){\inv}$.   It follows that
$l'a'f'=\omega(gx)$.  Since $L$ is compact and $F$ is finite, to
prove the proposition it suffices to show that
$\|a\|\asymp\|a'\|$. Since $\chi_i$ for $1{\leq}i{\leq}l$ form a
basis for $X(T)\bigotimes{\mathbb R}$, it is enough to show that
$|\chi_i(a)|\asymp|\chi_i(a')|$ for all $i$.

We apply $gx$ to $f{\inv}w_i$ for each $\rho_i$.  We obtain that
$\|(glaf)(f{\inv}w_i)\|=\|glaw_i\|{\asymp}|\chi_i(a)|$.  We also
have $\|(glaf)(f{\inv}w_i)\|=\|(l'a'f'{\gamma})(f{\inv}w_i)\|$.
For the latter we have
$\|(l'a'f'{\gamma})(f{\inv}w_i)\|{\asymp}|\chi_i(a')|$ if
$f'{\gamma}f{\inv}w_i$ is proportional to $w_i$ and
$\|(l'a'f'{\gamma})(f{\inv}w_i)\|{\succ}|\alpha_i(a')|^{-d}|\chi_i(a')|$
otherwise.  Therefore $|\chi_i(a')|\prec|\chi_i(a)|$.  Replacing
$f{\inv}w_i$ with $(f'{\gamma}){\inv}w_i$ and arguing in the same
manner yields $|\chi_i(a)|\prec|\chi_i(a')|$.  Therefore we have
$|\chi_i(a)|\asymp|\chi_i(a')|$ which suffices to prove the
proposition.  (We note that the constants implicit in the signs
$\asymp$ and $\prec$ used here depend on the compact set $M$.)
\end{proof}

\begin{corollary}
\label{corollary:cocyclebound} There exists a fundamental domain
$X$ for $\Gamma$ in $G$ such that for any compact set
$M{\subset}G$ there is a constant $C_M$ such that
$\ln^+\|\beta(g,x)\|{\leq}C_1(M)\ln^+\|x\|+C_2(M)$.
\end{corollary}

\begin{proof}
We write $gx=\omega(gx)\tau(gx)$ and recall that
$\beta(g,x)=\tau(gx){\inv}$ is the strict cocycle
$\beta:G{\times}G/{\Gamma}{\rightarrow}{\Gamma}$ corresponding to
$X$. This implies that $\beta(g,x)=x{\inv}g{\inv}\omega(gx)$.
Therefore
$$\ln^+\|\beta(g,x)\|{\leq}\ln^+\|x{\inv}\|+\ln^+\|g{\inv}\|+\ln^+\|\omega(gx)\|.$$
Since $G$ is algebraic there exist positive constants $c$ and $d$
such that $\|x{\inv}\|<c\|x\|^d$.  Combined with Proposition
\ref{proposition:fundamentaldomain} this implies that the previous
equation can be rewritten
$$\ln^+\|\beta(g,x)\|{\leq}dC_M\ln^+\|x\|+\ln^+(c)+\ln^+\|g\|.$$
Letting $C_1(M)=dC_M$ and
$C_2(M)=\supp_{g{\in}M}\ln^+\|g\|+\ln^+(c)$, we are done.
\end{proof}

\begin{proof}[Proof of Proposition
\ref{proposition:integrability}] We have assumed that the cocycle
$\alpha$ is $\Gamma$-integrable. If we let $K$ be a finite
generating set for $\Gamma$, and $K^j$ the set of all words in $S$
of length less than $j$, then $\Gamma$-integrability implies that
$\int_S\ln^+\|\alpha(\gamma,x)\|{\leq}Cj$ almost everywhere for
some constant $C$ (not depending on $j$) and all $\gamma{\in}K^j$.
Recall that the cocycle ${\tilde \alpha}$ is defined by ${\tilde
\alpha}(g,[g_0,s])=\alpha(\beta(g,[g_0]),s)$. To show that
${\tilde \alpha}$ is quasi-integrable, it therefore suffice to
show that the word length of $\beta(g,[g_0])$ is an $L^1$ function
on $G/\Gamma$.  We choose a fundamental domain to define $\beta$
as in Proposition \ref{proposition:fundamentaldomain}. For such a
domain it follows that $\ln^+\|x\|$ is in $L^1(X)$.  In what
follows, we identify $G/\Gamma$ with the fundamental domain $X$
and consider $\beta:G{\times}X{\rightarrow}\Gamma$ written as
$\beta(g,x)$.  To finish the argument, one then uses a theorem of
Lubotzky, Mozes and Ragunathan.  Define a distance function on $G$
by choosing a right $G$ and left $K$ invariant Riemannian metric
on $G$.   Then Lubotzky, Mozes, and Ragunathan show that the word
length metric on $\Gamma$ is bilipschitz equivalent to the induced
metric as a subset of $G$ \cite{LMR}. This result, combined with a
simple computation \cite[Proof of Proposition 7.9]{F}, shows that
the word length of $\beta(g,x)$ is bounded by a multiple of
$\ln^+\|\beta(g,x)\|$ plus a constant. One then applies Corollary
\ref{corollary:cocyclebound} to see that
$\ln^+\|\beta(g,x)\|<C_1(M)\ln^+\|x\|+C_2(M)$ for any $g{\in}M$
where $M{\subset}G$ is pre-compact.

Letting $M{\subset}G$ be any pre-compact set, writing
$\|\gamma\|_{\Gamma}$ for the word length of $\gamma$ and
collecting inequalities, we have:
$$\int_{S{\times}X}Q_{M,{\tilde \alpha}}(s,x)=\int_{S{\times}X}\supp_{g{\in}M}{\tilde
\alpha}(g,s,x)$$
$${\leq}C\int_X{\supp_{g{\in}M}}\|\beta(g,x)\|_{\Gamma}$$
$${\leq}C'\int_X{\supp_{g{\in}M}}\ln^+\|\beta(g,x)\|+B$$
$${\leq}{C_1(M)}C'\int_X\ln^+\|x\|+B+C_2(M).$$
This shows that ${\tilde \alpha}$ is $G$-integrable whenever $G$
is an algebraic group.

When $G$ or a simple factor of $G$ is not algebraic, we need to
extend the reduction theory arguments to fundamental domains for
$G/\Gamma$ in this setting. This is done at the end of section
7.2. of \cite{F}.  A key step in the argument there is showing
that $\Gamma{\cap}Z(G)<Z(G)$ is a subgroup of finite index. This
allows one to choose  a fundamental domain for $\Gamma$ in $G$
which is contained in a finite union of connected components of
pre-images of Siegel sets.  In this context there is a choice of
the norm on $G$, since $G$ is not a linear group. For the choice
made in \cite{F}, the extension of the results in \cite{LMR} is
obvious. We note that some statements in \cite{F} inherit the
inaccuracy of \cite[VIII.1.2]{M2}.  All of these inaccuracies can
be fixed easily using Proposition
\ref{proposition:fundamentaldomain} above.
\end{proof}

\subsection{\bf Uniqueness of the superrigidity homomorphism}
\label{subsection:uniqueness}

In this section, we show that the homomorphism $\pi$ appearing in
the formulation of Theorems \ref{theorem:Gsuperrigidity} and
\ref{theorem:Gammasuperrigidity} is unique up to conjugacy.  In
fact, we prove a more general fact that requires no assumption on
the rank of $G$.  In this subsection, $G$ will be as in subsection
\ref{subsection:superrigiditywarmup}, but with no assumption on
the rank of $G$, but still assuming $G$ has no compact factors.
As usual, $\Gamma<G$ is a lattice.

\begin{theorem}
\label{theorem:uniqueness} Let $D=G$ or $\Gamma$. Assume $D$ acts
on $S$ preserving $\mu$. For $j=1,2$, let $\pi_j:G{\rightarrow}L$
be continuous homomorphism, let $Z_k=Z_L(\pi_j(G))$ and let
$c_j:D{\times}S{\rightarrow}C_j$ be a cocycle over the $D$ action
taking values in a compact subgroup $C_j<Z_j$.  Let
$\alpha_j:D{\times}S{\rightarrow}L$ be the cocycle over the $D$
action defined by $\alpha_j(d,x)=\pi_j(d)c_j(d,x)$.  Then if
$\alpha_1$ is cohomologous to $\alpha_2$, the homomorphism $\pi_1$
and $\pi_2$ are conjugate.
\end{theorem}

Before proving the theorem, we recall some terminology and
notation. For any element $g$ of $GL_n(\mathbb R)$, there is a
unique decomposition of $g=us=su$ where $u$ is unipotent and $s$
is semisimple. Further, we have a unique decomposition $s=cp=pc$
where all eigenvalues of $p$ are positive and all eigenvalues of
$c$ have modulus one. We refer to $p$ as the {\em polar part} of
$g$ and denote it by $\pol(g)$.
  For any
subset $\Omega{\subset}GL_n(\mathbb R)$
 we define
$$P(\Omega)={\{}\pol(h):h{\in}\Omega{\}}.$$
In general $P(\Omega)$ is not a subset of $\Omega$, but if
$\Omega$ is a semisimple subgroup without compact factors, then
the Zariski closure of $P(\Omega)$ is $\Omega$.

For non-Archimedean fields $k$, the situation is more complicated.
We fix a uniformizer $\pi$ for the field $k$ and define a {\em
polar element} of $GL_n(k)$ to be an element $p$ all of whose
eigenvalues are powers of $\pi$. We call an element $c$ in
$GL_n(k)$ {\em compact} if $c$ generates a bounded subgroup in
$GL_n(k)$. Each element of $GL_n(k)$ can be written uniquely as
$su$ where $u$ is unipotent and $s$ is semisimple.  We call an
element $g$ {\em quasi-polar} if $g=su$ as above and $s$ can be
written as $s=cp=pc$ where $p$ is polar and $c$ is compact. Once
$\pi$ is fixed, this decomposition is unique.  We note that if $s$
is semisimple then $s^{n!}$ is quasi-polar. As above for $h$
quasi-polar, we denote the polar part by $\pol(h)$. For a subset
$\Omega$ in $GL_n(k)$, we define
$$P(\Omega)=\{\pol(h):h{\in}\Omega \hskip 2pt \text{and} \hskip 2pt h \hskip 3pt \text{ quasi-polar}\}.$$
As before, in general $P(\Omega)$ is not a subset of $\Omega$, but
if $\Omega$ is a semisimple subgroup without compact factors, then
the Zariski closure of $P(\Omega)$ is $\Omega$.

Recall from subsection \ref{subsection:characteristicmaps} that if
$M$ is a linear transformation and we let $\Omega(M)$ be the set
of all eigenvalues of $M$, we call the numbers $d=|\lambda|$ for
$\lambda{\in}\Omega(M)$ characteristic numbers of $M$. If
$W_\lambda(M)$ is the eigenspace corresponding to
$\lambda{\in}\Omega(M)$, we let
 $W_d(M)=[\oplus_{\ln|\lambda|=d}W_{\lambda}(M)]_k$ be the characteristic subspace of $M$ with
 characteristic number $d$.

{\noindent}{\bf Remark:} The key fact about polar elements is that
a polar element is completely determined by it's characteristic
numbers and subspaces.  We note that under any rational
homomorphism $\pi:\Ga(k){\rightarrow}GL_n(k)$ the image of a polar
element is a polar element.  In fact, for $g$ quasi-polar,
$\pol(\pi(g))=\pi(\pol(g))$.  This implies that we can define
polar and quasi-polar elements of a linear algebraic group $\Ga$
and that the definition is independent of the realization $\Ga$.
However, the set of polar elements of $\Ga$ does depend on the
choice of uniformizer for $k$.

For each $i{\in}I$ we fix an almost faithful representation of
$\Ga_i(k_i)$ in $GL_n(k)$.  We will call an element $g$ of $G$
polar if $\tau_i(g)$ is polar whenever it is non-trivial.  We call
a subgroup $F<G$ {\em Zariski dense} if $\tau_i(F)$ is Zariski
dense in $\Ga_i$ for each $i$.

\begin{lemma}
\label{lemma:determinedbypolar} There exists a finite collection
of quasi-polar element $g_1,\ldots,g_l{\in}{\Gamma}$ such that the
group $F$ generated by $pol(g_1),\ldots,pol(g_l)$ is Zariski dense
in $G$.
\end{lemma}

\begin{proof}  This is similar to \cite[Lemma 4.5]{MQ}.
Let $Z$ be the Zariski closure of
$<\pol(\gamma)|\gamma{\in}\Gamma>$.
 The proof follows
from the fact $Z$ is invariant under conjugation by elements of
$\Gamma$ and so also by elements of $G$ by the Borel-Wang density
theorem \cite[Theorem II.4.4]{M2}.  Hence $G{\cap}Z$ is a normal
subgroup of $G$.  That it is all of $G$ follows from results of
Mostow and Prasad-Ragunathan \cite{M,PR} which show that there is
a maximal split torus $A<G$ such that $A{\cap}{\Gamma}$ is a
lattice in $A$.  This implies that $A{\cap}{\Gamma}$ projects to a
non-compact subgroup of semisimple elements of each simple factor
of $G$.

Since  $<\pol(\gamma)|\gamma{\in}\Gamma>$ is Zariski dense in $G$
and algebraic groups satisfy an ascending chain condition, it
follows that there is a finite collection
$\gamma_1,{\ldots},\gamma_l$ such that
$<\pol(\gamma_1),{\ldots},\pol(\gamma_l)>$ is Zariski dense in
$G$.

Though the results in \cite{M,PR} are only stated for the case of
real algebraic $G$, the interested reader may generalize the proof
of \cite{PR} to the more general $G$ considered in the statement
of our theorems.
\end{proof}

\begin{proof}[Proof of Theorem \ref{theorem:uniqueness}]
If $D=\Gamma$, we apply Lemma \ref{lemma:determinedbypolar} to
obtain the group $F$ and $g_1,{\ldots},g_l$. If $D=G$, there
exists a Zariski dense finitely generated subgroup $F$ generated
by polar elements $g_1,{\ldots},g_l$.  In either case, we fix $F$
and $g_1,{\ldots},g_l$ for the remainder of the proof.

We define associated actions of $D$ on $S{\times}{\fl}$  by the
formula  $\rho_j(d)(x,v)=(dx,\Ad_{\fl}{\circ}\alpha_j(d,x))$.
Oseledec' multiplicative ergodic theorem implies that there are
characteristic exponents and subspaces for any $\Za$ action
defined by powers of an element $d$ in $D$.   Let
$\psi:S{\rightarrow}L$ be the measurable function such that
$\alpha_1(d,x)=\psi(dx){\inv}\alpha_2(d,x)\psi(x)$.  It follows
easily from Lemma \ref{lemma:recurrenceestimate} that for any
$d{\in}D$ the characteristic numbers for $\rho_2(d)$ and
$\rho_1(d)$ are equal. In fact, that lemma shows that if $\lambda$
is a characteristic number for $\rho_1(d)$ with characteristic
subspace $\W{1}{d}{\lambda}(x)$ then
$(\psi(x){\inv})\W{1}{d}{\lambda}(x)$ is a characteristic subspace
for $\rho_2(d)$ with characteristic number $\lambda$.  Since each
$\alpha_j$ is the product of a constant cocycle and a compact
valued cocycle, it follows from the discussion immediately
proceeding the proof of Theorem \ref{theorem:characteristicmap}
that the characteristic numbers and subspaces for $\rho_i(d)$ are
just the characteristic numbers and subspaces for the linear
representation $\Ad_{\fl}{\circ}\pi_j$.  This implies that
$\W{j}{d}{\lambda}(x)=W_{\lambda}(\pi_j(d))$. Combining these two
facts, we see that
$\psi(x)W_{\lambda}(\pi_1(d))=W_{\lambda}(\pi_2(d))$ for every
$d$, every $\lambda$, and almost every $x$.

We let $\{\la{j}{m}{k}\}$ be the characteristic numbers of
$\pi_j(g_m)$ for each $1{\leq}m{\leq}l$.  Then by the argument
above, we see that $\la{1}{m}{k}=\la{2}{m}{k}$ and also that
$\psi(x)W_{\la{1}{m}{k}}(\pi_1(g_m))=W_{\la{2}{m}{k}}(\pi_2(g_m))$
for almost every $x$ and $1{\leq}m{\leq}l$.  Fixing $x$ for which
the equation holds for all $m$, and letting $l=\psi(x)$, the
definition of the polar part of an element implies that
$\pol(\pi_1(g_m))=l{\inv}\pol(\pi_2(g_m)l$.  This implies that
$(\pi_1(\pol(g_m)))=l{\inv}(\pi_2(\pol(g_m))l$.  Since the group
$F$ generated by $\pol(g_1),{\ldots},\pol(g_l)$ is Zariski dense
in $G$ and since by Lemma \ref{lemma:liftingreps} each $\pi_j$
factors through a rational homomorphism of some $\Ga_i$ this
implies that $\pi_1=l{\inv}{\pi_2}l$.
\end{proof}

{\noindent}{\bf Remark:} If one of the cocycles $\alpha_j$ is
simply a homomorphism, it is possible to give a simpler proof
based on the Borel density theorem.

\subsection{\bf Cocycles with prescribed projections}
\label{subsection:funnycocycles}

In this subsection we state and prove some variants on the cocycle
superrigidity theorems.  Though these variants hold in many
settings, we only state variants of Theorems
\ref{theorem:Gsuperrigidity} and \ref{theorem:Gammasuperrigidity}.
Therefore, throughout this subsection $G$ will be as in the
introduction, i.e. with the assumption that $G$ has no rank $1$
simple factors.  The variants stated below are needed for our
applications to local rigidity of affine action and are used in
the proof of Theorem \ref{theorem:localcocyclerigidity}, which in
turn is used to prove Theorem \ref{theorem:semiconjugacy}.

Throughout this subsection $\Aa$ and $\Ha$ will be algebraic
$k$-subgroups of $\La$ such that $\La=\Aa{\ltimes}\Ha$. We further
assume that $\Aa$ is a connected semisimple $k$-group. We will
denote the $k$ points $\Ha(k)=H$ and $\Aa(k)=A$. We fix
homomorphisms $\pi_A^G:G{\rightarrow}A$ and
$\pi_A^{\Gamma}:\Gamma{\rightarrow}A$ with Zariski dense image. We
let $p_A:L{\rightarrow}A$ denote the standard projection.

\begin{theorem}
\label{theorem:Gsuperrigidityvar}  Assume $G$ acts ergodically on
$S$ preserving $\mu$. Let $\alpha:G{\times}S{\rightarrow}L$ be a
$G$-integrable Borel cocycle such that
$p_A{\circ}{\alpha}=\pi_A^G$. Further assume that $L$ is the
algebraic hull of the cocycle.  Then there is a measurable map
$\phi:S{\rightarrow}H$ such that
$\beta=\phi(gx){\inv}\alpha(g,x)\phi(x)=\pi(g)c(g,x)$ where
$\pi:G{\rightarrow}L$ is a continuous homomorphism and
$c:G{\times}S{\rightarrow}C$ is a cocycle taking values in a
compact subgroup $C<Z_L(\pi(\Gamma))$.  The fact that
$\phi(S)\subset{H}$ implies that $p_A{\circ}\beta=p_A^{\Gamma}$.
\end{theorem}

\begin{theorem}
\label{theorem:Gammasuperrigidityvar}  Assume $\Gamma$ acts
ergodically on $S$ preserving $\mu$.  Let
$\alpha:\Gamma{\times}S{\rightarrow}L$ be a $\Gamma$-integrable,
Borel cocycle such that $p_A{\circ}{\alpha}=\pi_A^{\Gamma}$.
Further assume $L$ is the algebraic hull of the cocycle.  Then
there is a measurable map $\phi:S{\rightarrow}H$ such that
$\beta=\phi({\gamma}x){\inv}\alpha(\gamma,x)\phi(x)=\pi(\gamma)c(\gamma,x)$
where $\pi:G{\rightarrow}L$ is a continuous homomorphism and
$c:\Gamma{\times}S{\rightarrow}C$ is a cocycle taking values in a
compact subgroup $C<Z_L(\pi(\Gamma))$.  The fact that
$\phi(S)\subset{H}$ implies that $p_A{\circ}\beta=p_A^{\Gamma}$.
\end{theorem}

These variants are proven from Theorems
\ref{theorem:Gsuperrigidity} and \ref{theorem:Gammasuperrigidity}
using Theorem \ref{theorem:uniqueness}, the following lemma and
some facts about the structure of algebraic groups.

\begin{lemma}
\label{lemma:cocyclecomputation} Let $D,R,F$ and $A$ be groups.
Assume $A{\times}F$ acts on $R$ by a (possibly trivial)
homomorphism into $Aut(R)$. We will write an element
$g{\in}(A{\times}F){\ltimes}R$ as $(g_A,g_F,g_R)$ where
$g_A{\in}A,g_F{\in}F$ and $g_R{\in}R$. Let
$\pi_{A}:D{\rightarrow}A$ be a homomorphism and $p_A$ the
projection of from $(A{\times}F){\ltimes}R$ to $A$. Let $D$ act on
a set $X$. Let
$\alpha:D{\times}X{\rightarrow}(A{\times}F){\ltimes}R$ be a
cocycle over the $D$ action on $X$.  Assume:
\begin{enumerate}
\item  $p_A{\circ}\alpha(d,x)=\pi_A(d)$;

\item Let $C<R$ be a subgroup such that $A$ and $F$ normalize $C$;

\item $\alpha$ is cohomologous to a cocycle $\beta$ taking values
in $(A{\times}F){\ltimes}{C}$.

\end{enumerate}

Assume $\lambda(x)=(\lambda_A(x), \lambda_F(x), \lambda_R(x))$ is
the function such that
$$\lambda(dx){\inv}{\alpha(d,x)}{\lambda(x)}=\beta(d,x)$$
Let $\lambda'(x)=(1_A,\lambda_F(x),\lambda_R(x))$. Then
$$\lambda'(gx){\inv}\alpha(g,x)\lambda'(x)=(\pi(d), \beta_F(d,x), \beta'_R(d,x)$$
where $\beta'_R(d,x)={}^{\lambda_A(dx)}\beta_R(d,x)$ and
$\beta'_R(D{\times}X){\subset}C$.
\end{lemma}

\begin{proof}
This is checked directly by multiplication.  We note that the
multiplication in $(A{\times}F){\ltimes}R$ can be written as
$(g_A,g_F,g_R)(h_A,h_F,h_R)=(g_Ah_A,g_Fh_F,
g_R{}^{(g_A,g_F)}(h_R))$ where ${}^{(g_A,g_F)}(h_R)$ is the image
of $h_R$ under the automorphism of $R$ given by $(g_A,g_F)$. Then
$\lambda'(dx){\inv}\alpha(d,x)\lambda'(x)=\lambda_A(dx)\beta(d,x)\lambda_A(x){\inv}$.
Then
$$(\lambda_A(dx),1_F,1_R)(\beta_A(d,x),\beta_F(d,x),\beta_R(d,x))(\lambda_A(x){\inv},1_F,1_R)$$
$$=(\lambda_A(dx)\beta_A(d,x)\lambda_A(x){\inv},\beta_F(d,x),{}^{\lambda_A(dx)}\beta_R(d,x)).$$
Defining $\beta'(d,x)={}^{\lambda_A(dx)}\beta_R(d,x)$, we then
have that $\beta'(D{\times}X){\subset}C$ since
$\beta_R(D{\times}X){\subset}C$ and $A$ normalizes $C$.
\end{proof}

\begin{proof}[Proof of Theorem \ref{theorem:Gsuperrigidityvar}]
We can write $\Ha=\Fa{\ltimes}\Ua$ where $\Fa$ is semisimple and
$\Ua$ is solvable. We first show that we can change $\Aa$ so that
it commutes with $\Fa$. If not, then since $\Aa$ is connected and
semisimple $\Aa$ acts on $\Fa$ by inner automorphisms and
therefore $\Aa$ is an almost direct product $\Aa_1\Aa_2$ where
$\Aa_1$ commutes with $\Fa$ and $\Aa_2$ is (virtually) a subgroup
of $\Fa$.  Let $\Delta{\inv}(\Aa_2)$ be the antidiagonal embedding
of $\Aa_2$ in $\Aa{\ltimes}{\Fa}$.  We replace $\Aa$ by
$\Aa'=\Aa_1\Delta{\inv}(\Aa_2)$.  A simple computation shows that
$\Fa$ and $\Aa'$ commute in the semidirect product
$\Aa'{\ltimes}\Ha$.  In what follows we replace $\Aa$ by $\Aa'$.

If we apply Theorem \ref{theorem:Gsuperrigidity} to $\alpha$ we
see that $\alpha$ is cohomologous to a cocycle $\beta$ where
$\beta$ is of the form $\pi{\cdot}c$ for a homomorphism
$\pi:G{\rightarrow}L$ and a cocycle $c:G{\times}S{\rightarrow}K$
where $K<L$ is compact.  We let $C=K{\cap}\Ua(k)$.  It follows
from Theorem \ref{theorem:uniqueness} that $p_A{\circ}\beta$ is
cohomologous to a homomorphism conjugate to $\pi_A^G$ and that $A$
is the algebraic hull of $p_A{\circ}{\alpha}$.  Since the
algebraic hull of $\alpha$ contains the Zariski closure of $C$, it
follows from Theorem \ref{theorem:algebraichull} that $\Aa$ and
$C$ commute.

 We now apply Lemma \ref{lemma:cocyclecomputation} with $A=A$,
$F=\Fa(k), R=\Ua(k)$ and $C=C$.
\end{proof}

\begin{proof}[Proof of Theorem
\ref{theorem:Gammasuperrigidityvar}] In the case when
$\pi_A^{\Gamma}$ is the restriction to $\Gamma$ of a continuous
homomorphism of $G$, the proof above applies verbatim.  Though
Theorem \ref{theorem:algebraichull} is only stated for $G$ actions
and cocycles, the analogous statement for $\Gamma$ actions and
cocycles is easily proven by inducing cocycles and actions.

When $k=\mathbb R$ and $\pi_A^{\Gamma}$ does not extend, the
argument is even simpler.  We let $\Ha=\Fa{\ltimes}\Ua$ where
$\Fa$ is reductive and $\Ua$ is unipotent.  As above we modify
$\Aa$ so that $\Aa$ commutes with $\Fa$.  This is possible with
$\Fa$ reductive since $\Aa{\ltimes}\Fa$ is reductive and therefore
$\Aa$ commutes with the torus in $\Aa{\ltimes}\Fa$ which contains
the torus in $\Fa$. We let $F=\Fa(\Ra)$ and $U=\Ua(\Ra)$. It
follows from Theorem \ref{theorem:Gammasuperrigidity} and the fact
that $U$ contains no compact subgroups that $\beta(g,x)$ takes
values in $A{\times}F$.  The theorem now follows from Lemma
\ref{lemma:cocyclecomputation} applied to $A,F,R=U$ and $C=1_U$.

The case of $k$ non-Archimedean and $\pi_A^{\Gamma}$ not
extendable is the most complicated and the only place where we
need the full strength of both Theorem \ref{theorem:uniqueness}
and Lemma \ref{lemma:cocyclecomputation}. Let
$\Ha={\Fa}{\ltimes}{\Ua}$ and $\La=({\Aa{\times}\Fa}){\ltimes}\Ua$
as in the proceeding paragraph. Then $L=(A{\times}F){\ltimes}U$
where $U=\Ua(k)$ and $F=\Fa(k)$. We can write $\Aa$ and $A$ as
almost direct products $\Aa=\Aa_1\Aa_2$ and $A=A_1A_2$ where
$A_1=\Aa_1(k),A_2=\Aa_2(k)$ and $p_{A_1}{\circ}\pi_A$ extends to a
continuous homomorphism of $G$ and $p_{A_2}{\circ}\pi_A$ has
bounded image. By Theorem \ref{theorem:algebraichull} applied as
above, $\Aa_1$ commutes with $\Ua$.

Note that $\Ua(k)$ is the injective limit of it's compact
subgroups. This follows from the same fact for the additive group
$\Qa_p$, where the compact subgroups are subgroups with bounded
denominators. Apply Theorem \ref{theorem:Gammasuperrigidity} to
$\alpha$ to obtain a cocycle $\beta$ cohomologous to $\alpha$,
where $\beta(\gamma,x)=\pi(\gamma)c(\gamma,x)$ where $\pi$ is a
continuous homomorphism of $G$ and $c$ is a cocycle taking values
in a compact subgroup $C<Z_L(\pi(G))$.  Write
$\beta(\gamma,x)=(\beta_A(\gamma,x),\beta_F(\gamma,x),\beta_U(\gamma,x))$
using coordinates as in Lemma \ref{lemma:cocyclecomputation}. Let
$C_1$ be the smallest compact subgroup of $U$ containing
$\beta_U(\Gamma{\times}S)$. Note that $A_1$ commutes with $C_1$.

Let $p_A{\circ}{\beta}=\beta_A$.  This is cohomologous to
$(p_A{\circ}{\pi}){\cdot}(p_A{\circ}c)$.  By Theorem
\ref{theorem:uniqueness}, $p_A{\circ}{\pi}$ is conjugate to
$\pi_A^E$ by an element $a{\in}A$.  Let $\beta'$ be the conjugate
of $\beta$ by $a$. Writing
$\beta'(\gamma,x)=(\beta'_A(\gamma,x),\beta'_F(\gamma,x),\beta'_U(\gamma,x))$
as above, it is clear that $\beta'_F=\beta_F$ and that
$\beta'_U(\Gamma{\times}X){\subset}C_2={}^aC_1$ and that $C_2$ is
a compact subgroup of $U$.  Note that $A_1$ commutes with $C_2$.
It also follows that
$\beta'_A(\gamma,x)=\pi^E_A(\gamma)c_A(\gamma,x)$ where
$c_A(\Gamma{\times}X){\subset}C_A$ where $C_A<A_2$ is a compact
subgroup and $\pi_A^E(\Gamma)$ is contained in $A_1$.  We note
that $C_A$ acts on $U$ by automorphisms, and let $K$ be the set of
all images of $C_2$ under the action of $C_A$.  It is clear that
$K$ is compact and is the union of subgroups of $U$.  Since $U$ is
the projective limit of it's compact subgroups, there is a compact
subgroup $C<U$ with $K{\subset}C$.  We now apply Lemma
\ref{lemma:cocyclecomputation} with $A=A, F=F, R=U$  and  $C=C$.
\end{proof}

\section{\bf Orbits in the space of representations}
\label{section:reporbits}

In this section we prove an independent result that is used in the
proof of our results on local rigidity of constant cocycles.  This
result appears to be known, but we include a proof for
completeness. In this section $D$ will be any finitely presented
group and $H={\mathbb H}(k)$ will be the $k$ points of a algebraic
$k$-group $\Ha$ where $k$ is a local field of characteristic $0$.
We will fix a realization $\Ha<GL(W)$.  We will call a
homomorphism $\rho$ from a group $D$ to $\Ha$ completely reducible
if the representation on $W$ given by $\rho(D)<\Ha<GL(W)$ is
completely reducible.  We let $\Hom(D,\Ha)$ be the space of
homomorphisms of $D$ into $\Ha$ which has a natural structure as
an algebraic subvariety of $\Ha^m$ where $m$ is the number of
generators of $D$.  We note that the structure of $\Hom(D,\Ha)$ as
a variety is independent of the presentation of $D$ and that
$\Hom(D,H)$ is the set of $k$-points of the variety. Note that
$\Ha$ (resp. $H$) acts on the space $\Hom(D,\Ha)$ (resp.
$\Hom(D,H)$) by conjugation.

\begin{theorem}
\label{theorem:closedorbits} Let $D$ and $H$ be as above.  Let
$\pi:D{\rightarrow}H$ be any completely reducible homomorphism.
Then
\begin{enumerate}
\item the $\Ha$ orbit of $\pi$ in $\Hom(D,\Ha)$ is Zariski closed and
\item the $H$ orbit of $\pi$ in $\Hom(D,H)$ is Hausdorff closed.
\end{enumerate}
\end{theorem}

Point $(2)$ of Theorem \ref{theorem:closedorbits} follows from
point $(1)$ and  a result of Bernstein and Zelevinsky \cite{BZ}.
The result of Bernstein and Zelevinsky  says that, given an action
of a $k$ group $G$ on a $k$ variety $V$, the $k$ points of any
orbit are locally closed in the Hausdorff topology. An examination
of the proof shows that for Zariski closed orbits, the orbit is
also Hausdorff closed. For an accessible proof in characteristic
zero see \cite[proof of Theorem 6.1]{AB}. Let $K$ be the algebraic
closure of $k$ and consider ${\mathbb H}<GL(W)$ where $W=K^n$.  We
prove part $1$ of Theorem \ref{theorem:closedorbits} from the
following result.

\begin{theorem}
\label{theorem:closedorbitslinear} Let $D$ be as above and let
$\pi:D{\rightarrow}GL(W)$ be any completely reducible
representation.  Then the $GL(W)$ orbit of $\pi$ in
$\Hom(D,GL(W))$ is Zariski closed.
\end{theorem}

\begin{proof}[Proof of Theorem \ref{theorem:closedorbitslinear}]
We let $K[D]$ be the group ring of $D$.  The representation $\pi$
defines a representation $\tilde \pi$ of $K[D]$.  This representation factors
through a finite dimensional quotient $A=K[D]/\ker(\tilde \pi)$ and since $\pi$ is
completely reducible $A$ is a semisimple algebra.  (See for example
\cite[XVII.6]{La}, particularly Theorem $6.1$.)

If $\tilde \pi'$ is in the closure of the $GL(W)$ orbit of $\pi$
then $\tilde \pi'$ also vanishes on $\ker(\tilde \pi)$ and so
$\tilde \pi'$ is also a representation of $A$. We recall that two
representations of a semisimple algebra $A$ are conjugate if and
only if they have the same character, see for example
\cite[Theorem XVII.3.8]{La}. This implies that the $GL(W)$ orbit
of $\pi$  is closed.
\end{proof}

The proof of Theorem \ref{theorem:closedorbits} is modelled on the
proof that conjugacy classes of semisimple elements in algebraic
groups are closed.

\begin{proof}[Proof of $(1)$ in Theorem \ref{theorem:closedorbits}]
The space $\Hom(D,{\mathbb H})$ is an algebraic variety over $k$.
Assume that the ${\mathbb H}$ orbit of $\pi$ is not closed.  Then
the orbit closure $\Zar({\mathbb H}{\cdot}\pi)$ is the union of
${\mathbb H}{\cdot}\pi$ with a collection of subvarieties of
strictly smaller dimension.

Given any homomorphism $\sigma:D{\rightarrow}G$ where $G$ is any
group, the orbit of $\sigma$ in $\Hom(D,G)$ is naturally
identified with $G/{Z_{G}(\sigma(D))}$.  So given a homomorphism
$\bar \pi{\in}\Zar({\mathbb H}{\cdot}\pi){\backslash}({\mathbb
H}{\cdot}\pi)$, we have that $\dim(Z_{\Ha}(\bar
\pi(D)))>\dim(Z_{\Ha}(\pi(D)))$.  We now show that
$\dim(Z_{\Ha}(\bar \pi(D)))=\dim(Z_{\Ha}(\pi(D)))$  for any
completely reducible homomorphisms $\pi$ and $\bar \pi$ with $\bar
\pi{\in}\Zar(\Ha{\cdot}{\pi})$.

Let $\fh$ be the Lie algebra of $\mathbb H$ and $\Ad_{\fh}$ the
adjoint representation of $\mathbb H$ on $\fh$.  Let $\fz(\pi)$
(resp. $\fz(\bar \pi)$) be the $\Ad_{\fh}{\circ}{\pi}(D)$ (resp.
$\Ad_{\fh}{\circ}{\bar \pi}(D)$) invariant vectors in $\fh$. Since
the characteristic of $K$ is zero, we have that $\fz(\pi)$ is the
Lie algebra of $Z_{\Ha}(\pi(D))$ and $\fz(\bar \pi)$ is the Lie
algebra of $Z_{\Ha}(\bar \pi(D))$ (see \cite{B}, II.7). By
construction $\Ad_{\fh}{\circ}{\bar \pi}$ is in the closure of the
$\Ha$ orbit of $\Ad_{\fh}{\circ}{\pi}$ in $\Hom(D,GL(\fh))$.  By
Theorem \ref{theorem:closedorbitslinear} $\Ad_{\fh}{\circ}{\bar
\pi}$ is conjugate to $\Ad_{\fh}{\circ}{\pi}$ by an element of
$GL(\fh)$ and therefore $\dim(\fz(\pi))=\dim(\fz(\bar \pi))$. This
implies that the $\Ha$ orbit of $\pi$ is closed in the Zariski
topology on $\Hom(D,\Ha)$.
\end{proof}

Theorem \ref{theorem:closedorbits} would suffice to prove Theorems
\ref{theorem:localcocyclerigidityG} and
\ref{theorem:localcocyclerigidityGamma}.  However, to prove
Theorem \ref{theorem:localcocyclerigidity}, we need the following.

\begin{corollary}
\label{corollary:smallergrouporbits} Let ${\mathbb L}={\mathbb
A}{\ltimes}{\mathbb H}$ where all groups are $k$-algebraic.  Let
$D$ be a finitely generated group.  If $H={\mathbb H}(k)$ and
$L={\mathbb L}(k)$, then $H$ orbits of completely reducible
homomorphisms in $\Hom(D,L)$ are Hausdorff closed.
\end{corollary}

\begin{proof}
It follows from part $1$ of Theorem \ref{theorem:closedorbits} that $\La$
 orbits in $\Hom(D,\La)$ are  Zariski closed.  We let $U$ be an $\La$ orbit in
$\Hom(D,\La)$. Then for any $u{\in}U$, then
${\La}u={\cup}_{a{\in}\Aa}a{\Ha}u$. Since $\Ha$ is normal in
$\La$, $a{\Ha}u={\Ha}au$. The $\Ha$ action on $U$ is algebraic, so
the closure of an orbit must consist of the orbit plus sets of
strictly lower dimension. Since all the sets $a{\Ha}u$ have the
same dimension, this means that the closure of each ${\Ha}au$ in
${\La}u$ is ${\Ha}au$. Since ${\La}u$ is Zariski closed, so is
${\Ha}au$ and in particular ${\Ha}u$. That $Hu$ is Hausdorff
closed follows from the (proof of) the result of
Bernstein-Zelevinsky as in the proof of Theorem
\ref{theorem:closedorbits}.
\end{proof}

\section{\bf Proof of Local Rigidity for Cocycles}
\label{section:generalresults}

In this section we prove Theorems
\ref{theorem:localcocyclerigidityG} and
\ref{theorem:localcocyclerigidityGamma}. Actually we prove Theorem
\ref{theorem:localcocyclerigidity} below which implies Theorems
\ref{theorem:localcocyclerigidityG} and
\ref{theorem:localcocyclerigidityGamma}.  As mentioned in the
introduction Theorem \ref{theorem:localcocyclerigidity} is need to
prove Theorem \ref{theorem:semiconjugacy}. We first fix some
notations for the entire section.  In the first subsection we
formulate Theorem \ref{theorem:localcocyclerigidity} and prove the
theorem modulo a result on perturbations of cocycles taking values
in compact groups.  This last result Theorem
\ref{theorem:Tlocalrigidity}, which holds for any compact valued
cocycle over an action of a group with property T, is proven in
the second subsection.

\subsection{\bf Theorem \ref{theorem:localcocyclerigidity} and
proof} \label{subsection:localcocyclerigidity}
 We will let $k$ be a local field of characteristic zero,
$\La$  an algebraic $k$-group and $\Ha,\Aa<\La$  $k$-algebraic
subgroups such that $\La=\Aa{\ltimes}\Ha$. We further assume that
$\Aa$ is a connected semisimple $k$-group.  We let $L,A$ and $H$
denote the $k$ points of $\La,\Aa$ and $\Ha$ respectively. We will
denote by $D$ a group that is either $G$ or $\Gamma$ where $G$ and
$\Gamma$ are as in the introduction. If $\pi_1$ and $\pi_2$ are
two homomorphisms of a group $B$ into a group $C$ whose images
commute, we will denote by $(\pi_1,\pi_2)$ the homomorphism
defined by $(\pi_1,\pi_2)(b)=(\pi_1(b),\pi_2(b))$ for all
$b{\in}B$. We also fix a continuous homomorphism
$\pi_0:D{\rightarrow}L$ such that $\pi_0=(\pi_A,\pi_H)$ where
$\pi_A:D{\rightarrow}A$ and $\pi_H:D{\rightarrow}H$ are continuous
homomorphisms such that $\pi_A(D)$ and $\pi_H(D)$ commute in $L$.
If $D=\Gamma$, recall that $\pi_0$ is superrigid, which means that
$\pi_0=(\pi^E_0,\pi^K_0)$ where $\pi^E_0$ is the restriction to
$\Gamma$ of a continuous homomorphism from $G$ to $L$, and
$\pi^K_0$ is homomorphism from $\Gamma$ to $L$ with bounded image,
and $\pi_0^E(\Gamma)$ and $\pi_0^K(\Gamma)$ commute. Similar
statements hold for $\pi_A$ and $\pi_H$.  If $D=G$, we abuse
notation by writing $\pi^E_0$ for $\pi_0$ and $\pi^K_0$ for the
trivial homomorphism, and similarly for $\pi_H$ and $\pi_A$. We
note that our assumption that $\pi_0=(\pi_A,\pi_H)$ actually
follows from the structure theory of algebraic group when $D$ is
$G$ or the superrigidity theorems when $D$ is $\Gamma$, though in
both cases possibly only after replacing $\Aa$ by an isomorphic
subgroup of $\La$.

\begin{theorem}
\label{theorem:localcocyclerigidity} Let $(S,\mu)$ be a standard
probability measure space,  $\rho$ a measure preserving action of
$D$ on $S$, and let $\alpha_{\pi_0}:D{\times}S{\rightarrow}L$ be
the constant cocycle over the action $\rho$ given by
$\alpha_{\pi_0}(d,x)=\pi_0(d)$. Furthermore, let
$\alpha:D{\times}S{\rightarrow}L$ be a Borel cocycle over the
action $\rho$ such that:
\begin{enumerate}

\item the cocycles $\alpha_{\pi_0}$ and $\alpha$ are $L^{\infty}$
close

\item the projection of $\alpha$ to $A$ is $\pi_A$.

\end{enumerate}
Then there exist measurable maps $\phi:S{\rightarrow}H$ and
$z:D{\times}S{\rightarrow}Z$ where $Z=Z_L(\pi_0^E(D)){\cap}H$ such
that
\begin{enumerate}

\item we have $\alpha(d,x)=\phi(dx){\inv}(\pi_A(d),
\pi_H^E(d))z(d,x)\phi(x)$;

\item $\phi:S{\rightarrow}H$ is small in $L^{\infty}$;

\item the map $(\pi_A^K(d),z(d,x))$ is a cocycle and is
$L^{\infty}$ close to the constant cocycle defined by $\pi^K_0$.

\item the cocycle $(\pi_A^K,z)$ is measurably conjugate to a
cocycle taking values in a compact subgroup $K$ of $Z$ and $K$ is
contained in a small neighborhood of $\pi^K_0(D)$.
\end{enumerate}

\noindent Furthermore if $S$ is a topological space,
$\supp(\mu)=S$ and $\alpha$ and $\rho$ are continuous then both
$\phi$ and $z$ can be chosen to be continuous.
\end{theorem}

\noindent {\bf Remark:}  Theorem
\ref{theorem:localcocyclerigidity} is an immediate consequence of
Theorem \ref{theorem:Tlocalrigidity} if $\pi_0$ has bounded image,
so we assume throughout that $\pi_0$ has unbounded image.

\noindent {\bf Remark:} If $D=G$ or $k$ is Archimedean or
$\pi_A^K$ is trivial, the map $z$ is a cocycle.  More generally,
it is a {\em twisted cocycle}, as discussed in \cite[Section
7.1]{F}.

As in the proof of Theorem 3.1 of \cite{MQ}, we deduce our result
from the stability of partially hyperbolic vector bundle maps,
though the details of the argument are quite different.  To make
the line of the argument clear, we outline it briefly for $G$
actions assuming that the group $A$ above is trivial and that the
$G$ action on $S$ is ergodic. (In other words we sketch a proof of
Theorem \ref{theorem:localcocyclerigidityG} from the introduction,
with the additional assumption that the $G$ action is ergodic.) It
follows from Theorem \ref{theorem:Gsuperrigidity} that
$\alpha(g,x)=\phi(gx){\inv}\pi'(g)c(g,x)\phi(x)$ where
$\pi':G{\rightarrow}L$ is a continuous homomorphism and
$c:G{\times}S{\rightarrow}C$ is a cocycle taking values in a
compact group $C<Z_L(\pi'(G))$ and $\phi:X{\rightarrow}L$ is a
measurable map. Since the set of polar elements of $G$ is Zariski
dense in $G$ (see subsection \ref{subsection:funnycocycles}) we
can find a finitely generated Zariski dense subgroup
$<g_1,{\ldots},g_l>=F<G$ where each generator $g_i$ has the
property that it's image under any rational homomorphism
$\pi:G{\rightarrow}L$ is uniquely determined by the characteristic
numbers and subspaces of $\pi(g_i)$.  We realize $L$ as a subgroup
of $GL(n,k)$ and study the dynamics of the skew product actions
$\rho_{\alpha_{\pi_0}}$ and $\rho_{\alpha}$ on $S{\times}k^n$
determined by $\alpha_{\pi_0}$ and $\alpha$.  It follows from
standard arguments on stability of partially hyperbolic vector
bundle maps, see Lemma \ref{lemma:hyperbolicdynamics} below, that
for any compact set $K<G$, if we pick $\alpha$ close enough to
$\alpha_{\pi_0}$, the characteristic numbers and subspaces for
$\rho_{\alpha_{\pi_0}}(g)$ can be made arbitrarily close to those
of $\rho_{\alpha}(g)$.  This allows us to show that for almost any
$x{\in}S$ and any $i$, the image of ${}^{\phi(x)}\pi'(g_i)$ is
close to the image of $\pi_0(g_i)$.  It then follows from Theorem
\ref{theorem:closedorbits} that $\pi_0$ and $\pi'$ are conjugate
as homomorphisms of $F$, and therefore as homomorphisms of $G$
since $F<G$ is Zariski dense.  Most of the remaining conclusions
of the theorem are deduced by a more careful analysis of the data
coming from Lemma \ref{lemma:hyperbolicdynamics} and Theorem
\ref{theorem:closedorbits}.  That $z$ is measurably conjugate to a
cocycle taking values in a compact group contained in a
neighborhood of the identity uses Theorem
\ref{theorem:Tlocalrigidity}.  The general case follows more or
less the same outline, using Theorems
\ref{theorem:Gsuperrigidityvar} and
\ref{theorem:Gammasuperrigidityvar} in place of Theorem
\ref{theorem:Gsuperrigidity} and requiring somewhat more care due
to the presence of many ergodic components.  For $\Gamma$ actions
and cocycles there is an additional nuance since we cannot choose
$F<\Gamma$ and here we use Lemma \ref{lemma:determinedbypolar}.

If $D=\Gamma$ then by Lemma \ref{lemma:determinedbypolar} of
subsection \ref{subsection:uniqueness} there are elements
$g_1,{\ldots},g_l{\in}D$ and such that the group $F$ generated by
their polar parts is Zariski dense in $G$.  If $D=G$, we pick a
collection $g_1,{\ldots},g_l$ of polar elements in $G$ such that
the group $F$ generated by $g_1,{\ldots},g_l$ is Zariski dense in
$G$.  In either case, we fix $F$ and $g_1,{\ldots},g_l$ for the
remainder of the section.  The reader is referred to subsection
\ref{subsection:uniqueness} for a discussion of polar elements.

\noindent {\bf Remark:} Under any rational homomorphism
$\pi:\Ga(k){\rightarrow}GL_n(k)$ the image of a polar element is a
polar element and $\pol(\pi(g))=\pi(\pol(g))$.

As discussed above, we will use a dynamical argument to show that
$\pol(\tilde \pi_i(g_j))$ is close to $\pol(\pi_0(g_j))$ for a
finite collection of $g_j$ in $\Gamma$ and then use this to
conclude that ${\tilde \pi}_i^E$ and $\pi_0^E$ are close in the
compact open topology on homomorphisms.

We fix an almost faithful representation
$\sigma:L{\rightarrow}GL(n,k)$. We can associate to any action
$\rho$ of $D$ on a space $X$ and any $L$ valued cocycle $\alpha$
over $\rho$, an action $d_{\alpha}{\rho}$ of $D$ on the trivial
bundle $X{\times}k^n$ via $g(x,v)=({\rho}(g)x,
\sigma({\alpha(g,x))v})$. We use this to define two actions
$d_{\alpha_{\pi_0}}\rho$ and $d_{\alpha}{\rho}$, both of which are
linear extensions of $\rho$.

For any linear map $A:k^n{\rightarrow}k^n$, there is a finite
field extension $k'$ of $k$ and a decomposition
$\rn=\bigoplus_{j=1}^{l}W_j$ such that the eigenvalues of
$A|_{W_j{\otimes}k'}$ have the same absolute value $\lambda_j$ for
each $j$, and $\lambda_1>\lambda_2>{\ldots}\lambda_l$. We call
$W_j$ the {\em characteristic subspace} of $A$ with {\em
characteristic number} $\lambda_j$.

Let $S, \mu$ be a finite measure space, $T$ a measure preserving
transformation of $S$, and $\alpha:{\mathbb
Z}{\times}S{\rightarrow}GL(n,k)$ a cocycle.  We let $$W_{\alpha,
\varepsilon, \lambda}(x)=\{w{\in}W|{\hskip
11pt}\limsup_{m{\rightarrow}{\infty}}|\frac{1}{m}\log\|\alpha(m,x)w\|-\lambda|<\varepsilon\}{\cup}\{0\}.$$

\begin{lemma}
\label{lemma:hyperbolicdynamics} Let $A{\in}GL(n,k)$ be linear
transformation,  $(S,\mu)$ a finite measure space, $T$ a measure
preserving transformation of $(S,\mu)$ and $\alpha_0$ be the
constant cocycle over the $T$ action defined by $A$. Then there
exists $\varepsilon_0$ depending only on $A$, such that for any
$\varepsilon<\varepsilon_0$ and any cocycle $\alpha:{\mathbb
Z}{\times}S{\rightarrow}GL(n,k)$ over the $T$ action that is
sufficiently $L^{\infty}$ close (depending on $\varepsilon$) to
$A$ the spaces $W_{A,\lambda_j}$ and $W_{\alpha,\varepsilon,
\lambda_j}(x)$ are $L^{\infty}$ close. Furthermore, the subspaces
$W_{\alpha,\varepsilon, \lambda_j}(x)$ are measurable functions of
$x{\in}S$ and if $T$ and $\alpha$ are continuous, then
$W_{\alpha,\varepsilon,\lambda_i}(x)$ is defined for all $x$ and
depends continuously on $x,T$ and $\alpha$.
\end{lemma}

\noindent {\bf Remark:} The number $\varepsilon_0$ is explicitly
known and can be taken to be one third of
$\min_{1{\leq}j<l}(\lambda_j-\lambda_{j-1})$.

\noindent This lemma follows from standard arguments on the
stability of invariant distributions for partially hyperbolic
vector bundle maps.  There are many possible sources for such
arguments, which go back at least as far as Anosov \cite{An}.  For
a proof that is easily adapted to this setting see the proofs of
Theorem 1 and Lemma 1 in \cite{P}.

Let $g_1,{\ldots},g_k$ be as above.  For $\pi_0$, we let
$\lambda_j^l$ be the characteristic numbers of $\pi_0(g_l)$ and
let $W_j^l$ be the corresponding subspaces.  For $\alpha$ as in
Theorem \ref{theorem:localcocyclerigidity}, we apply Lemma
\ref{lemma:hyperbolicdynamics} to $g_1,{\ldots},g_k$. We denote
the resulting subspaces of $W$ by $W_{\alpha, \varepsilon, g_l,
\lambda_j^l}(x)$.

Before proceeding with the proof of Theorem
\ref{theorem:localcocyclerigidity}, we record the information
given by applying the results of section
\ref{section:superrigidity} to $\alpha$.

\begin{proposition}
\label{proposition:superrigidreduction}  Let $D,S,\rho,
\mu,\alpha$ and $\pi$ be as above.   On each ergodic component
$\mu_i$ of the measure $\mu$ the cocycle
$\alpha:D{\times}S{\rightarrow}L$ is cohomologous to the product
of a constant cocycle with a compact valued cocycle. More
precisely there exist $\mu_i$-measurable maps
$\phi_i:S{\rightarrow}L$ such that
$$\alpha(d,x)=\phi_i(dx)^{-1}(\pi_A,\pi_i)(d)c_i(d,x)\phi_i(x)$$
for all $d{\in}D$ and $\mu_i$ almost every x.  Here $\pi_i$ is the
restriction to $D$ of a continuous homomorphism
$\pi_i:G{\rightarrow}H$ and $c_i:D{\times}S{\rightarrow}C_i$ are
measurable maps taking values in compact subgroups $C_i<L$
commuting with $(\pi^E_A,\pi_i)(D)$.  Furthermore, we can assume
that there is a finite set $\Pi$ of homomorphisms of $G$ in $H$
containing all $\pi_i$ and that each $\phi_i$ takes values in $H$.
\end{proposition}

\begin{proof}[Proof of Proposition \ref{proposition:superrigidreduction}]
The proposition follows by applying Theorems
\ref{theorem:Gsuperrigidityvar} and
\ref{theorem:Gammasuperrigidityvar} to $\alpha$. To apply those
theorems, we need to see that $\alpha$ is $D$-integrable. Since
$\alpha$ is $L^{\infty}$ close to a constant cocycle, it follows
that $\ln^+Q_{M,{\alpha}}(x)$ is essentially bounded for any
precompact $M$ and so in $L^1$.  Therefore the cocycle is
$D$-integrable for almost every ergodic component $\mu_i$ of
$\mu$.

Theorems \ref{theorem:Gsuperrigidityvar} and
\ref{theorem:Gammasuperrigidityvar} show that $\alpha$ is
cohomologous to $\beta$ where $p_A{\circ}\beta=\pi_A$ and
$\beta={\bar \pi_i}c_i$ where $\bar \pi_i:G{\rightarrow}L$ is a
continuous homomorphism and $c_i:D{\times}S{\rightarrow}C$ is a
cocycle taking values in a compact group $C_i<Z_L({\bar
\pi_i}(D))$. Furthermore, we have that $p_A{\circ}\beta=\pi_A$ and
that the cohomology $\phi_i$ takes values in $H$.  To prove the
proposition, we need to see that we can write $\bar
\pi_i=(\pi_A^E,\pi_i)$.

Let $\Ha=\Fa{\ltimes}\Ua$ be Levi decomposition, with $\Fa$
reductive and $\Ua$ unipotent.  Since $\Aa$ is connected, it
follows that $\Fa=\Fa_1\Fa_2$ an almost direct product, where
$\Aa$ and $\Fa_2$ commute as subgroups of $\La$ and $\Aa$ acts on
$\Fa_1$ by automorphisms via a homomorphism
$\Aa{\rightarrow}\Fa_1$ composed with the adjoint action of
$\Fa_1$ on itself.  We write
$\La=((\Aa{\ltimes}\Fa_1)\Fa_2){\ltimes}\Ua$ and as usual denote
$F_1=\Fa_1(k),F_2=\Fa_2(k)$ and $U=\Ua(k)$.  The map
$\Delta(a,f)=(a,a{\inv}f)$ is an isomorphism between
$\Aa{\ltimes}\Fa_1$ and $\Aa{\times}\Fa_1$.  By replacing
$\Aa<\La$ by $\Delta(\Aa)$, we have that $\Aa$ and $\Fa$ commute.
This does not effect $\pi_0$, since we have assumed that $\pi_A$
and $\pi_H$ commute which forces $\pi_H$ to take values in $F_2$.
After conjugating by some element of $U$, we may assume $\bar
\pi_i$ takes values in $\Aa{\times}\Fa$.  Writing
 ${\bar \pi_i}(d)=\pi_A^E(d)\pi_i(d)$ where $\pi_i$ takes values in
$\Fa$, it follows from the fact that $\Fa$ and $\Aa$ commute that
$\pi_i$ is a homomorphism.

That all $\pi_i$ are contained in a finite collection $\Pi$
follows from the fact that there are only finitely many (conjugacy
classes of) homomorphisms of $G$ into $L$.
\end{proof}

\noindent We denote by $\tilde \pi_i$ the homomorphism
$(\pi_A,\pi_i)$ and call $\Pi'$ the collection of all $\tilde
\pi_i$. We fix the collection $\Pi'$ of homomorphisms of $G$ into
$L$ (or equivalently $\pi_A$ and the collection $\Pi$ of
homomorphisms of $G$ into $H$) for the remainder of the section.

Let $\tilde \pi_i(g_l)$
have characteristic numbers $\A{i}{l}1,...,\A{i}{l}{s_i}$ and
$\W{i}{l}{j}$ be the characteristic subspaces corresponding to
$\A{i}{l}{j}$. Then for any $w{\in}W-\{0\}$ and almost every
$x{\in}\supp{\mu_i}$, there is a $j$ such that
$$\limsup_{m{\rightarrow}{\infty}}\frac{1}{m}\log\|\alpha(g_l^m,x)w\|=\A{i}{l}{j}.$$
Since the set $\{\lambda_m^l,\A{i}{l}{j}\}_{i,l,j,m}$ is finite,
if we choose $\alpha$ close enough to $\pi_0$ and $\varepsilon$
small enough, after re-indexing we have

\begin{equation}
\label{equation:star}\phi_i(x)\W{i}{l}{j}=W_{\alpha, \varepsilon,
g_l, \lambda_m^l}(x)
\end{equation}

\noindent for $\mu_i$ almost all $x$ and all $i$. Furthermore, we
have that for each $j$ there is an $m$ such that
$\A{i}{l}{j}=\lambda_m^l$ and that
$\dim(\W{i}{l}{j})=\dim(W_m^l)$. This proof of equation
\ref{equation:star} is essentially contained in \cite{MQ}
discussion preceding Lemma 3.4 or \cite{QZ} proof of Theorem A.

\begin{proof}[Proof of theorem \ref{theorem:localcocyclerigidity}]

First we show that $\tilde \pi_i=(\pi_A,\pi_H^E)$ for all $i$.
Recall that $\pi_0=(\pi_0^E,\pi_0^K)$ and that $\tilde
\pi_i=({\tilde \pi}_i^E,{\tilde \pi}_i^K)$ where $\pi_0^E$ and
$\pi_0^K$ (resp. ${\tilde \pi}_i^E$ and ${\tilde \pi}_i^K$)
commute and $\pi_0^E,{\tilde \pi}_i^E$ are restrictions of
continuous homomorphisms of $G$ and ${\tilde \pi}_i^K,\pi_0^K$
have bounded image. Note that for any $g{\in}\Gamma$ we have
$\pol({\tilde \pi_i}(g))=\pol({\tilde \pi}_i^E(g))$ and that the
same is true for $\pi_0$.  Since $\pi_0=(\pi_A,\pi_H)$ and $\tilde
\pi_i=(\pi_A,\pi_i)$ and $\pi_i$ is the restriction of a rational
homomorphism from $\Ga$ to $\La$, it suffices to show that
${\tilde \pi_i^E}=\pi_0^E$.

By Lemma \ref{lemma:hyperbolicdynamics},by choosing $\varepsilon$
small enough and $\alpha$ close enough to $\pi_0$,  the space
$W_{\alpha,\varepsilon,g_l, \lambda_m^l}(x)$ can be made
arbitrarily close to the space $W_m^l$ for almost every $x{\in}S$.
By equation \ref{equation:star}, we have that
$\phi_i(x)\W{i}{l}{j}=W_{\alpha,\varepsilon,g_l, \lambda_m^l}(x)$
for $\mu_i$ almost every $x{\in}S$, so $W_m^l$ is close to
$\phi(x)\W{i}{j}{l}$ for almost every $x$.  Furthermore,  by the
remark following equation \ref{equation:star}, for each $j$ there
is an $m$ such that the action of the polar part of
${}^{\phi(x)}{\tilde \pi_i}(\gamma_l)$ on ${\phi(x)}\W{i}{l}{j}$
is the same as the action of the polar part of $\pi_0(\gamma_l)$
on $W_m^l$. This implies that, for $\alpha$ close enough to
$\pi_0$, we can make $\pol({}^{\phi_i(x)}{\tilde
\pi_i}(\gamma_l))={}^{\phi_i(x)}{\tilde \pi_i}(\pol(\gamma_l))$
arbitrarily close to $\pol(\pi_0(\gamma_l))=\pi_0(\pol(\gamma_l)$
for every $l$ and $i$ and $\mu_i$ almost every $x$.  Therefore,
for almost every $x{\in}S$, the homomorphisms ${}^{\phi(x)}{\tilde
\pi}_i^E$ and $\pi_0^E$ can be made arbitrarily close as
homomorphisms of $F$ by choosing $\alpha$ close enough to $\pi_0$.
Since there are only finitely many homomorphisms in $\Pi'$ and $H$
orbits in $\Hom(F,L)$ are closed by Corollary
\ref{corollary:smallergrouporbits}, this implies that ${\tilde
\pi}_i^E$ and $\pi_0^E$ are in the same $H$ orbit in $\Hom(F,L)$
and so are conjugate by an element of $H$.  Since $F$ is Zariski
dense in $G$ and ${\tilde \pi}_i^E$ and $\pi_0^E$ are restrictions
of rational homomorphisms, it follows that ${\tilde \pi_i}^E$ is
conjugate to $\pi_0^E$ as homomorphisms of $G$.

By relabelling, we now have

$$\alpha(g,x)={\phi_i}(gx){\inv}(\pi_A,\pi_H^E)(g)c_i(g,x)\phi_i(x)$$

\noindent for each $i$ and $\mu_i$ almost every $x$ where each
$c_i$ is a cocycle taking values in a compact subgroup of
$Z=Z_L((\pi_0^E)(D)){\cap}H$.

We define a map $\bar \phi_h:S{\rightarrow}\Hom(F,L)$ by taking
$x$ to ${}^{\phi_i(x)}{\tilde \pi_i}|_F$.  By Lemma
\ref{lemma:hyperbolicdynamics} and the discussion above, this map
is measurable and $L^{\infty}$ small and has image contained in a
single $H$ orbit in $\Hom(F,L)$.  Since $F$ and $D$ are Zariski
dense in $G$, and $\pi_0^E$ is rational, it follows that
$Z_L(\pi_0^E(F))=Z_L(\pi_0^E(G))=Z_L(\pi_0^E(D))$ and so we may
identify the $H$ orbit in $\Hom(F,L)$ with $H/Z$ where
$Z=Z_L(\pi_0^E(D)){\cap}H$ as before. Therefore, choosing a point
$w$ in the image, we can define a measurable map
$\phi_h:S{\rightarrow}H$ such that $\bar \phi_h(x)=\phi_h(x)w$.
Furthermore, we can choose $\phi_h$ to be $L^{\infty}$ small. This
uses that we can choose a Borel section of $H{\rightarrow}H/Z$
which is continuous in a small neighborhood of $[1_H]$ and does
not increase norms.

It follows from the definitions of $\phi_i$ and $\phi_h$ that for
each $i$, we have that $\phi_i(x)=\phi_h(x)\phi_i^z(x)$ where
$\phi_i^z$ takes values in $Z$. Then a simple computation shows
that:
$$\alpha(g,x)=\phi_h(gx)\phi_i^z(gx)(\pi_A,\pi_H^E)(g)c_i(g,x)\phi_i^z(x){\inv}\phi_h(x){\inv}$$
$$=\phi_h(gx)(\pi_A,\pi_H^E)(g)z(g,x)\phi_h(x){\inv}$$
where
$z(g,x)=\pi_A^K(g){\inv}\phi_i^z(gx)\pi_A^K(g)c_i(g,x)\phi_i^z(x){\inv}$
takes values in $Z$.

Since $\alpha$ is measurable and $L^{\infty}$ close to $\pi_0$ and
$\phi_h$ is measurable and $L^{\infty}$ small, it follows that
$z(g,x)=\phi_h(gx)\alpha(g,x)\phi_h(x){\inv}(\pi_A,\pi_i^E)(g){\inv}$
is measurable and $L^{\infty}$ close to $\pi_H^K$.  It is clear
that $(\pi_A^K(g),z(g,x))$ is a cocycle taking values in
$Z_L(\pi_0^E)$, which implies that $z$ is a cocycle when $\pi_A^K$
is trivial. Theorem \ref{theorem:Tlocalrigidity} in the next
section shows that $(\pi_A^K,z(g,x))$ is measurably conjugate to a
cocycle taking values in a group that is contained in a small open
neighborhood of $(\pi_A^K, \pi_H^K)(D)$.  When $k$ is Archimedean,
this implies that $(\pi_A^K,z(g,x))$ is measurably conjugate into
the closure of $(\pi_A^K, \pi_H^K)(D)$ by Lemma
\ref{lemma:closesubgroupsconjugate}. This implies that $\pi_A^K$
and $z$ take values in groups which commute, which by Lemma
\ref{lemma:ratioiscocycle} implies that $z$ is a cocycle.

When $\rho$ and $\alpha$ are continuous, we deduce continuity of
$\phi_h$ from the continuity of  $W_{\alpha, \varepsilon, g_l,
\lambda_m^l}(x)$ which follows from Lemma
\ref{lemma:hyperbolicdynamics}. Since we know that
$\pi_i^E=\pi_0^E$, we can rewrite equation \ref{equation:star} as
$\phi(x)W^{l}_{m}=W_{\alpha, x, g_l, \lambda_m^l}$ for each $l$.
This implies that the polar part of ${}^{\phi(x)}\pi_0^E(g_l)$
depends continuously on $x$ for each $l$. Combined with Lemma
\ref{lemma:determinedbypolar} this implies that
${}^{\phi(x)}\pi_0^E$ depends continuously on $x$. Therefore
$\phi(x)$ is continuous modulo the stabilizer of $\pi_0^E$ in $H$,
or $\phi_h$ is continuous.  That $z$ is then continuous as well
follows from the formula defining $z$ two paragraphs above.
\end{proof}

\subsection{\bf Local rigidity of compact valued cocycles}
\label{subsection:Tlocalrigidity}

In this subsection $D$ will denote a locally compact group with
property T of Kazhdan, $A$ will be a locally compact group and
$C<A$ will be a compact subgroup.  As before $(S,\mu)$ will be a
standard measure space and we will assume $D$ acts on $S$
preserving $\mu$.  We will say that a group $C'$ is close to $C$
if there is a small neighborhood $U$ of $C$ such that
$C'{\subset}U$. Throughout $\nu$ will denote right Haar measure on
$A$.  We will say that two cocycles
$\alpha,\alpha':D{\times}S{\rightarrow}A$ are close if there is a
compact generating set $K$ for $D$ and a small neighborhood $U$ of
$1$ in $A$ such that $\alpha(d,s)\alpha'(d,s){\inv}$ is in $U$ for
all $d{\in}K$.

\begin{theorem}
\label{theorem:Tlocalrigidity} Let
$\alpha_0:D{\times}S{\rightarrow}C$ be a cocycle over the $D$
action. Any cocycle $\alpha:D{\times}S{\rightarrow}A$ which is
close to $\alpha_0$ is conjugate to a cocycle into a compact group
$C'$ that is close to $C$.
\end{theorem}

\noindent {\bf Remark:} The proof of Theorem
\ref{theorem:Tlocalrigidity} is simpler if one assumes the $D$
action on $S$ is ergodic.

We will need the following proposition in order to find $C'$
satisfying the conclusions of the theorem.

\begin{proposition}
\label{proposition:closesubgroups} Let $C<A$ be a compact
subgroup. Then given any small enough neighborhood $U$ of $C$,
there exists a compact group $C'<A$, contained in $U$, such that
any subgroup $C''<A$ contained in $U$ is conjugate to a subgroup
of $C'$ by a small element of $A$.
\end{proposition}

\begin{lemma}
\label{lemma:closesubgroupsconjugate} Let $C<A$ be a compact
subgroup such that $A/C$ is a manifold. Then any compact subgroup
$C'$ sufficiently close to $C$ is conjugate to a subgroup of $C$
by a small element of $A$.
\end{lemma}

\begin{proof}
This follows from a barycenter argument that is similar to the one
that shows that any two maximal compact subgroups in a semisimple
group are conjugate.  We let $X=A/C$ and take the $C'$ orbit
$\mathcal{O}$ of the coset of the identity $[1_A]$. Since $C'$ is
close to $C$, it follows that $\mathcal{O}$ is contained in a
small neighborhood of $[1_A]$. We can then take the barycenter or
center of gravity for $\mathcal{O}$.  This is defined as the
unique minimum of the function
$d_{\mathcal{O}}(x)=\int_{\mathcal{O}}d(x,y)^2d{\mu}(y)$ where
$\mu$ is the push-forward of Haar measure on $C'$ to
$\mathcal{O}$. That a barycenter exists and is unique can be
proven from convexity of the distance function on a small enough
neighborhood $U_{[1_A]}$. Convexity of the distance function on
this neighborhood can be proven by comparison with the sphere
whose sectional curvature is the maximum of the sectional
curvatures of two planes in $T(A/C)$, using the fact that for
small enough neighborhoods on the sphere, the distance function is
convex, see for example  \cite[Exercises 2.3(1),p.176]{BH}. The
barycenter is then a fixed point for the $C'$ action and is close
to $[1_A]$ since $c'[1_A]$ is close to $[1_A]$ for all
$c'{\in}C'$. This implies that $C'<aCa{\inv}$ for $a{\in}A$ small.
\end{proof}

\begin{proof}[Proof of Proposition
\ref{proposition:closesubgroups}] We let $A^0$ be the connected
component of the identity in $A$, and $p:A{\rightarrow}A/A^0$ be
the projection. Then, since $A/A^0$ is totally disconnected, there
is an open subgroup $\bar C$ containing $p(C)$.  Since $\bar C$ is
open, if $U$ is a sufficiently small neighborhood of $C$, then
$p(U)$ is contained in $\bar C$.  We will find $C'<p{\inv}(\bar
C)$ and so replace $A$ by $A'=p{\inv}(\bar C)$.

Given any open set $U$ containing the identity in $A'$ there is a
compact normal subgroup $N{\subset}U$ such that $A'/N$ has no
small subgroups, and is therefore a manifold.  This is an
extension by Glushkov \cite{Gl} of results due to Gleason,
Montgomery and Zippin, see  \cite{K} for further discussion,
particularly Theorem 18 and the remark following on page 137. We
let $C'=CN$. Then $A/C'=(A/N)/(C/(C{\cap}N)$ and so is a manifold.
It then follows from Lemma \ref{lemma:closesubgroupsconjugate}
that if $U$ is small enough, any subgroup contained
$C''{\subset}U$ is conjugate to a subgroup of $C'=CN$.
\end{proof}

Given a non-negative, integrable function $h$ on $A$ and a unitary
representation $\rho$ of $A$ on a Hilbert space $\fH$, we define
$\rho(h)=\int_Ah(a)\rho(a)d\nu(a)$. If $\int_Ah(a)d\nu(a)=1$, then
$\|\rho(h)\|{\leq}1$ as verified in \cite[III.1.0]{M2}.

We recall that a locally compact group $D$ has property T of
Kazhdan if the trivial representation is isolated in the unitary
dual.  This has the following consequence, which can be seen as an
effective version of the standard statement ``any
$(\varepsilon,K)$-invariant vector is close to a $D$ invariant
vector".

\begin{lemma}
\label{lemma:propertyTis} Let $D$ be a locally compact group with
property T.  Then for any compact generating set $K$ for $A$ there
is a non-negative continuous function $h$ with support contained
in $K^2$ and $\int_Afd{\mu}_A=1$ and a constant $B=B(K,h)$, such
that for any unitary representation $\rho$ of $D$ on a Hilbert
space $\fH$ and for any vector $v{\in}\fH$, we have
\begin{enumerate}

\item $\lim_{n{\rightarrow}{\infty}}\rho(h)^nv=v_F$ exists;

\item $v_F$ is fixed by $D$;

\item $d(v_F,v){\leq}B\supp_{k{\in}K}d(kv,v)$.
\end{enumerate}
\end{lemma}

\begin{proof}
This follows from the definition of property T and from
\cite[III.1.3]{M2}, which shows that $\rho(h)$ is a contraction on
the orthogonal complement of the $D$ fixed vectors in $\fH$.  That
$h$ can be chosen with support in $K^2$ follows from the
construction of $h$ in the proof of \cite[III.1.1]{M2}.
\end{proof}

\begin{proof}[Proof of Theorem \ref{theorem:Tlocalrigidity}]
We let $\fH=L^2(A)$ and write the natural action of $A$ coming
from the right regular representation as $v{\rightarrow}va{\inv}$.
Fix a vector $v_0{\in}\fH$ whose stabilizer is $C$. We define two
representations $\rho_{\alpha_0}$ and $\rho_{\alpha}$ of $D$ on
$L^2(S,\fH,\mu)$ by
$(\rho_{\alpha_0}(d)f)(x)=f(d{\inv}x)\alpha_0(d,x){\inv}$ and
$(\rho_{\alpha}(d)f)(x)=f(d{\inv}x)\alpha(d,x){\inv}$.  Then the
function $f_0:S{\rightarrow}\fH$ defined by $f_0(x)=v_0$ for all
$x$ is $\rho_{\alpha_0}$ invariant.  It is easy to see that
$\supp_{k{\in}K}d(\rho_{\alpha}(k)f_0,f_0)<\varepsilon$ where
$\varepsilon$ only depends on how close $\alpha$ is to $\alpha_0$.

First assume that the action of $D$ on $S$ is ergodic.  By $2$ and
$3$ of Lemma \ref{lemma:propertyTis}, there is a function
$f{\in}L^2(S,\fH)$ such that $f$ is $\rho_{\alpha}(D)$ invariant
and $\|f-f_0\|_2$ is small.  By the proof of \cite[Lemma
9.1.2]{Z1} one sees that the $A$ action on $\fH$ is tame and so
ergodicity of the $D$ action on $S$ implies that $f$ takes values
in a single $A$ orbit $\mathcal{O}$ in $\fH$. This then implies by
\cite[Lemma 5.2.11]{Z1} that $\alpha$ is equivalent to a
representation into the stabilizer $A_v$ of some vector
$v{\in}{\mathcal{O}}$.  It is easy to verify that $A_v$ is
compact. Since $f$ is $L^2$ close to $f_0$, we can choose $v$ to
be close to $v_0$. This immediately implies that the stabilizer of
$A_v$ is Hausdorff close to the stabilizer of $v_0$.

If the action is not ergodic, we cannot conclude that $f$ takes
values in a single orbit.   If one traces through the above
argument, one sees that, on each ergodic component $\mu_i$ of
$\mu$, $f$ takes values in a single $A$ orbit we can view $f$ as a
$\mu$ measurable map $S{\rightarrow}S{\times}H/C_i$ where $C_i<A$
is a compact subgroup depending on $\mu_i$. That we can find a
$\mu$ measurable function conjugating $\alpha$ to a cocycle
$\alpha'$ taking values in $C_i$ for $\mu_i$ almost every $x$
follows from the existence of a Borel section for the map
$S{\times}H{\rightarrow}S{\times}H/C_i$. The existence of such a
section can be deduced from \cite[Theorem A.5]{Z1}. However, since
we only know that $f$ is close to $f_0$ in $L^2(S,\fH,\mu)$, it
only follows from this argument that $C_i$ is close to $C$ on
"most" ergodic components. We want to see that $f$ is actually
close to $f_0$ in $L^2(S,\fH,\mu_i)$ for almost every ergodic
component.  This is deduced from $1$ of Lemma
\ref{lemma:propertyTis}, since
$f=\lim_{n{\rightarrow}{\infty}}\rho(h)^nf$ and this equation
holds in both $L^2(S,\fH,\mu)$ and $L^2(S,\fH,\mu_i)$.  That $f$
is close to $f_0$ in $L^2(S,\fH,\mu_i)$ for almost every ergodic
component then follows from $3$ of Lemma \ref{lemma:propertyTis}.
This implies that each $C_i$ is contained in a small neighborhood
of $C$, and so is conjugate into a subgroup $C'<A$ contained in a
small neighborhood of $C$ by an element $a_i$.  Conjugating by a
map $\phi:S{\rightarrow}A$ such that $\phi(s)=a_i$ for $\mu_i$
almost every $s$, we have that $\alpha$ is conjugate to a cocycle
taking values in $C'$.
\end{proof}

\section{\bf Affine actions, perturbations and cocycles}
\label{section:cocycles}

In this section we prove Corollary \ref{corollary:qz} and Theorem
\ref{theorem:semiconjugacy}.  In order to do so, we need a more
detailed description of the actions in  Definition
\ref{definition:affine} when the group acting is $G$ or $\Gamma$
as above.  The section is divided into two subsections, the first
giving an algebraic description of affine actions, the second
proving Corollary \ref{corollary:qz} and Theorem
\ref{theorem:semiconjugacy}.

\subsection{\bf Description of affine actions}
\label{subsection:affineactions}

Let $H$ be a connected Lie group and $\Lambda<H$ a discrete
cocompact subgroup. We will let $\Aff(H/{\Lambda})$ be the group
of affine transformation of $H/{\Lambda}$.  Any affine
transformation $f$ of $H/{\Lambda}$ has a lift of the form
$h_f{\circ}L_f$ where $h_f{\in}H$ and $L_f{\in}\Aut(H)$ where
$\Aut(H)$ denotes the group of continuous automorphisms of $H$.
Let $N_{\Aut(H)}(\Lambda)$ be the group of elements
$L{\in}\Aut(H)$ such that $L{\cdot}{\Lambda}{\subset}{\Lambda}$.
Then $L_f{\in}N_{\Aut(H)}(\Lambda)$ and we have a map
$\phi:N_{\Aut(H)}(\Lambda){\ltimes}H{\rightarrow}\Aff(H/{\Lambda})$.
There is a map
$\Delta{\inv}:{\Lambda}{\rightarrow}N_{\Aut(H)}(\Lambda){\ltimes}{\Lambda}$
given by
$(\lambda){\rightarrow}{\Ad((\lambda){\inv}),(\lambda))}$. We
denote the image of this map by $\Delta{\inv}(\Lambda)$. (If $H$
has trivial center, then the image is in fact an anti-diagonal
embedding of $\Lambda$.)

\begin{proposition}
\label{proposition:affinegroup} The kernel of the map
$\phi:N_{\Aut(H)}(\Lambda){\ltimes}H{\rightarrow}\Aff(H/{\Lambda})$
is $\Delta{\inv}(\Lambda)$.
\end{proposition}

\begin{proof}
First note that any diffeomorphism $f$ of $H/\Lambda$ gives rise
to an element $f_*$ of $\Out(\pi_1(H/\Lambda))$.  If $f$ is
trivial, then $f_*$ must be trivial as well.    If
$f=\phi((a,h_0))$ for $a{\in}{N_{\Lambda}}(\Aut(H))$ and
$h_0{\in}H$ and $f_*$ is trivial, then $a$ must be an inner
automorphism of $H$ preserving $\Lambda$. This implies that
$a=z\Ad(\lambda)$ for some $\lambda{\in}\Lambda$ and
$z{\in}Z_{\Aut(H)}(\Lambda)$, the centralizer in $\Aut(H)$ of
$\Lambda$. But then
$(a,h_0)[h]=[z{\cdot}({\lambda}h_0h\lambda{\inv})]=[z{\cdot}({\lambda}(h_0)h)]$
for all $h{\in}H$.  So $(a,h_0)[h]=[h]$ for all $h{\in}H$ if and
only if $z{\cdot}(\lambda{h_0}){\inv}z{\cdot}h=h$ for all
$h{\in}H$. Since $z$ centralizes $\lambda$ this is equivalent to
$\lambda(z{\cdot}h_0)(z{\cdot}h)=h$ for all $h{\in}H$.  Picking
$h=1$ this forces $z{\cdot}h_0=\lambda{\inv}$.  This implies that
$z{\cdot}h=h$ for all $h{\in}H$ which implies that $z=1$ and so
$h_0=\lambda{\inv}$.
\end{proof}

 Most of the
difficulty in proving the theorems we need describing affine
actions derive from tori in the reductive component of $H$. To
deal with this difficulty we replace $H$ and $\Lambda$ by groups
$H'$ and $\Lambda'$ such that the respective quotients are
diffeomorphic and the affine groups are the same.  First we note
the a simple fact about covers.

\begin{lemma}
\label{lemma:covers} Let $H$ be a real Lie group and $\Lambda<H$ a
cocompact lattice.  Let $p:H'{\rightarrow}H$ be a covering map and
$\Lambda'=p{\inv}(\Lambda)$.  Then

\begin{enumerate}
\item $H/{\Lambda}$ is diffeomorphic to $H'/{\Lambda'}$

\item $\Aff(H/{\Lambda})<\Aff(H'/\Lambda')$
\end{enumerate}
\end{lemma}

\begin{proof}
The first claim of the lemma is immediate.  To see the second, we
note that any continuous automorphism $A$ of $H$ lifts to a
continuous automorphism $A'$ of $H'$. This uses the fact that the
fundamental group of $H$ is abelian and so any cover is a normal
cover.  If $A{\cdot}{\Lambda}={\Lambda}$ then
$A'{\cdot}{\Lambda'}={\Lambda'}$.  Also, given any element of $h$,
we can choose an element $h'$ in $H'$ projecting to $h$.  It is
easy to verify that $hA$ and $h'A'$ induce the same diffeomorphism
of $H/{\Lambda}=H'/{\Lambda}'$.
\end{proof}

We now show how to replace $H$ by a cover $H'$, though we need to
use an algebraic structure on $H'$ so that the cover is not a
rational map.

\begin{proposition}
\label{proposition:Hreplacement} Given a real algebraic group $H$
and a cocompact lattice $\Lambda$ there is a cover
$p:H'{\rightarrow}H$ and a realization of $H'$ as $\Ha'(\Ra)$ for
a connected $\Ra$ algebraic group $\Ha'$, such that
\begin{enumerate}

\item there is a finite index subgroup $\Aut^A(H')<\Aut(H')$ such
that all elements of $\Aut^A(H')$ are rational automorphisms of
$H'$ and

\item $\Aut^A(H'){\ltimes}H'$ is the real points of a real
algebraic group which we denote by $\Aut^A(\Ha')\ltimes\Ha'$.
\end{enumerate}
\end{proposition}

\begin{proof}
We first define the group $H'$.  By definition $H=\Ha(\Ra)$ where
$\Ha$ is an algebraic $\Ra$-group. We take a Levi decomposition
$\Ha=\La{\ltimes}\Ua$ where $\La$ is reductive and $\Ua$ is
unipotent. We first pass to a finite central extension $\tilde
\Ha$ so as to be able to assume that $\La$ is a direct product of
a torus $\Ta$ and a simply connected semisimple group $\Ja$.  We
let $\sigma:\La{\rightarrow}\Aut(\Ua)$ be the representation
defining the semidirect product.  We let
$\Ta_1=\ker(\sigma){\cap}\Ta$. This is a finite extension of a
connected group $\Ta_1^0$, and $\Ta_1^0<Z(\tilde \Ha)$ and so
$\tilde \Ha=\Ta_1^0{\times}\Ha^*$. The universal cover of
$\Ta_1^0(\Ra)$ is isomorphic to $\Ra^n$ for $n=\dim(\Ta_1^0)$, and
we can realize $\Ra^n$ as the real points of a unipotent algebraic
group which we denote by $\Ua^*$. We replace $\Ha$ by
$\Ha'=\Ha^*{\times}\Ua^*$ and $H$ by $H'=\Ha'(\Ra)$.  There is a
covering map $p_1^*:\Ua^*(\Ra){\rightarrow}\Ta_1^0(\Ra)$ which
defines a covering map $p_1:\Ha'(\Ra){\rightarrow}{\tilde
\Ha}(\Ra)$ which we compose with the covering map $p_2:{\tilde
\Ha}(\Ra){\rightarrow}\Ha(\Ra)$ to define a covering map
$p:\Ha'(\Ra){\rightarrow}\Ha(\Ra)$.  We let
$\Lambda'=p{\inv}(\Lambda)$.

To continue the proof, we will need a Levi decomposition of
$\Ha'=\La'{\ltimes}\Ua'$.  By the discussion above, we can write
$\La'=\Ja'{\times}{\Ta'}$ where $\Ja'$ is semisimple and $\Ta'$ is
a torus.  We have constructed $\Ha'$ such that the homomorphism
$\Ta'{\rightarrow}\Aut(\Ua')$ has finite kernel. (The attentive
reader will note that $\Ja'$ is isomorphic to $\Ja$ above and that
$\Ta$ above is $\Ta_1{\times}{\Ta'}$, but we will not need this in
the discussion that follows.)  As usual,
$L'=\La'(\Ra),U'=\Ua'(\Ra),J'=\Ja'(\Ra)$ and $T'=\Ta'(\Ra)$.

 The group $\Aut^A(H')$ will consist of those
automorphisms of $H'$ which project to inner automorphisms of $L$.
We first show this has finite index in $\Aut(H)$.  The group of
outer automorphisms of $J$ is finite, so it suffices to show that
the group of automorphisms of $T'$ that extend to automorphisms of
$H'$ is finite.  In fact, the group $\Xi$ of automorphisms of $T'$
which extend to $T'{\ltimes}U'$ is finite.  Any such automorphism
must induce a permutation of the finite collection $\Delta$ of
weights defining the representation $\sigma$ of $\Ta$ on $\Ua$, so
by passing to a subgroup $\Xi'<\Xi$ of finite index, we may assume
that $\Xi'$ fixes $\Delta$ pointwise.  Since the kernel of
$\sigma$ is finite, $\Delta$ forms a basis for the group of
characters of $\Ta'$ which vanish on $\ker{\sigma}$.  Therefore
$\Xi'$ acts trivially on a subgroup of finite index in the group
of characters of $\Ta'$ and, since $\Ta'$ is connected, acts
trivially on $\Ta'$.

We can write any element of $\phi{\in}\Aut^A(H')$ as a composition
of three elements.  First we translate by an element $u$ of $U'$
so that $u{\circ}\phi(L)=L$.  Then we conjugate by an element $l$
of $L$ so that $\Ad(l){\circ}{u}{\circ}{\phi}$ is trivial on $L$.
The automorphism $\Ad(l){\circ}{u}{\circ}{\phi}=a$ is clearly an
automorphism of $U$ which commutes with the action of $L$ on $U$.
We write $\phi=alu$.  Viewing $a$ as belonging to
$Z_{\Aut(U)}(L)$, $l$ as an element of the adjoint group of $\bar
L$ of $L$, and $u$ as an element of $U/Z_U(L)$, this decomposition
is unique, and clearly makes $\Aut^A(H')$ the set of $\Ra$ point
of an $\Ra$-variety.  Writing the multiplication on
$\Aut^A(H'){\cong}Z_{\Aut(U)}(L){\times}{\bar
L}{\times}(U/Z_U(L))$ it is clear that all factors commute
pairwise except the last two. The product ${\bar
L}{\ltimes}(U/Z_U(L))$ is clearly a quotient of the adjoint group
 of $J{\ltimes}U$ by the image of $Z_U(L)$ in the adjoint group of
 ${J{\ltimes}U}$, and so the real points of an algebraic group
defined over $\Ra$, and so $\Aut^A(H')$ is the real points of an
algebraic group defined over $\Ra$.

It also follow easily from our description of $\Aut^A(H')$ that
every element of $\Aut^A(H')$ is the restriction of a rational
automorphism of $\Ha'$. This follows from the fact that
$\Aut(\Ua)$ acts rationally on $\Ua$ which follows from the fact
that $\exp$ and $\ln$ are rational diffeomorphisms between $U$ and
$\fu$ by the Baker-Campbell-Hausdorff formula, and that any
automorphism of $\fu$ is linear and therefore rational.

To show that $\Aut^A(H'){\ltimes}H'$ is the real points of an
algebraic group defined over $\Ra$ only requires that we show that
$\Aut^A(H'){\ltimes}H'{\rightarrow}H'$ is the restriction of a
rational map. Using the coordinates on $\Aut^A(H')$ described
above this reduces to showing that the map
$Z_{\Aut(\Ua)}(\La){\ltimes}\Ua{\rightarrow}\Ua$ is rational. This
follows from the fact that $\Aut(\Ua){\ltimes}\Ua{\rightarrow}\Ua$
is rational which follows from the fact that $\Aut(\Ua)$ is
defined as an algebraic subgroup of $GL(\fu)$ and from rationality
of automorphisms of $U$ discussed above.
\end{proof}

\begin{theorem}
\label{theorem:describingactionsG} Let $G$ be as in section
\ref{subsection:superrigiditywarmup}, but with no assumption on
the rank of $G$. Assume $H$ is a connected real algebraic group
and $\Lambda$ a cocompact discrete subgroup. Let $\rho$ be an
affine action of $G$ on $H/{\Lambda}$. Then the action $\rho$ is
given by $\rho(g)[h]=[\pi_0(g)h]$ where $\pi_0:G{\rightarrow}H$ is
a continuous homomorphism.
\end{theorem}

\begin{theorem}
\label{theorem:describingactionsgamma} Let $G$ be as in section
\ref{subsection:superrigiditywarmup} and let $\Gamma<G$ be a
weakly irreducible lattice. Let $H$ and $\Lambda$ be as above. Let
$\rho$ be an affine action of $\Gamma$ on $H/{ \Lambda}$. Then
there is a finite index subgroup $\Gamma'<\Gamma$ such that,
possibly after replacing $H$ by $H'$ as in Proposition
\ref{proposition:Hreplacement}, the $\Gamma'$ action on
$H/{\Lambda}$ is given by
$\rho(\gamma)[h]=[\pi_H(\gamma){\cdot}\pi_A(\gamma)h]$. Here
$\pi_H:\Gamma'{\rightarrow}H'$ and
$\pi_A:\Gamma'{\rightarrow}\Aut(H')$ are homomorphisms whose
images commute as subgroups of $\Aut(H'){\ltimes}H'$. Furthermore,
we can assume that $(\pi_A,\pi_H)(\Gamma')$ is contained in
$\Aut^A(H'){\ltimes}H'$, an algebraic group.
\end{theorem}

{\noindent} {\bf Remark:} Using the results of
\cite{C,GS,Rg1,Rg2,St} in combination with the arguments in
\cite{M2}, one can assume only that $\Gamma$ projects to a dense
subgroup of a rank one simple factors not locally isomorphic to
$F_4^{-20}$ or $Sp(1,n)$.

It is obvious from this description that the action of $G$ or
$\Gamma$ on $H/\Lambda$ lifts to $H$ on a subgroup of finite
index.

\begin{proof}[Proof of Theorem \ref{theorem:describingactionsG}]
The action $\rho$ defines a continuous homomorphism
$$\pi:G{\rightarrow}(N_{\Aut(H)}(\Lambda){\ltimes}H)/(Z(H){\cap}\Lambda).$$
Since the target is a Lie group, any simple factor $\Fa(k)$ of $G$
which is defined over a non-Archimedean field $k$ has trivial
image, since is totally disconnected and topologically almost
simple.  Therefore it suffices to consider the case where $G$ is a
connected Lie group.   We replace $H$ by a group $H'$ and
$\Lambda$ by $\Lambda'$ as in Proposition
\ref{proposition:Hreplacement}. Let $\Aut^A(H')<\Aut(H')$ be the
subgroup of finite index such that every element of $\Aut^A(H')$
is rational. Then since $\pi(G)$ is connected, $\pi(G)$ must be
contained in the image of
$(\Aut^A(H'){\ltimes}H'){\cap}(N_{\Aut(H')}(\Lambda){\ltimes}H')$.

It follows from generalizations of Borel's density theorem that
$\Lambda'$ is Zariski dense in a cocompact normal subgroup of
$H'$, see for example \cite{D} or \cite[Theorem 1.1]{Sh}.
Therefore by Lemma \ref{proposition:Hreplacement} any automorphism
of $H'$ that fixes $\Lambda'$ pointwise factors through an
automorphism of a compact quotient $\bar H$ of $H$.  It is easy to
see that $\Aut(\bar H)$ is a discrete extension of a compact
group, since $\bar H$ is an almost direct product of compact
simple groups and compact torii. Since $\Lambda'$ is discrete
$N_{\Aut(H')}(\Lambda')/Z_{\Aut(H')}(\Lambda')$ is discrete and
therefore $N_{\Aut(H')}(\Lambda')$ is a discrete extension of a
compact group.  As remarked above $\pi(G)$ is contained in the
connected component of $\Aff(H/\Lambda')$. Letting $Z^0$ be the
connected component of $Z_{\Aut(H')}(\Lambda')$, the connected
component of $\Aff(H/{\Lambda})$ is $\phi(Z^0{\ltimes}H')$.  Now
$Z^0{\ltimes}H'{\cap}\ker(\phi)=Z(H'){\cap}\Lambda'$ so
$\phi(Z^0{\ltimes}H'){\cong}(Z^0{\ltimes}H')/(Z(H'){\cap}\Lambda')$.
Since $G$ has no compact factors and $Z^0$ is compact, the map
$\pi:G{\rightarrow}(Z^0{\ltimes}H')/(Z(H'){\cap}\Lambda')$ takes
values in $H/(Z(H){\cap}{\Lambda'})$.

Since $G$ is either simply connected or simply connected as an
algebraic group, we can lift $\pi$ to a homomorphism $\tilde
\pi:G{\rightarrow}H'$.  This also define a homomorphism
$\pi:G{\rightarrow}H$, and it is easy to verify that $\tilde \pi$
and $\pi$ define the same affine action on
$H/\Lambda{\cong}H'/{\Lambda'}$.
\end{proof}

The proof for $\Gamma$ actions is more complicated and requires
the use of the superrigidity theorems.  In addition to a direct
application, we will also use the following consequence of the
superrigidity theorems.

\begin{lemma}
\label{lemma:liftinggammareps} Let $G$ and $\Gamma$ be as above.
Let $\pi:\Gamma{\rightarrow}D$ be any Zariski dense homomorphism
into a real algebraic group $D$.  Let $\tilde D$ be a real
algebraic group and $\tilde D{\rightarrow}D$ an isogeny. Then
there is a finite index subgroup $\Gamma'<{\Gamma}$ and a
homomorphism $\tilde \pi:\Gamma'{\rightarrow}{\tilde D}$ such that
$p{\circ}{\tilde \pi}=\pi$ where $p:{\tilde D}{\rightarrow}D$ is
the natural covering map.
\end{lemma}

\begin{proof}[Proof of Lemma \ref{lemma:liftinggammareps}]
Since by \cite[IX.5.8]{M2}, the image of any homomorphism from
$\Gamma$ into a real algebraic group has semisimple algebraic
closure, it suffices to consider the case where $D$ is semisimple.
Since it also suffices to consider the case where $\Da$ is simply
connected as an algebraic group and simply connected semisimple
algebraic groups are direct products of simple groups by
\cite[I.1.4.10]{M2}, it suffices to consider the case where $\Da$
is simple and simply connected as an algebraic group.

We first assume $G=\prod_I{G_i}$ where each $G_i$ is algebraic and
$G_1=\Ga_1(\Ra)$. We have a homomorphism
$\pi:\Gamma{\rightarrow}{\mathbb D}({\mathbb R})$. It is clear
that $\pi(\Gamma)<{\mathbb D}(k)$ where $k$ is a finite extension
of $\mathbb Q$, and we let $\bar k$ be the algebraic closure of
$k$.  Since $\Gamma$ is Zariski dense in $D$, it follows that $D$
is defined over $k$.  A corollary of the superrigidity theorems,
see \cite[Theorems VII.6.5 and VII.6.6]{M2}, shows that, after
passing to a subgroup of finite index, there is an embedding
$\sigma$ of $\bar k$ in $\mathbb C$, and a $\mathbb C$-rational
map $\eta:\Ga_1{\rightarrow}{}^{\sigma}\Da$ such that
$\pi(\gamma)=\sigma{\inv}(\eta(\gamma))$.  Here ${}^{\sigma}\Da$
is the $\sigma(k)$ algebraic group defined by the image under
$\sigma$ of the equations defining $\Da$.  Since $\Ga_1$ is simply
connected as an algebraic group, it follows that we can lift
$\eta$ to a map to $\tilde \Da$, and then define the lift of $\pi$
by the same equation.

For $G_1$ a topological cover of a real algebraic group $\bar
G_1$, the argument above  gives the same conclusion concerning
$\pi, \eta$ and $\sigma$ where $\eta$ is a continuous homomorphism
of $G_1$ which factors through $\bar G_1$.
\end{proof}

\begin{proof}[Proof of Theorem
\ref{theorem:describingactionsgamma}]
 The action $\rho$ is described by a homomorphism
$\pi:{\Gamma}{\rightarrow}(N_{\Aut(H)}(\Lambda){\ltimes}H)/\Delta{\inv}(\Lambda)$.
We observe that a finite index subgroup in
$(N_{\Aut(H)}(\Lambda){\ltimes}H)/{\Delta}{\inv}(\Lambda)$ maps
into $({\Aut^A(H)}{\ltimes}H)/{\Delta}{\inv}(\Lambda)$ which maps
onto $({\Aut^A(H)}{\ltimes}H)/{\Delta}{\inv}(H)$. Passing to a
subgroup of finite index in $\Gamma$ and composing $\pi$ with this
inclusion and surjection, we get a map $\bar
\pi:{\Gamma}{\rightarrow}({\Aut^A(H)}{\ltimes}H)/\Delta{\inv}(H)$.
Recall that $({\Aut^A(H)}{\ltimes}H)$ is an algebraic group, and
note that $\Delta{\inv}(H)$ is an algebraic subgroup, so the
quotient $({\Aut^A(H)}{\ltimes}H)/\Delta{\inv}(H)$ is an algebraic
group. Applying the superrigidity theorems to $\pi$, we see that
$\bar \pi=(\bar \pi_E,\bar \pi_K)$ where $\bar \pi_E$ is the
restriction of a continuous homomorphism of $G$ and $\bar \pi_K$
has bounded image. By \cite[IX.5.8]{M2}, we know that $\bar
\pi(\Gamma)$ has semisimple Zariski closure $\Ja$. Here
$\Ja=\Ja_1\Ja_2$ where $\Ja_1$ is isotropic over $\Ra$ and $\Ja_2$
is anisotropic over $\Ra$.

We let ${\Aut^A(H)}{\ltimes}H=\La_1(\Ra){\ltimes}\Ua_1(\Ra)$ where
$\La_1$ is reductive algebraic group and $\Ua_1$ is the unipotent
radical.  Similarly, we let
$\Delta{\inv}(H)=\La_2(\Ra){\ltimes}\Ua_2(\Ra)$ and
$({\Aut^A(H)}{\ltimes}H)/\Delta{\inv}(H)=\La_3(\Ra){\ltimes}\Ua_3(\Ra)$
where $\La_2$ and $\La_3$ are reductive and $\Ua_2$ and $\Ua_3$
are unipotent.  Then $\La_3=\La_1/\La_2$ and these are all
reductive, we can find a subgroup $\tilde \La_3$ such that
$\La={\La_2}{\tilde \La_3}$ is an almost direct product and the
map ${\tilde \La_3}{\rightarrow}\La_3$ is an isogeny. By Lemma
\ref{lemma:liftinggammareps}, there is a finite index subgroup
$\Gamma'<\Gamma$ on which we can lift $\bar \pi$ to a
representation into $\tilde \pi$ into $\La_3(\Ra)$ and therefore
into ${\Aut^A(H)}{\ltimes}H$.  It is easy to verify that $\bar
\pi$ and $\pi$ define the same affine action of $\Gamma'$.

That $\pi$ can be chosen to be of the form $(\pi_A,\pi_H)$
requires a supplementary argument.  Let $\Ha=\La{\ltimes}\Ua$ be
Levi decomposition. It follows from the proof of Proposition
\ref{proposition:Hreplacement} that the Levi complement of
$\Aut^A(H)$ a direct product of the adjoint group $\bar \La$ of
$\La$ and a reductive subgroup $\La'$ of $\Aut(\Ua)$ that commutes
with $\bar \La$.  In the description above, one can take
$\La_2=\Delta{\inv}(\La)$ and $\La_1=\La'{\times}({\bar
\La}{\ltimes}{\La})$.  Therefore $\La_3$ is isomorphic to
$\La'{\times}{\bar \La}$, and if one chooses $\tilde \La_3$ to be
$\La'{\times}\La<\La'{\times}({\bar \La}{\ltimes}{\La})$ one has
the desired conclusion.
\end{proof}

\subsection{\bf Applications}
\label{subsection:whatthewholedamnthingwasfor} We now proceed to
prove the applications listed in the introduction.  For the
remainder of this section $G$ is as defined in the second
paragraph of the introduction and $\Gamma$ is a lattice in $G$.
For any affine action of $G$ or $\Gamma$ on $H/{\Lambda}$, or any
associated quasi-affine or generalized affine action, we assume
that $H$ satisfies the conclusions of Proposition
\ref{proposition:Hreplacement} and so we can describe the action
by Theorem \ref{theorem:describingactionsG} or
\ref{theorem:describingactionsgamma}.  We call the homomorphism
defining the action $\pi_0=(\pi_A,\pi_H)$ and let $\Aa$ be the
Zariski closure of $\pi_A(\Gamma')$ in $\Aut^A(\Ha)$.  If we are
concerned with $\Gamma$ actions, $\pi_0$ only defines the action
on a subgroup of finite index $\Gamma'$.  For the proof of Theorem
\ref{theorem:semiconjugacy}, we replace $\Gamma'$ by a subgroup of
finite index to assure that $\Aa$ is connected.

\begin{proof}[Proof of Corollary \ref{corollary:qz}]
Let $D=G$ or $\Gamma$ and $\rho$ be a generalized standard affine
action of $D$ on a manifold $M$.  Since the entropy of an element
is determined by the entropy of it's $k$th power, it suffices to
prove the corollary for a subgroup $D'$ of finite index.

We first prove the corollary for $\rho$ affine and $M=H/{\Lambda}$
and then describe the modifications for the general case.  For an
affine action we have $TM=H/{\Lambda}{\times}{\fh}$,  and by
Theorems \ref{theorem:describingactionsG} and
\ref{theorem:describingactionsgamma} there is a subgroup of finite
index $D'<D$ and a homomorphism
$\pi_0:D{\rightarrow}Aut(H){\ltimes}H$ such that the derivative
cocycle of the $D'$ action is $\tilde
\pi_0=(\Ad_{\Aut(H){\ltimes}H}{\circ}\pi_0)|_\fh$.

Let $\rho'$ be a $C^2$ action $C^1$ close to $\rho$.  By a result
of Seydoux, $\rho'$ preserves a measure that is in the same
measure class as Lebesgue measure \cite{S}.  Since the derivative
cocycle $\alpha_{\rho'}$ of $\rho'$ is $C^0$ close to the constant
cocycle given by $\tilde \pi_0$ (which is the derivative cocycle
for $\rho$), it follows from Theorem
\ref{theorem:localcocyclerigidity} that $\alpha_{\rho'}$ is
cohomologous to $\tilde \pi_0^E{\cdot}c$ where $\tilde \pi_0^E$ is
the extendable part of $\tilde \pi_0$, and $c$ is cocycle over the
$D$ action taking values in a compact group that commutes with
$\tilde \pi_0^E(D')$. The corollary now follows from the fact that
the entropies $h_{\rho}(d)$ and $h_{\rho'}(d)$ can both be
computed in the same manner from the eigenvalues of $\tilde
\pi_0^E(d)$. This follows from Pesin's formula relating entropy to
Lyapunov exponents as observed by Furstenberg in \cite{Fu}, see
also \cite[Chapter 9]{Z1}.

We now pass to the case of a generalized affine action $\rho$ of
$D$ on $K{\backslash}H/{\Lambda}$.  By definition $\rho$ is the
quotient of an affine action $D$ on $H/{\Lambda}$, so on a
subgroup of finite index, $\rho$ is given by a homomorphism
$\pi_0:D{\rightarrow}\Aut(H){\ltimes}H$.  Since $K<H$ commutes
with $D$, and since $K$ is compact and the Zariski closure of
$\pi_0(D)$ in $\Aut(H){\ltimes}H$ is semisimple, there is a
splitting $\fh=\fk{\oplus}\fm$ invariant under
$\Ad_{\Aut(H){\ltimes}H}$ restricted to both $K$ and $\pi_0(D)$.
The tangent bundle to $K{\backslash}H/{\Lambda}$ can be identified
with $K{\backslash}H/{\Lambda}{\times}\fm$ and, on a subgroup
$D'<D$ of finite index, the derivative cocycle is
$(\Ad_{\Aut(H){\ltimes}H}{\circ}\pi_0)|_{\fm}$.   The remainder of
the proof follows as before.
\end{proof}

\noindent We remark that the proof of the corollary only uses part
of Theorems \ref{theorem:localcocyclerigidityG} and
\ref{theorem:localcocyclerigidityGamma}, namely the conclusion
that the derivative cocycle for the perturbed action is
cohomologous to $\pi_0^E{\cdot}c$ where $c$ takes values in a
compact group.


We now prove Theorem \ref{theorem:semiconjugacy} from the
introduction.  Actually, we prove a slightly stronger statement.
Let $H$ be a connected real algebraic group and $\Lambda<H$ a
discrete cocompact subgroup. Let $\rho$ be a standard affine
action of $D=G$ or $\Gamma$ on $H/{\Lambda}{\times}M$.  By
Theorems \ref{theorem:describingactionsG} and
\ref{theorem:describingactionsgamma} above, there is a finite
index subgroup $D'<D$ such that the action of $D'$ on
$H/{\Lambda}{\times}M$ is given by
$\rho(d)([h],m)=(\pi_0(d)h,\iota(d,h)m)$ where
$\pi_0:D{\rightarrow}{\Aut(H)}{\ltimes}H$ is a homomorphism and
$\iota:D{\times}H/{\Lambda}{\rightarrow}\Isom(M)$ is a cocycle. It
is clear from this description that the action of $D'$ lifts to
$H{\times}M$. For $G$ actions, we let $Z=Z_H(\pi_0(G))$. For
$\Gamma$ actions the description of $Z$ is more complicated. Let
$\Gamma'$ be the subgroup of finite index given by Theorem
\ref{theorem:describingactionsgamma}. Recall that
$\pi_0=(\pi_A,\pi_H)$ and let $A$ be the set of real points of the
Zariski closure of $\pi_A(\Gamma')$, and define $L=A{\ltimes}H$.
We view $\pi_0$ as taking values in $L$. We let
$Z=Z_L(\pi_0^E(\Gamma)){\cap}H$.

\begin{theorem}
\label{theorem:semiconjugacydetailed}   Let $\rho$ be a standard
affine action of $D=G$ or $\Gamma$ on $H/{\Lambda}{\times}M$ as
above. Let $D'$ and $Z$ be as above.  Given any action $\rho'$
sufficiently $C^1$ close to $\rho$, there is a cocycle
$z:D{\times}H/{\Lambda}{\times}M{\rightarrow}Z$ and a continuous
map $f:H/{\Lambda}{\times}M{\rightarrow}H/{\Lambda}$, such that
\begin{enumerate}
\item $f$ is $C^0$ close to the natural projection

\item for any $d{\in}D'$ and any
$([h],m){\in}H/{\Lambda}{\times}M$ we have
$f(\rho'(d)([h],m))=(\pi_A(d),\pi_H^E(d)z(d,([h],m)))f([h],m)$.
\end{enumerate}
\end{theorem}

To prove Theorem \ref{theorem:semiconjugacydetailed}, we recall,
and generalize slightly, the construction, from \cite{MQ} of a
cocycle from a perturbation of an affine action.  We will prove
Theorem \ref{theorem:semiconjugacydetailed} by applying Theorem
\ref{theorem:localcocyclerigidity} to this cocycle.

Once again, we let $D$ denote our acting group.  For the
construction of the cocycle, $D$ can be more general than $G$ or
$\Gamma$ above, as long as the $D$ action is as in the conclusions
of Theorems \ref{theorem:describingactionsG} or
\ref{theorem:describingactionsgamma}. First we define the cocycle
for actions by left translations. Let the $D$ action $\rho$ on
$H/{\Lambda}$ be defined via a homomorphism
$\pi_0:D{\rightarrow}H$. Let $\rho'$ be a perturbation of $\rho$.
If $D$ is connected it is clear that the action lifts to $\tilde
H$ and therefore to $H$. If $D$ is discrete, this lifting still
occurs, since the obstacle to lifting is a cohomology class in
$H^2(D, \pi_1(H/{\Lambda}))$ which  does not change under a small
perturbation of the action. (A direct justification without
reference to group cohomology can be found in \cite{MQ} section
2.3.)  Write the lifted actions of $D$ on $H$ by $\tilde \rho$ and
${\tilde \rho}'$ respectively. We can now define a cocycle
$\alpha:D{\times}H{\rightarrow}H$ by
$${\tilde \rho}'(g)x={\alpha}(g,x)x$$
for any $g$ in $D$ and any $x$ in $H$. It is easy to check that
this is a cocycle and that it is right $\Lambda$ invariant, and so
defines a cocycle $\alpha:D{\times}H/{\Lambda}{\rightarrow}H$. See
\cite{MQ} section 2 for more discussion.

Similarly if $D$ acts on $H/{\Lambda}$ affinely we can obtain a
cocycle as well. Note this only occurs if $D$ is discrete. For a
finite index subgroup $D'<D$ the action is described by a
homomorphism $\pi_0$ into a group $A{\ltimes}H=L$ as described
before the statement of Theorem
\ref{theorem:semiconjugacydetailed}. This $D'$ action lifts to $H$
as noted above.  We define the cocycle only after replacing $D$
with $D'$. As before the  $D$ action defined by $\rho'$ also lifts
to $H$. We let $\pi_A$ be the projection of $\pi_0$ on $A$. Recall
that $\pi_A(D)$ must normalize $\Lambda$ and that $\pi_A(D)$ is
Zariski dense in $A$. Here we consider $H/{\Lambda}$ as the space
$\pi_A(D){\ltimes}H/\pi_A(D){\ltimes}{\Lambda}$, and then define
the cocycle by the same equation as before, noting that we get a
cocycle
$\alpha:D{\times}H/{\Lambda}{\rightarrow}\pi(D){\ltimes}H$.
Observe that there is a map
$\alpha':D{\times}H/{\Lambda}{\rightarrow}H$ such that
$\alpha(d,x)=(\pi_A(d), \alpha'(d,x))$. It is possible to define
$\alpha'$ directly by the same equation as $\alpha$, but it is not
a cocycle.  The map $\alpha'$ is called a twisted cocycle in
\cite{F}, further discussion can be found there and in \cite{MQ}
section 2, example 2.2.

Whenever we have the cocycle $\alpha$ defined as taking values in
$\pi_A(D){\ltimes}H$, we extend this to a cocycle into
$L=A{\ltimes}H$ where $A$ is as above.

If $D$ acts on $H/{\Lambda}{\times}M$ we can also define a cocycle
as above. For simplicity we discuss the definition of this cocycle
when the $D$ action on $H/{\Lambda}$ is given by a homomorphism
$\pi_0:D{\rightarrow}H$. Let
$p:H/{\Lambda}{\times}M{\rightarrow}{H/{\Lambda}}$ be the
projection. In this case we define
${\alpha}:D{\times}H/{\Lambda}{\times}X{\rightarrow}H$ via the
formula
$$\alpha(g,x)p(x)=p({\tilde \rho}'(g)x).$$
The verification that this is a cocycle and defined on
$D{\times}H/{\Lambda}{\times}M$ follows exactly as in \cite{FW}.
The case of more general $\pi_0$ follows exactly as above and
again results in a cocycle into $A{\ltimes}H$. Note that to define
this cocycle, we do not need to know that the skew product action
is isometric on $M$.

In all three cases, we can define a cocycle $\alpha_{\pi_0}$
corresponding to the action $\rho$ it is clear from the definition
that ${\alpha_{\pi_0}(g,x)}={\pi_0}(g)$ in all cases.

\begin{proof}[Proof of Theorem \ref{theorem:semiconjugacydetailed}]
Given $\rho'$ a perturbation of $\rho$, we construct a cocycle
$\alpha:D{\times}H/{\Lambda}{\rightarrow}A{\ltimes}H$ as described
above. It is clear from the construction that $\alpha$ is
continuous and close to the constant cocycle defined by $\pi_0$.
Applying the result of Seydoux, we see that $\rho'$ preserves a
measure $\mu$ that is in the same measure class as Lebesgue
measure. Therefore $\mu$ has full support. Applying Theorem
\ref{theorem:localcocyclerigidity} to $\alpha$, we have that there
is a map $\phi_h:H/{\Lambda}{\times}M{\rightarrow}H$ and a cocycle
$z:H/{\Lambda}{\times}M{\rightarrow}H$ such that
\begin{enumerate}
\item
$\alpha(g,x)=\phi_h(gx){\inv}(\pi_A(g),\pi_H^E(g)z(g,x))\phi_h(x)$

\item $\phi_h$ is continuous and $C^0$ small and depends
continuously on $\rho'{\in}\Diff^1(H/{\Lambda}{\times}M)$

\item $z$ is continuous and $C^0$ close to the constant cocycle
defined by $\pi_H^K$
\end{enumerate}

\noindent We define
$f:H/{\Lambda}{\times}M{\rightarrow}H/{\Lambda}$ by
$f(x)=p({\phi}_h(x)x)$ where we write
$p:{H/{\Lambda}}{\times}M{\rightarrow}H/{\Lambda}$ for the natural
projection. It then follows from the definition of $\alpha$ and
our conclusions about $\phi$ that:
$$f(\rho'(g)x)=p(\phi_h(gx)\rho'(g)x)=p(\phi_h(gx)\alpha(g,x)x)$$
$$=p(\phi_h(gx)\phi_h(gx)^{-1}(\pi_A(g),\pi_H^E(g)z(g,x))\phi_h(x)x)$$
$$=(\pi_A(g),\pi^E_H(g)z(g,x))p(\phi_h(x)x)$$
$$=(\pi_A(g),\pi^E_H(g)z(g,x))f(x).$$
\end{proof}

{\noindent}Theorem \ref{theorem:semiconjugacy} is an immediate
consequence of Theorem \ref{theorem:semiconjugacydetailed} and the
definition of $Z$.

\bigskip
\noindent
David Fisher\\
Department of Mathematics and Computer Science\\
Lehman College - CUNY\\
250 Bedford Park Boulevard W\\
Bronx, NY 10468\\
dfisher@lehman.cuny.edu\\
\\
G.A. Margulis\\
Department of Mathematics\\
Yale University\\
P.O. Box 208283\\
New Haven, CT 06520\\
margulis@math.yale.edu\\

\end{document}